\documentclass[11pt]{article}

\usepackage[margin=1in]{geometry}

\usepackage{enumitem}
\usepackage[utf8]{inputenc}
\usepackage{amsmath,amssymb,amsthm,tabu}
\usepackage{ifthen}
\usepackage{color}
\usepackage{calc}

\usepackage{tikz}
\usetikzlibrary{positioning,arrows,decorations.pathreplacing,shapes,decorations.pathmorphing}
 \usetikzlibrary{calc}

\usepackage{multirow}
\usepackage{boxedminipage}
 \usepackage{xifthen}
 \usepackage{tabularx}

\usepackage{newtxtext,newtxmath}

\usepackage{algorithmic}
\usepackage[linesnumbered,vlined,ruled]{algorithm2e}
\usepackage{hyperref}
\newcommand{\footremember}[2]{%
    \footnote{#2}
    \newcounter{#1}
    \setcounter{#1}{\value{footnote}}%
}

\newlength{\bibitemsep}
\newlength{\bibparskip}
\let\oldthebibliography\thebibliography
\renewcommand\thebibliography[1]{%
	\oldthebibliography{#1}%
	\setlength{\parskip}{\bibitemsep}%
	\setlength{\itemsep}{\bibparskip}%
}

\newcommand{\ALG}{\textsc{Alg}}
\newcommand{\State}{\mathbf{S}}
\newcommand{\state}{\mathbf{s}}
\newcommand{\OPT}{\textsc{Opt}}

\newcommand{\R}{\mathbf{R}}

\newcommand{\Pol}{\mathbf{P}}
\newcommand{\NP}{\mathbf{NP}}

\newcommand{\x}{\mathbf{x}}

\newcommand{\X}{\mathbf{X}}

\DeclareMathOperator*{\E}{\mathbf{E}}
\DeclareMathOperator{\risk}{\mathrm{Risk}}
\DeclareMathOperator{\Prob}{\mathbf{P}}

\DeclareMathOperator{\cost}{cost}

\DeclareMathOperator{\argmin}{argmin}

\newcommand{\halmos}{\qed}

\newtheorem{claim}{Claim}
\newtheorem{example}{Example}[section]

\newtheorem{theorem}{Theorem}[section]
\newtheorem{proposition}[theorem]{Proposition}
\newtheorem{corollary}[theorem]{Corollary}
\newtheorem{lemma}[theorem]{Lemma}

\newtheorem{innercustomgeneric}{\customgenericname}
\providecommand{\customgenericname}{}
\newcommand{\newcustomtheorem}[2]{%
	\newenvironment{#1}[1]
	{%
		\renewcommand\customgenericname{#2}%
		\renewcommand\theinnercustomgeneric{##1}%
		\innercustomgeneric
	}
	{\endinnercustomgeneric}
}

\newcustomtheorem{customthm}{Theorem}
\newcustomtheorem{customlemma}{Lemma}
\newcustomtheorem{customproposition}{Proposition}
\newcustomtheorem{customclaim}{Claim}

\let\proof\relax
\let\endproof\relax
\newenvironment{proof}[1]{\noindent {\bf #1} \ignorespaces}{}

\title{Adaptive Bin Packing with Overflow}

\author{
	Sebastian Perez-Salazar\footremember{2}{Georgia Institute of Technology, \texttt{sperez@gatech.edu}}
	\and Mohit Singh\footremember{3}{Georgia Institute of Technology, \texttt{mohit.singh@isye.gatech.edu}}
	\and Alejandro Toriello\footremember{4}{Georgia Institute of Technology, \texttt{atoriello@isye.gatech.edu}}
}

\allowdisplaybreaks

\begin{document}

\maketitle

\begin{abstract}
Motivated by bursty bandwidth allocation \cite{kleinberg2000allocating} and by the allocation of virtual machines to servers in the cloud \cite{gupta2015lagrangian}, we consider the online problem of packing items with random sizes into unit-capacity bins. Items arrive sequentially, but upon arrival an item's actual size is unknown; only its probabilistic information is available to the decision maker. Without knowing this size, the decision maker must irrevocably pack the item into an available bin or place it in a new bin. Once packed in a bin, the decision maker observes the item's actual size, and overflowing the bin is a possibility. An overflow incurs a large penalty cost and the corresponding bin is unusable for the rest of the process. In practical terms, this overflow models delayed services, failure of servers, and/or loss of end-user goodwill. The objective is to minimize the total expected cost given by the sum of the number of opened bins and the overflow penalty cost. We present an online algorithm with expected cost at most a constant factor times the cost incurred by the optimal packing policy when item sizes are drawn from an i.i.d.\ sequence of unknown length. We give a similar result when item size distributions are exponential with arbitrary rates. We also study the offline model, where distributions are known in advance but must be packed sequentially. We construct a soft-capacity PTAS for this problem, and show that the complexity of computing the optimal offline cost is $\#\Pol$-hard. Finally, we provide an empirical study of our online algorithm's performance. 
\end{abstract}

\section{Introduction}

Bin Packing is one of the oldest problems in combinatorial optimization, and has been studied by multiple communities in a variety of forms. In the classical online formulation, $n$ items with sizes in $[0,1]$ arrive in an online fashion, and the objective is to pack the items into the fewest possible number of unit-capacity bins. The model has wide applicability in areas including cargo shipping~\cite{wilson1999principles}, assigning virtual machines to servers~\cite{wilcox2011solving}, a variety of scheduling problems \cite{coffman1978application,garey1976resource,vijayakumar2013dual}, and so on. In many of these applications, the items' sizes may be uncertain, with this uncertainty often modeled via probability distributions. 
In much of the stochastic bin packing literature, an item's size is observed before it must be packed, e.g.\ \cite{gupta2015lagrangian,shor1991pack}. Nevertheless, in many applications this assumption is unrealistic. For instance, in bandwidth allocation, connection requests are often bursty and deviate from their typical utilization. If the utilization of the request is higher than expected, it can jeopardize the stability of other connections sharing the same channel. Moreover, the only way to observe the actual traffic required by the connection is to first allocate the request and then observe the traffic pattern.

Motivated by these considerations, we introduce an online adaptive bin packing problem that takes into account the following ingredients:
\begin{enumerate}[leftmargin=*, itemsep=0em]
\item Arrivals are adversarial distributions and the length of the item sequence is unknown to the decision maker.
\item In contrast to existing work in the online and/or stochastic bin packing literature, when an item arrives, the decision maker only observes a probability distribution of its size. 
\item The decision maker observes the item's actual size only after irrevocably placing it in a bin; therefore, overflowing a bin is possible. 
\item An overflowed bin incurs a penalty and renders the bin unusable from that point on. The objective is to minimize the expected cost given by the sum of the number of open bins and overflow penalty.
\end{enumerate}
%


\subsection{Motivating Applications}

The online adaptive bin packing problem captures the uncertainty introduced by the online nature of the problem, and also the uncertainty introduced by learning the size of an item \emph{after} it is packed in a bin. While the variant of the bin packing problem we consider is general and widely applicable, the following examples give some concrete applications: 

\begin{description}[leftmargin=*, itemsep=0em]
	\item[Bandwidth Allocation] An operator is in charge of assigning sequentially arriving independent connection requests. The operator can open new fixed-capacity connections (bins) of unit cost or try to use one of the available connections to pack the incoming request. Traffic on a connection may be bursty, requiring more than the available bandwidth. In this case, the connection suffers from the overflow of the channel, which could represent a monetary penalty or extra work involved in reassigning the request(s) to other connection(s). See also \cite{kleinberg2000allocating}.
	
	\item[Freight Shipping] A dispatcher in a fulfillment center is in charge of packing items into trucks for delivery. Truckloads must comply with a maximum weight limit, and our model applies when the dispatcher assigns items into trucks before their final weighing. An overweight truck incurs a penalty representing additional labor or possible fines.  
	See also \cite{kleywegt_papastavrou98,kleywegt_papastavrou01,papastavrou_etal96}. 
	
	\item[Cloud Computing] A controller is in charge of assigning virtual machines (VM) to servers. The controller has statistical knowledge of the amount of resource a VM will utilize (CPU, RAM, I/O bandwidth, energy, etc.), learned via historical data. The actual resource usage is observed once the VM runs in a server. Excessive consumption of a resource by the VM could compromise the stability of the server and negatively affect other VM's sharing the same infrastructure. See also \cite{gupta2015lagrangian}.
	
	\item[Operation Room Scheduling] In hospitals, an administrator is in charge of assigning incoming surgeries to different operation rooms. There may be a statistical estimation of a procedure's duration, but the real time spent in the room is only learned once the operation has finished. Over-allocating a room could incur economic penalties and loss of patients' good will. See also \cite{denton2010optimal,dexter1999algorithm}.
\end{description}

\subsection{The Model}

We consider the problem of sequentially packing items arriving in an online fashion into homogeneous bins of unit capacity. 
The input consists of a sequence of $n$ nonnegative independent random variables $X_1,\dotsc,X_n$, observed sequentially one at a time. Similar to the bin packing literature, we refer to items interchangeably either by their index $i$ or their corresponding random variable $X_i$. At iteration $i$, random variable $X_i$ arrives and we observe its distribution but not its outcome. We decide irrevocably to pack $X_i$ into an available bin with nonnegative remaining capacity (if any), or to place $X_i$ in a new bin and pay a unit cost. Once packed, we observe the outcome of the random variable $X_i=x_i$, and the chosen bin's capacity is reduced by this amount. A bin overflows when the sizes of items packed in it sum to more than one; when this happens, we incur in an additional cost $C\geq 1$ and the overflowed bin becomes unavailable for future iterations.

We measure the performance of an algorithm $\mathcal{P}$ based on the expected overall cost incurred and denote it $\cost(\mathcal{P})$. 
Because of the online nature of the problem, we cannot expect to compute the optimal cost for an arbitrary sequence of distributions. Even if we knew all distributions in advance, computing the minimum-cost packing is still computationally challenging; the deterministic version reduces to the $\NP$-hard offline bin packing problem. To quantify the quality of an online algorithm, we compare the expected cost incurred by the algorithm against the expected cost incurred by an optimal adaptive packing policy that knows all distributions in advance. This benchmark knows all size distributions in advance but not their outcomes, and must pack the items sequentially in the same order as the online algorithm\footnote{See Section~\ref{sec:thealgorithm} for a more detailed description of policies.}.
This measure of quality differs from the traditional online competitive ratio, cf.\ \cite{albers2003online,borodin2005online}. In the latter, we would compare the performance of an online algorithm against the performance of an extremely powerful optimal offline algorithm that knows all item sizes in advance.

\begin{example}\label{ex:fully_offline_bad}
Consider $n$ i.i.d.\ random variables, 
where $ X_i = 1 $ with probability $1/C$, and $ X_i = 1/n $ with the remaining probability.
We expect $n/C$ random variables to realize to $1$. Therefore, the expected cost of an offline solution that observes the sizes is at most $n/C+1$. In contrast, the cost incurred by any online algorithm (or even an offline algorithm that observes distributions but not sizes) is at least $n$.
\end{example}

Therefore, when measured against the more powerful benchmark, no online algorithm can have a bounded competitive ratio, which motivates us to use a more refined benchmark that knows distributions but not outcomes before the items are packed. In terms of complexity, we show that computing the cost of the optimal offline policy is $\#\Pol$-hard (Theorem~\ref{thm:main_hardness_offline}).

It is worth mentioning that simple greedy strategies based only on a bin's used capacity can perform poorly compared to the optimal offline policy. One such strategy is the \emph{Greedy Algorithm} that compares the instantaneous expected cost of packing the incoming item in an available bin, $C\cdot \Prob(X_i \text{ overflows bin})$, versus the unit cost of opening a new bin, 
selecting the cheapest available choice. This strategy performs poorly in general, even for i.i.d.\ input sequences.

\begin{example}
Consider $n$ i.i.d.\ items, with $X_i \sim \mathrm{Bernoulli}(1/C)$. The optimal policy incurs an expected cost of at most $n/C + 1$: This corresponds to the policy that stops utilizing a bin after observing an item of size $1$. On the other hand, Greedy incurs an expected cost of at least $n/2$, since it will keep trying to pack items in a bin until breaking it. Intuitively, in a sequence of Bernoulli trials the expected time to observe two items of size $1$ is $2/C$; therefore, every $2/C$ items (in expectation), Greedy pays a penalty, incurring an expected cost of roughly $ n/2$.
\end{example}

Another simple choice for a heuristic packing policy is a \emph{Threshold Algorithm}, which establishes a threshold $\alpha\in (0,1)$ such that a bin filled to more than $\alpha$ of its capacity is not used again. Notice that for any $\alpha\in(0,1)$, the optimal policy and the threshold policy incur roughly the same cost for the i.i.d.\ input $X_i\sim \mathrm{Bernoulli}(1/C)$. We now argue that these policies can perform poorly.

\begin{example}
Assume that $\alpha \leq 1/2$ (the case $\alpha>1/2$ is handled similarly) and consider the i.i.d.\ input
\[
X_{i} = \begin{cases}
	0 & \text{w.p.\ } 1-1/C \\
	\alpha & \text{w.p.\ } 1/2C \\
	1-\alpha/2 & \text{w.p.\ } 1/2C .
\end{cases}
\]
The optimal policy incurs an expected cost of at most $n/C+1$, since the policy that stops using a bin upon observing a positive outcome incurs  at most this cost. On the other hand, the Threshold Algorithm incurs an expected cost of at least $n/24-C$; we sketch an argument here to obtain this bound, ignoring the $-C$ term for the sake of clarity: The expected number of positive outcomes is $n/C$. A bin is overflowed by the Threshold Algorithm when an item of size $\alpha$ is followed by another of size $1-\alpha/2$ (regardless of the number of items of size $0$ in between). Focusing solely on the positive outcomes, the number of expected disjoint triplets of the form $(1-\alpha/2,\alpha, 1-\alpha/2)$ is at least a fraction $(1/8) \times (1/3)=1/24$ of these positive outcomes, from which the bound follows.
\end{example}

We include a brief discussion of threshold policies for i.i.d.\ input sequences in Appendix~\ref{sec:threshold_iid}. If the common distribution of the input sequence is finite, a threshold policy can be computed as a function of the distribution, with expected cost a constant factor of the optimal expected cost.


Until now, we have presented examples in which the optimal policies do not break any bin. To not give the false impression that optimal policies do not risk breaking bins, we present the following example.
\begin{example}
Consider $n$ i.i.d.\ items, where $X_i= 1$ with probability $ 1/C^2 $ and $ X_i = 1/n $ with the remaining probability.
The optimal policy has expected cost no more than $n/C+1$, far less than the policy that does not break any bins, which incurs an expected cost of $n$.
\end{example}

In deterministic bin packing problems, one of the most useful bounds for the number of used bins is the sum of the item sizes. It is known that this value is at least half the number of bins used by any greedy algorithm~\cite{man1996approximation}. In our stochastic setting, the expected sum of item sizes could be far from the number of bins used. Indeed, for the random variables considered in Example~\ref{ex:fully_offline_bad}, we have $\sum_{i=1}^n \E[X_i] = n\left( \frac{1}{n}\left( 1- \frac{1}{C}\right) + \frac{1}{C}   \right)= (n-1)/C+ 1$, while the expected cost of any policy is at least $n$ for this input sequence.

\subsection{Our Results and Contributions}


We propose a heuristic algorithm called Budgeted Greedy and denoted $\ALG$ (Algorithm~\ref{alg:BA}). 
Budgeted Greedy uses a \emph{risk budget} in each bin as a way to control the risk of overflowing the bins. If we consider packing item $i$ in bin $j$, this action's risk is equal to the probability of overflowing the bin; Budgeted Greedy maintains a bin's risk below its risk budget. 
At every step, similar to the bin's capacity, when an item is packed in a bin, the bin's risk budget is reduced by the probability of the current item overflowing the bin. If no currently opened bin has enough risk budget left, then a new bin is opened. Observe that the risk of packing item $i$ into any available bin depends on the realized sizes of items $1,\ldots, i-1$ and these items' assignments.

The risk as defined above can be calculated for any policy. While there are instances where the optimal policy incurs a large risk for certain bins it opens, our first structural result shows that \emph{any} policy can be converted to one with budgeted risk with at most a constant factor loss.
\begin{theorem}\label{thm:key_1_informal}
	Let $X_1,\ldots,X_n$ be an arbitrary sequence of independent nonnegative random variables (not necessarily identically distributed). For any $\gamma > 0$ and for any policy $\mathcal{P}$ that sequentially packs $X_1,\dotsc,X_n$, there exists a risk-budgeted policy $\mathcal{P}'$ packing the same items, such that no bin surpasses the risk budget $ \gamma/C $, and with expected cost \[\cost(\mathcal{P}')\leq  ( 1 + 2/\gamma )   \cost(\mathcal{P}).\footnote{If there are items $X_i$ with $\Prob(X_i > 1)> \gamma/C$, these are packed individually. Bins not containing these items have risk bounded by $ \gamma/C$. See Section~\ref{sec:iid_analysis}.}\]
\end{theorem}
Theorem~\ref{thm:key_1_informal} is obtained by updating policy $\mathcal{P}$'s decision tree whenever the risk budget is violated by opening a new bin. The extra cost of the new opened bins is paid by a delicate charging argument. Notice that as $\gamma \to \infty$, we recover the original cost of the policy.

%

While the cost of any policy involves two terms, the expected number of open bins and the expected penalty for overflowed bins, we show (Lemma~\ref{lem:bound}) that for a \emph{budgeted policy}, the cost of overflowed bins is at most the number of opened bins in expectation. This allows us to exclusively focus on the number of bins opened by the budgeted policy. A consequence of these structural results is the following.
\begin{theorem}\label{thm:iid_main_thm}
	If the input sequence $X_1,\dotsc,X_n$ is i.i.d., Budgeted Greedy with $\gamma=\sqrt{2}$ minimizes the expected number of opened bins among all budgeted policies. As a consequence, $\cost(\ALG)\leq (3 + 2 \sqrt{2} )\cost(\OPT)$, where $\OPT$ denotes the optimal policy that knows $n$ in advance.
\end{theorem}

This i.i.d.\ model can be interpreted in the following manner. Suppose there is a probability distribution over the nonnegative real numbers. There are $n$ item sizes independently drawn from this distribution, $x_1,\dotsc,x_n$. For each $i = 1,\dotsc,n$, we are asked to pack the $i$-th item without observing its size. This is indeed a model for basic allocation systems where only a population distribution is known about the item's size, which is a typical occurrence in practical applications if more granular information is not available.

As a consequence of Theorem~\ref{thm:iid_main_thm}, we can also show the existence of instance-dependent threshold policies with similar guarantees as Budgeted Greedy.
\begin{corollary}
	If the input sequence $X_1,\dotsc,X_n$ is i.i.d.\ with finite support, there is a threshold $\alpha\in [0,1]$ that depends on the common distribution of the $X_i$, such that the threshold policy $\mathcal{P}_{\alpha}$, which stops using bins when their capacity exceeds $\alpha$, satisfies $\cost(\mathcal{P}_{\alpha}) \leq (3+ 2 \sqrt{2})\cost(\OPT) + 1$, where $\OPT$ denotes the optimal policy that knows $n$ in advance.
\end{corollary}
The proof is based on the theory of discounted Markov decision processes (see~\cite{puterman2014markov}). We need the finiteness of the support of the distribution to show the existence of a fixed point, which is crucial for the Bellman recursion in the discounted setting. The proof of this corollary appears in the Appendix~\ref{sec:threshold_iid}.

As a second contribution, we show that for arbitrary exponential distributions, i.e.\ a sequence of random variables $X_1,\dotsc,X_n$ with $\Prob(X_i > x)= e^{-\lambda_i x}$, Budgeted Greedy incurs a cost that is at most a factor $\mathcal{O}(\log C)$ times the benchmark cost. Moreover, if the exponential random variables are sufficiently small, this factor can be reduced to a constant.
\begin{theorem}\label{thm:exp_main_bound_thm}
If each $X_i$ is exponentially distributed with rate $ \lambda_i > 0 $, Budgeted Greedy satisfies \[\cost(\ALG)\leq \mathcal{O}(\log C)\cost(\OPT).\] Furthermore, if $\lambda_i\geq 2\log C$ for all $i = 1, \dotsc, n $, $\cost(\ALG)\leq \mathcal{O}(1) \cost(\OPT)$.
\end{theorem}



We show that Budgeted Greedy opens bin $i+1$ if it either packs at least $\Omega(1/\log C)$ in bin $i$, or the risk in bin $i$ is bounded by a small constant, which we obtain via an auxiliary non-convex maximization problem. With this, Budgeted Greedy's cost is bounded by $\mathcal{O}(\log C) \sum_i \E[X_i] \leq \mathcal{O}(\log C) \cost(\OPT)$. 
When the exponential random variables are small enough, 
Budgeted Greedy opens bin $i+1$ if a constant amount of mass in bin $i$ is packed, thereby reducing the $\log C$ factor to a constant. We also give a $\Omega\bigl(\sqrt{\log C}\bigr)$ lower bound for Budgeted Greedy's competitive ratio in the case of exponentially distributed sizes.

\paragraph{Offline Model} Although our motivation for studying the bin packing model is an online application, the offline sequential version of the problem is interesting in its own right, as it interpolates the online setting and the completely offline setting, where items are packed in an arbitrary order. In the offline sequential version of the problem, an ordered list of random variables is given to the decision maker, and the objective is to design a sequential policy to minimize expected cost, in time polynomial in the number of items and possibly $\log C$. As in the online model, right after the decision maker packs a random variable (item) into a bin, the actual size is revealed to her. The optimal offline expected cost computed here corresponds to the benchmark we consider in the online setting.

In this offline framework, we present two main contributions. Following the resource augmentation literature~\cite{fu2018ptas,li2013stochastic}, the first contribution states that there is a polynomial-time approximation scheme (PTAS) for computing a policy when the capacity of the bins is extended by $\varepsilon$.
\begin{theorem}\label{thm:main_soft_ptas_offline}
For any $0< \varepsilon \leq \sqrt{6}(\sqrt{15}- 3)$, there is an algorithm running in $\mathcal{O}\left( n^{2(6/\varepsilon)^5} / \varepsilon^{10} \right)$ time that computes a polynomial-size policy $\mathcal{P}$ packing items into bins of size $1+\varepsilon$, and incurring an expected cost of at most $(1+\varepsilon) \cost(\OPT)$, where $\OPT$ is the optimal policy packing items into bins with unit capacity.
\end{theorem}

The algorithm uses a discretization of possible item sizes similar to~\cite{li2013stochastic}. This allows us to find an optimal policy for discretized outcomes via dynamic programming in polynomial time. The cost of this policy is almost the original optimal cost. We recover a policy for the original items by a tracking argument simulating the discretized policy in parallel. The policy follows the discretized policy's decision to pack items in a bin $j$ as long as the error between the sizes in $j$ and its discrete version remains small. When this fails, the policy opens a new copy of $j$ and keeps following the discretized policy as before. This tracking is enough to guarantee similar cost between the two policies.

Our second result for the offline model relates the complexity of computing the optimal value to counting problems. Specifically, we show that computing the optimal offline cost is $\#\Pol$-hard---hence, the optimal \emph{online} benchmark is also $\#\Pol$-hard to compute.
\begin{theorem}\label{thm:main_hardness_offline}
	It is $\# \Pol$-hard to minimize $\cost(\mathcal{P})$.
\end{theorem}
The proof of this result is divided into two parts. First, we show that counting solutions of symmetric logic formulas in 4CNF\footnote{Logic formula in conjunctive normal form with $4$ literals in each clause.} is $\#\Pol$-hard (Theorem~\ref{thm:sym_sat_hard}). From a symmetric 4CNF formula we construct a stochastic input of the stochastic bin packing problem, where $\min_{\mathcal{P}}\cost(\mathcal{P})$ allows us to count the solutions of the 4CNF formula. The proof resembles the reduction from the Partition problem to the Bin Packing problem. Intuitively, randomized items model outcomes of variables in the 4CNF formula, one item for each positive and negative literal. The main step in the proof is to correlate the outcomes of the positive/negative literals corresponding to the same variable. The proof of Theorem~\ref{thm:main_hardness_offline} is deferred to Appendix~\ref{sec:hardness_offline}.

\subsection{Organization}

The rest of the paper is organized as follows. We follow this introduction with a brief literature review. In Section~\ref{sec:thealgorithm}, we present the Budgeted Greedy algorithm and introduce the necessary notation for the rest of the paper. Section~\ref{sec:iid_analysis} focuses on the i.i.d. case, including the proofs of Theorem~\ref{thm:key_1_informal} and Theorem~\ref{thm:iid_main_thm}. In Section~\ref{sec:expo_analysis} we turn to exponentially distributed item sizes, with the proof of Theorem~\ref{thm:exp_main_bound_thm} and the construction of the corresponding lower bound. Section~\ref{sec:soft_ptas_offline} discusses the offline case, including the proof of Theorem~\ref{thm:main_soft_ptas_offline}. In Section~\ref{sec:numbers} we present a numerical study of our algorithms, comparing it with natural benchmarks. 
\section{Related Work}\label{sec:litrev}

In the classic one-dimensional bin packing problem, $n$ items with sizes $x_1,\ldots,x_n$ in $[0,1]$ must be packed in the fewest unit-capacity bins without splitting any item into two or more bins. 
This is a well-studied $\NP$-complete problem spanning more than sixty years of work~\cite{de1981bin,garey1976resource,gilmore1961linear,hoberg2017logarithmic,johnson1974fast,karmarkar1982efficient,rothvoss2013approximating}. For excellent surveys see~\cite{christensen2016multidimensional,man1996approximation}. In the online version, the list of items $L=(x_1,\ldots,x_n)$ is revealed online one item at a time. At round $t$, we observe item $x_t$ and we need to decide irrevocably and without knowledge of future arrivals whether to pack the item in an open bin with enough remaining space, or to open a new unit-capacity bin at unit cost. 
It is standard to measure an online algorithm's performance via its \emph{(asymptotic) competitive ratio}~\cite{albers2003online,borodin2005online,man1996approximation} $\limsup_{|L|\to \infty} \cost_{\text{alg}}(L)/\cost_{\text{OPT}}(L)$, where $\cost_\text{alg} (L)$ is the cost incurred by the online algorithm with input $L$, and $\cost_\text{OPT}(L)$ is the cost incurred by the optimal offline solution that knows $L$ in advance. The best known competitive ratio is 1.57829 \cite{balogh2018new}, and the best current lower bound is 1.5403~\cite{balogh2012new}, see also~\cite{van1992improved,yao1980new}.

In several real-world applications, exact item sizes are unknown to the decision maker at the time of insertion \cite{della2016probabilistic,vijayakumar2013dual}. This uncertainty is typically modeled via probability distributions on the items' size. Several online and offline bin packing models introducing stochastic components have been studied \cite{coffman2000bin,csirik2006sum,goel1999stochastic,gupta2012online,gupta2020interior,kleinberg2000allocating,li2013stochastic,rhee1988optimal,shor1991pack,shor1986average}. These stochastic models have revealed connections with balls-into-bins problems \cite{shor1991pack}, sums of squares \cite{csirik2006sum}, queuing theory \cite{coffman1980stochastic}, Poisson approximation \cite{li2013stochastic}, etc. For the online case, common to all these models is the assumption that the item size is observed before packing it. Nevertheless, observing the item size is unrealistic in many scenarios. For instance, in cloud computing, before running a job in a cluster, we may have some statistical knowledge of the amount of resource the job will utilize. However, the only way to observe the real utilization is to start the job. In this work, we propose a new model variant where items' size distributions are revealed in an online fashion but each outcome is observed only after packing the item. We therefore relax the strict capacity constraint by allowing each bin to overflow at most once, at the expense of a penalty. 
Related to this kind of online input are the works \cite{coffman2000bin,csirik2006sum,gupta2012online,rhee1988optimal,shor1991pack,shor1986average}. 

Our model also shares similarities with adaptive combinatorial optimization, particularly stochastic knapsack models introduced in \cite{derman1978renewal}. Recent treatments began with \cite{dean2008approximating}; a large body of work has now studied this model from several perspectives 
\cite{balseiro_brown19,bhalgat2011improved,blado_etal16,blado_toriello19,blado_toriello20,fu2018ptas,gupta2011approximation,li2013stochastic,ma2018improvements}. Most of these works assume complete knowledge of the input distributions, and online treatments are scarcer in the literature, see \cite{alaei2013online,goyal2019online,mehta2014online}.

A related area of work is the extensible bin packing problem~\cite{coffman2001approximation,dell19981312}. Roughly speaking, a fixed number of bins are given and a set of items must be packed into them. The cost of bin $B$ corresponds to $\max\left\{ \sum_{i\in B} x_i, 1 \right\}$, a fixed unit cost and a linear excess cost. The objective is to design packings with small overall cost; even though we do not allow bins to be utilized after overflow, we could interpret our model as a nonlinear version of a stochastic extensible bin packing problem. For a generalization to different costs and bin capacities, see~\cite{levin2019approximation}. For a stochastic approach similar to our posterior observability, see~\cite{sagnol2018price}.

\section{The Algorithm}\label{sec:thealgorithm}

\subsection{Preliminaries}


The problem's input consists of $n$ independent nonnegative random variables $X_1,\ldots,X_n$. The (possible) bins to utilize are denoted by $B_1,B_2,\ldots,B_n$. A \emph{state} $\state$ for round $i\in[n+1]$ is a sequence $(x_1,1\to j_1)(x_2,2\to j_2)\cdots(x_{i-1},i-1\to j_{i-1})$, where $x_k$ is an outcome of $X_k$ for all $k < i$. The pair $(x_k,k\to j)$ represents round $k$, and refers to packing $X_k$ in bin $j$ and observing outcome $X_k=x_k$. A state for round $i$ represents the path followed by a decision maker packing items $X_1,\ldots,X_n$ sequentially into bins and the outcomes for each of these decisions until round $i-1$. States have a natural recursive structure: $\state = \state' (x_{i-1}, i-1\to j_{i-1})$, where $\state'$ is the state for round $i-1$. The \emph{initial state} $\state_0$ is the empty state. Bin $B_j$ is \emph{open} by state $\state$ if some $(x_k,k\to j)$ appears in $\state$. The items packed into bin $B_j$ by state $\state$ are $B_j(\state) = \{ k : (x_k,k\to j) \text{ appears in }\state  \}$. The number of bins opened by state $\state$ is $|\{ j : (x_k, k\to j) \text{ appears in }\state  \}|$. The \emph{usage} of bin $B_j$ at the beginning of round $i$ in state $\state$ is
\[
S_j^{i-1}(\state)=\sum_{\!\!\!\substack{k \leq i-1 \\ (x_k, k \to j)\in \state} } x_k ,
\]
the sum of sizes of items packed in bin $j$. A bin $B_j$ is \emph{broken} or \emph{overflowed} in $\state$ if $S_j^{i-1}(\state) > 1$. In our model, we stop using bins that overflow. A state $\state$ for round $i$ is \emph{feasible} if any overflowed bin by round $k$ is never used again after $k$, for any $k < i$. The \emph{state space} is the set of \emph{all feasible states}, denoted $\State$. The set of all feasible states for round $i\leq n$ is denoted by $\State_i$.

A \emph{policy} $\mathcal{P}$ is a function $\mathcal{P}: \State_n \to [n]$ such that for a feasible state $\state\in \State_n$ for round $i$, $\mathcal{P}(\state)=j$ indicates that item $i$ is packed into bin $j$; we write this as $i\to j$ when the policy and state are clear from the context. The policy is \emph{feasible} if $\state'=\state(x_i,i\to \mathcal{P}(\state))$ is a feasible state for any feasible $\state \in \State_n$ for round $i$ and outcome $x_i$ of $X_i$. From now on, we only consider feasible policies. A state $\state' \in \State$ is \emph{reachable} by the policy if $\state'=\state_0$ or $\state' =\state (x_i,i\to \mathcal{P}(\state))$ with $\state$ reachable, $\state$ for round $i$ and $x_i$ an outcome of $X_i$. For a reachable state $\state$ for round $i\in [n]$, we say that $\mathcal{P}$ \emph{opens} bin $j$ if $\mathcal{P}(\state)=j$ and $B_j$ is not open in $\state$. We say that the policy \emph{overflows} bin $B_j$ at state $\state$ if $B_j$ overflows for $\state'=\state(x_i,i\to \mathcal{P}(\state))$ but $B_j$ is not overflowed in $\state$. We set the \emph{cost} of a policy as
\[
\cost(\mathcal{P}) = \E[N_\mathcal{P}] + C\E[O_\mathcal{P}],
\]
where $N_\mathcal{P}$ is the number of bins opened and $O_\mathcal{P}$ is the number of bins broken by reachable states for round $n+1$. The randomness is over the items' outcomes. Notice that non-reachable states in $\State$ are unimportant for $\cost(\mathcal{P})$, hence we can always assume $\mathcal{P}(\state)=n$ for non-reachable $\state\in \State_n$. A policy specifies the actions to apply in any epoch of the sequential decision-making problem. Note that our states are typically considered \emph{histories} in the Markov decision processes literature~\cite{puterman2014markov}. We use our description of states to keep close track of policies' actions in the subsequent analysis.

Any policy $\mathcal{P}$ has a natural $(n+1)$-level decision tree representation $\mathcal{T}_\mathcal{P}$, which we call the \emph{policy tree}. The root, denoted $r$, is at level $1$ and represents item $X_1$ and state $\state_0$. A node at level $i\in [n]$ is labeled with $\mathcal{P}(i,\state)$ where $\state$ is the state of the system obtained by following the path from the root to the current node. There is a unique arc going out of the node for every possible outcome of $X_i$ directed to a unique node in level $i+1$. Nodes at level $n+1$ are \emph{leaves} denoting that the computation has ended. Nodes in levels $i\in [n]$ are called \emph{internal nodes}. To compute the $\cost(\mathcal{P})$ using the policy-tree $\mathcal{T}_\mathcal{P}$, we add two labels to the tree:
\begin{itemize}[leftmargin=*, itemsep=0em]
	\item For an internal node $u$, $\ell_u= 1$ if $\mathcal{P}$ opens a new bin in node $u$; $0$ otherwise. For leaves we define $\ell_u=0$.
	
	\item For arcs $a=(u,v)$, we define $c_{a}= C$ if the outcome of the random variable belonging to the level where $u$ is located overflows the bin chosen by the policy at node $u$; $0$ otherwise.
\end{itemize}
We refer to this tree as cost-labeled tree $\mathcal{T}_\mathcal{P}$ with cost vectors $(\ell,c)$, or simply cost-labeled tree $\mathcal{T}_\mathcal{P}$ if the costs are clear from the context. The tree structure gives us a recursive way of computing the cost of the policy. Let $\mathcal{T}_\mathcal{P}(u)$ be the cost-labeled sub-tree of $\mathcal{T}_{\mathcal{P}}$ rooted at node $u$; then
\[
\cost_{\ell,c}(\mathcal{T}_\mathcal{P}(u)) =\begin{cases}
 \ell_{u} + \E_{X_i}[c_{ (u,u_{X_i} ) } +  \cost_{\ell,c}(\mathcal{T}_{\mathcal{P}}(u_{X_i}))]  & \text{if } u \text{ is at level } i =1,\ldots,n \\
 0 &\text{if } u \text{ is at level }n+1
\end{cases},
\]
where $u_{X_i}$ is the node at level $i+1$ connected to $u$. Thus, $\cost(\mathcal{P}) = \cost_{\ell,c}(\mathcal{T}_\mathcal{P}(r))$. We define $
\OPT = \argmin_{\mathcal{P}} \cost(\mathcal{P})$ as the optimal policy for sequentially packing items $X_1,\ldots,X_n$. This policy might not exist in cases where the number of states is uncountable, for example, when $X_1,\dotsc,X_n$ have continuous distributions. In this case, the policy tree has uncountably many edges emanating from nodes, corresponding to all possible realizations of $X_i$. Nevertheless, a $\varepsilon$-optimal policy is guaranteed to exist, i.e.\ a policy $\mathcal{P}$ that ensures $\cost(\mathcal{P}) \leq \inf_{\mathcal{P}} \cost(\mathcal{P}) + \varepsilon$. We abuse notation by calling $\OPT$ the optimal policy (or an arbitrarily good approximation if it does not exist).

Note that we defined only \emph{deterministic} policies, since the action $\mathcal{P}(\state)$ is deterministic. If $\mathcal{P}(\state)$ were a probability distribution over $[n]$, then we would have a randomized policy. A standard result from Markov decision processes theory ensures that any randomized policy has a deterministic counterpart incurring the same cost; hence, we only focus on deterministic policies. For more details see~\cite{powell2007approximate,puterman2014markov}.

The following proposition characterizes the expected number of bins overflowed by a policy. The proof appears in Appendix~\ref{sec:Appendix}.
\begin{proposition}\label{prop:breakbinpenalty}
	Let $X_1,\ldots,X_n$ be nonnegative independent random variables, and let $\mathcal{P}$ be any policy that sequentially packs these items. The expected number of bins broken by the policy $\mathcal{P}$ is
	\[
	\E[O_\mathcal{P}] = \sum_{j=1}^n \E_{X_1,\ldots,X_n} \left[  \sum_{i=1}^n \Prob_{X_i}(X_i + S_{j}^{i-1} > 1)\mathbf{1}_{\{ i \to j \}}^{\mathcal{P}}  \right] ,
	\]
	where $S_{j}^{i-1}$ is the usage of bin $j$ at the beginning of iteration $i$ and $\mathbf{1}_{\{i\to j\}}^\mathcal{P}$ is the indicator random variable of the event in which $\mathcal{P}$ packs item $X_i$ into bin $j$.
\end{proposition}
If we interpret $\Prob_{X_i}( X_i + S_j^{i-1} > 1 )$ as the \emph{risk} that $X_i$ overflows bin $j$ if packed there, the result says that the number of overflowed bins is the expected aggregation of these risks. We define the \emph{risk} of a bin $j$ as $\risk(B_j)=\sum_{i=1}^n\Prob_{X_i}(X_i + S_{j}^{i-1} > 1)\mathbf{1}_{\{ i \to j \}}$. Then $\E[O_\mathcal{P}] = \sum_{j=1}^n \E[\risk(B_j)]$.
%
%
A policy $\mathcal{P}$ is \emph{risk-budgeted} or simply \emph{budgeted} with \emph{risk budget} $r>0$ if no bin incurs a risk larger than $r$, $\risk(B_j) \leq r$ for $ j \in [n] $.

A deterministic online algorithm induces a policy, with non-reachable states simply mapped to $\emptyset$. Since online algorithms are not aware of the number of items $n$, we label the $j$-th bin opened by an online algorithm as $B_j$ in this case. The cost of an online algorithm is naturally defined as the cost of the corresponding induced policy.

We use the notation $z(B)= \sum_{i\in B} z_i$ for a vector $z=(z_1,\ldots,z_n)$. If $X=(X_1,\ldots,X_n)$ is the vector of random variables and $B=B_j$, then $X(B)= S_j^n$ is the usage of bin $B_j$. The following propositions are probabilistic analogues of the well-known size lower bound for deterministic bin packing. We use them in Sections~\ref{sec:expo_analysis} and~\ref{sec:soft_ptas_offline}. The proofs are deferred to Appendix~\ref{sec:Appendix}.
\begin{proposition}\label{prop:sizeinbins}
	For any sequence of nonnegative i.i.d.\ random variables $X_1,\ldots,X_n$, for any bin $B=B_j$ and any policy $\mathcal{P}$, we have
	\[
	\E\left[  \sum_{i\in B} \E [X_i\wedge 1] \right] = \E\left[  \sum_{i\in B} ( X_i\wedge 1) \right]  \leq 2 \Prob(\mathcal{P}\text{ opens bin } B),
	\]
	where $X_i\wedge 1 = \min\{X_i,1\}$.
\end{proposition}
When all items sizes are aggregated, we can improve the factor of $2$ as follows.
\begin{proposition}\label{prop:sizesinbins_cost}
	For any sequence of nonnegative i.i.d.\ random variables $X_1,\ldots,X_n$, for any policy, we have
	\[
	\E\left[ \sum_{i=1}^n (X_i \wedge 1 )\right] \leq \cost(\mathcal{P}).
	\]
\end{proposition}

\subsection{The Budgeted Algorithm}\label{subsec:algorithm}


In the \emph{Budgeted Greedy} algorithm, we keep a risk budget for each bin that is initialized as $\gamma/C$, where $\gamma\geq 1$ is an algorithm parameter. We pack items in a bin as long as the usage of the bin is at most $1$ and its risk budget has not run out. More formally, when opening a bin, say bin $j$ at round $i$, we initialize its risk of overflow at $r_j^{i-1}=0$. At round $i$, when item $X_i$ arrives, we find a bin $j$ such that $r_j^{i-1} + p_i(S_j^{i-1})\leq \gamma/C$, where $r_j^{i-1}$ is the accumulated risk of overflowing the bin until $i-1$, $S_j^{i-1}$ is the usage of the bin $j$ until the previous round and $p_i(S_{j}^{i-1}) = \Prob_{X_i}( X_i + S_j^{i-1} > 1 )$ is the risk that $X_i$ overflows bin $j$. If that bin $j$ exists, we pack the incoming item into bin $j$, breaking ties arbitrarily, and we update the risk of overflow as $r_j^{i}= r_j^{i-1} + p_i(S_j^{i-1}) $ and $r_{j'}^{i}= r_{j'}^{i-1}$ for any $j'\neq j$. Such a bin may not exist, in which case we open a new bin $k$ with $r_k^i = p_i(0)$. Strictly speaking, Budgeted Greedy is not a budgeted policy with risk budget $\gamma/C$ unless all items satisfy $\Prob(X_i > 1)\leq \gamma/C$; items with $\Prob(X_i > 1) > \gamma/C$ are packed into individual bins. In Algorithm~\ref{alg:BA}, we formally present the description of Budgeted Greedy.

\begin{algorithm}[H]
	\small
	\setlength{\skiptext}{12pt}
	\setlength{\skiprule}{7pt}
	\SetAlgoLined
	Initialize: $I=\emptyset$.\\
	\For{$i=1\ldots,n$}{
		\uIf{$\exists j\in I$ such that $r_j^{i-1} + p_i(S_{j}^{i-1})\leq \gamma/C$}{
			$S_j^{i}= S_j^{i-1}+ X_i$.\\
			$r_j^{i}= r_j^{i-1} + p_i(S_{j}^{i-1})$.\\
		}
		\Else{
			Define $r_j^{i}=p_i(0)$ for $j$ such that $j=\inf\{ j\geq 0: j\notin I \}$.\\
			$S_j^i = X_i$.\\
			Update $I=I\cup\{j\}$.\\
		}
		\For{$j'\neq j$}{
			$S_{j'}^{i}= S_{j'}^{i-1}$.\\
			$r_{j'}^{i}= r_{j'}^{i-1}$.
		}
	}
	\caption{$\textsc{Budgeted-Greedy}(\gamma,X_1,\ldots,X_n)$}
	\label{alg:BA}
\end{algorithm}

\begin{lemma}\label{lem:bound}
	Let $\gamma\geq 1$ and assume that for all $i$, $\Prob(X_i > 1) \leq \gamma/C$. For any bin $j$, Algorithm~\ref{alg:BA} guarantees
	$$\Prob(\ALG \text{ breaks bin }j)\leq \frac{\gamma}{C}\Prob(\ALG \text{ opens bin }j).$$
\end{lemma}
\proof{Proof.}
	Using Proposition~\ref{prop:breakbinpenalty},
	\begin{align*}
	\Prob(\ALG \text{ breaks bin }j) & = \E\left[\left(\sum_{i=1}^n \Prob(X_i + S_j^{i-1}>1) \mathbf{1}_{\{ t\to j \}}^\ALG \right)\mathbf{1}_{\{ \ALG \text{ opens bin }j \}} \right]\\
	& = \E\left[\risk(B_j) \mathbf{1}_{\{ \ALG \text{ opens bin }j \}}\right]\leq \frac{\gamma}{C}\Prob(\ALG \text{ opens bin }j),
	\end{align*}
since once the bin has been opened, its risk never goes beyond $\gamma/C$.
	\hfill\halmos
\endproof

As a result, we have the following corollary, which implies that we only need to bound the expected number of bins opened by Budgeted Greedy in our analysis.
\begin{corollary}\label{cor:boundgamma}
	Under the same assumptions as Lemma~\ref{lem:bound}, $\cost(\ALG)\leq (1+\gamma)\E[N_\ALG]$.
\end{corollary}

\section{A Policy-Tree Analysis for I.I.D.\ Random Variables}\label{sec:iid_analysis}

In this section we prove Theorem~\ref{thm:key_1_informal} for general input distributions in Theorem~\ref{thm:key_1}. We use this result to prove Theorem~\ref{thm:iid_main_thm}, which gives the Budgeted Greedy guarantee for the i.i.d.\ case. Theorem~\ref{thm:key_1} states that any policy can be converted into a budgeted version, where a risk budget is never surpassed for any bin. This transformation can be carried out while only incurring a small multiplicative loss. The proof relies on a charging scheme in the cost paid by overflowing bins. Starting with the original policy tree, we increase the cost paid by overflowing bins by an amount $\delta >0$. The overall cost of the tree increases multiplicatively by at most $(1+\delta/C)$. We show that this additional $\delta$ allows us to pay for new bins whenever the risk of the bin goes beyond $\gamma/C$, for an appropriate choice of $\delta$ and $\gamma$.

%

\begin{theorem}\label{thm:key_1}
	
	Let $X_1,\dotsc,X_n$ be an arbitrary sequence of independent, nonnegative random variables that are not necessarily identical. Fix $\gamma > 0$. For any policy $\mathcal{P}$ that sequentially packs items $X_1,\dotsc,X_n$, there exists a policy $\mathcal{P}'$ for the same items such that:
	\begin{itemize}[leftmargin=*, itemsep=0em]
		\item $\mathcal{P}'$ packs items with $\Prob(X_i > 1) > \gamma/C$ into individual bins, and bins not containing these items never exceed the risk budget $\gamma/C$.
		
		\item $\mathcal{P}'$ satisfies $\cost(\mathcal{P}')\leq (1 + 2/\gamma ) \cost(\mathcal{P})$.
	\end{itemize}
	In particular, if all items satisfy $\Prob(X_i > 1)\leq \gamma/C$, $\mathcal{P}'$ is a risk-budgeted policy with risk budget $\gamma/C$.
\end{theorem}
\proof{Proof.}
	The proof follows two phases. In the first and longest phase, we show that we can modify the policy $\mathcal{P}$ in such a way that the risk of each bin exceeds $\gamma/C$ at most once. In the second phase, we show that the item surpassing the risk budget in each bin can be packed into an individual bin. At the end, no item with $\Prob(X_i > 1)\leq \gamma/C$ can exceed the risk $\gamma/C$.

	In the rest of the proof we utilize the tree representation of the policy. Let $\delta = C/\gamma > 0$. We proceed as follows:
	\begin{enumerate}[leftmargin=*, itemsep=0em]
		\item In the cost labeled tree $\mathcal{T}_\mathcal{P}$, increase the cost of overflowing the bins from $C$ to $C+2\delta$. That is, $\widehat{c}_{(u,v)} = C+2\delta$ if $c_{(u,v)} =C$ and $0$ otherwise for any arc $(u,v)$ in $\mathcal{T}_{\mathcal{P}}$. Then,
		\[
		\cost_{\ell,\widehat{c}}(\mathcal{T}_\mathcal{P}(r)) \leq \left(1+2\frac{\delta}{C} \right) \cost_{\ell ,c}(\mathcal{T}_\mathcal{P}(r))= \left( 1+ 2\frac{\delta}{C} \right)\cost(\mathcal{P}).
		\]
		
		\item Starting at the root of this new cost-labeled tree, find a node $u$ at level $i=1,\ldots,n$ where the policy $\mathcal{P}$ decides to open a new bin, say bin $j$. In each of the branches starting at node $u$ and directed to some leaf, find the sequence of nodes $u_1=u, u_2, \ldots, u_k$ where the policy packs items into bin $j$ and node $u_k$ corresponds to the first node in the branch where the risk budget $\gamma/C$ is surpassed for bin $j$. Define $u_k$ as a leaf if in the branch the risk budget is not surpassed for bin $j$. Let $i_1=i, i_2, \ldots, i_k$ be the items packed into bin $j$ in this branch; that is, node $u_\ell$ is at level $i_\ell$. Then, we have
		\[
		\sum_{m = 1}^{k-1} \Prob_{X_{i_m}}(X_{i_m} + S_{j}^{i_m-1} > 1 ) \leq \frac{\gamma}{C}, \quad \text{ and } \quad
		\sum_{m =1 }^{k} \Prob_{X_{i_m}}(X_{i_m} + S_j^{i_m - 1} > 1) > \frac{\gamma}{C}
		\]
		if node $u_k$ is not a leaf. Here $S_j^{i_m-1}$ represents the usage of bin $j$ at node $i_m$.
		
		Consider the following modifications to the cost-labeled tree $\mathcal{T}_{\mathcal{P}}$: We start with the same tree as $\mathcal{P}$ but in the subtree rooted at $u$, bin $j$ is utilized only in nodes $u_1,\ldots,u_k$ for the different branches. Any future utilization of bin $j$ after passing through node $u_k$ is moved to a new bin $j'$. Now, we update the cost labels as follows. For all the branches, we reduce the cost of $C+2\delta$ appearing in the arcs going out from nodes $u_1,\ldots,u_k$ to $C+\delta$. We label the first node appearing after node $u_k$ where the bin $j'$ is opened with a $1$. We reduce the labels of arcs going out of nodes using bin $j'$ if they do not overflow the bin $j'$ anymore (bin $j'$ has smaller usage than bin $j$). Formally, for any branch and nodes $u_1,\ldots,u_k$ defined as before,
		$$c_{(a,b)}' = \begin{cases}
		\frac{C+\delta}{C+2\delta}\widehat{c}_{(a,b)} & a=u_m \text{ for some branch starting at }u \\
		\widehat{c}_{a,b} & \text{otherwise}
		\end{cases},$$
		and for nodes,
		\[
		\ell_{a}'  = \begin{cases}
		1 & a\text{ is the first node packed into bin }j' \text{ in the subtree }\mathcal{T}_\mathcal{P}(u_k) \\
		\ell_a &\text{otherwise}
		\end{cases}.
		\]
		We denote this new policy by $\mathcal{P}'$. Figure~\ref{fig:policy_tree_modification} displays the modification process.

		\begin{figure}[h!]
			\centering
			\resizebox{\columnwidth}{!}{
			\begin{minipage}[t]{0.49\linewidth}
				\begin{tikzpicture}[thick, scale=0.45,
				vertex_style/.style={draw=blue!60,minimum size=0cm,scale=0.8,very thick,fill=blue!10!white},
				edge_style/.style={black}]
				
				\useasboundingbox  (0,0) rectangle (18,18);
				
				\path[-,black,draw] (0,0) -- (18,0) -- (9,15.58) -- cycle;
				\small
				
				\node(TP) at (3,16) {$\mathcal{T}_{\mathcal{P}}$};
				
				\node(root) at (9,15.58) [draw,circle,black,minimum size=1pt,inner sep=1pt,fill=white] {$~r~$};

				\begin{scope}[xshift=11cm,  yshift=12.12cm]
				\path[thick,-stealth ]
				(0,0) edge  node [left] {$0$} (-1,-1.5)
				(0,0) edge  node [right] {$0$} (1,-1.7)
				;
				\path[thick, dashed, bend right=30] (-1,-1.5) edge node [] {} (1,-1.7);
				\end{scope}
				
				\begin{scope}[xshift=9cm,  yshift=9.12cm]
				\path[thick,-stealth ]
				(0,0) edge  node [left] {$C+2\delta$} (-1,-1.7)
				(0,0) edge  node [right] {$0$} (1.15,-1.73)
				;
				\path[thick, dashed, bend right=30] (-1,-1.7) edge node [] {} (1,-1.7);
				\end{scope}
				
				\begin{scope}[xshift=11cm,  yshift=6.12cm]
				\path[thick,-stealth ]
				(0,0) edge  node [left] {$C+2\delta$} (-1,-2.06)
				(0,0) edge  node [right] {$0$} (1,-2.06)
				;
				\path[thick, dashed, bend right=30] (-1,-2.06) edge node [] {} (1,-2.06);
				\end{scope}
				
				\begin{scope}[xshift=4cm,  yshift=6.92cm]
				\path[thick,-stealth ]
				(0,0) edge  node [left] {$C+2\delta$} (-1,-1.7)
				(0,0) edge  node [right] {$0$} (1,-1.7)
				;
				\path[thick, dashed, bend right=30] (-1,-1.7) edge node [] {} (1,-1.7);
				\node(v) at (0,0) [circle,draw,minimum size=1pt,inner sep=1pt,fill=white] {$~v~$};
				\end{scope}

				\node(u1) at (11,12.12) [draw,circle,black,minimum size=1pt,inner sep=1pt,fill=white] {$u_1$};
				\node(u2) at (9,9.12) [draw,circle,black,minimum size=1pt,inner sep=1pt,fill=white] {$u_2$};
				\node(uk) at (11,6.12) [draw,circle,black,minimum size=1pt,inner sep=1pt,fill=white] {$u_k$};
				\node(a) at (11.5,3) [draw,circle,black,minimum size=1pt,inner sep=1pt,fill=white] {$~a~$};
				
				\def\x{14cm}
				\draw (\x,12.12) node(u11) [right = 0.7cm,rectangle,fill=white, draw] {$\substack{\ell_{u_1}=1\\i_1\to j}$};
				\draw (\x,9.12) node(u22) [right = 0.7cm,rectangle,fill=white, draw] {{$\substack{\ell_{u_2}=0\\i_2\to j}$}};
				\draw (\x,6.12) node(ukk) [right = 0.7cm, rectangle, draw, fill=white] {{$\substack{\ell_{u_k}=0\\i_k\to j}$}};
				\draw (\x,3) node(aa) [right = 0.7cm, rectangle, draw, fill=white] {{$\substack{\ell_{a}=0~\\~i\to j~}$}};
				
				\path[thick,-, shorten >= 1pt]
				(u1) edge [dashed] node {} (u2)	
				(u2) edge [dashed] node {} (uk)	
				(uk) edge [dashed] node {} (14,0)
				(uk) edge [dashed] node {} (8,0)
				(u1) edge [thin] node {} (u11)
				(u2) edge [thin] node {} (u22)
				(uk) edge [thin] node {} (ukk)
				(a)  edge [thin] node {} (aa)
				;
				
				\begin{scope}[xshift=11.5cm,  yshift=3cm]
				\path[thick,-stealth ]
				(0,0) edge  node [left] {\colorbox{white}{$C+2\delta$}} (-1,-1.7)
				(0,0) edge  node [right] {$0$} (1,-1.7)
				;
				\path[thick, dashed, bend right=30] (-1,-1.7) edge node [] {} (1,-1.7);
				\node(v) at (0,0) [circle,draw,minimum size=1pt,inner sep=1pt,fill=white] {$~v~$};
				\end{scope}

				\end{tikzpicture}
			\end{minipage}~
		\vline~
			\begin{minipage}[t]{0.49\linewidth}
				\begin{tikzpicture}[thick, scale=0.45,
				vertex_style/.style={draw=blue!60,minimum size=0cm,scale=0.8,very thick,fill=blue!10!white},
				edge_style/.style={black}]
				
				\useasboundingbox  (0,0) rectangle (18,18);
				
				\path[-,black,draw] (0,0) -- (18,0) -- (9,15.58) -- cycle;
				\small
				
				\node(TP) at (3,16) {$\mathcal{T}_{\mathcal{P}'}$};
				
				\node(root) at (9,15.58) [draw,circle,black,minimum size=1pt,inner sep=1pt,fill=white] {$~r~$};

				\begin{scope}[xshift=11cm,  yshift=12.12cm]
				\path[thick,-stealth ]
				(0,0) edge  node [left] {$0$} (-1,-1.5)
				(0,0) edge  node [right] {$0$} (1,-1.7)
				;
				\path[thick, dashed, bend right=30] (-1,-1.5) edge node [] {} (1,-1.7);
				\end{scope}
				
				\begin{scope}[xshift=9cm,  yshift=9.12cm]
				\path[thick,-stealth ]
				(0,0) edge  node [left] {$C+\delta$} (-1,-1.7)
				(0,0) edge  node [right] {$0$} (1.15,-1.73)
				;
				\path[thick, dashed, bend right=30] (-1,-1.7) edge node [] {} (1,-1.7);
				\end{scope}
				
				\begin{scope}[xshift=11cm,  yshift=6.12cm]
				\path[thick,-stealth ]
				(0,0) edge  node [left] {$C+\delta$} (-1,-2.06)
				(0,0) edge  node [right] {$0$} (1,-2.06)
				;
				\path[thick, dashed, bend right=30] (-1,-2.06) edge node [] {} (1,-2.06);
				\end{scope}
				
				\begin{scope}[xshift=4cm,  yshift=6.92cm]
				\path[thick,-stealth ]
				(0,0) edge  node [left] {$C+2\delta$} (-1,-1.7)
				(0,0) edge  node [right] {$0$} (1,-1.7)
				;
				\path[thick, dashed, bend right=30] (-1,-1.7) edge node [] {} (1,-1.7);
				\node(v) at (0,0) [circle,draw,minimum size=1pt,inner sep=1pt,fill=white] {$~v~$};
				\end{scope}

				\node(u1) at (11,12.12) [draw,circle,black,minimum size=1pt,inner sep=1pt,fill=white] {$u_1$};
				\node(u2) at (9,9.12) [draw,circle,black,minimum size=1pt,inner sep=1pt,fill=white] {$u_2$};
				\node(uk) at (11,6.12) [draw,circle,black,minimum size=1pt,inner sep=1pt,fill=white] {$u_k$};
				\node(a) at (11.5,3) [draw,circle,black,minimum size=1pt,inner sep=1pt,fill=white] {$~a~$};
				
				\def\x{14cm}
				\draw (\x,12.12) node(u11) [right = 0.7cm,rectangle,fill=white, draw] {$\substack{\ell_{u_1}'=1\\i_1\to j}$};
				\draw (\x,9.12) node(u22) [right = 0.7cm,rectangle,fill=white, draw] {{$\substack{\ell_{u_2}'=0\\i_2\to j}$}};
				\draw (\x,6.12) node(ukk) [right = 0.7cm, rectangle, draw, fill=white] {{$\substack{\ell_{u_k}'=0\\i_k\to j}$}};
				\draw (\x,3) node(aa) [right = 0.7cm, rectangle, draw, fill=white] {{$\substack{\ell_{a}'=1~\\~i\to j'~}$}};
				
				\path[thick,-, shorten >= 1pt]
				(u1) edge [dashed] node {} (u2)	
				(u2) edge [dashed] node {} (uk)	
				(uk) edge [dashed] node {} (14,0)
				(uk) edge [dashed] node {} (8,0)
				(u1) edge [thin] node {} (u11)
				(u2) edge [thin] node {} (u22)
				(uk) edge [thin] node {} (ukk)
				(a)  edge [thin] node {} (aa)
				;
				
				\begin{scope}[xshift=11.5cm,  yshift=3cm]
				\path[thick,-stealth ]
				(0,0) edge  node [left] {\colorbox{white}{$C+2\delta$}} (-1,-1.7)
				(0,0) edge  node [right] {$0$} (1,-1.7)
				;
				\path[thick, dashed, bend right=30] (-1,-1.7) edge node [] {} (1,-1.7);
				\node(v) at (0,0) [circle,draw,minimum size=1pt,inner sep=1pt,fill=white] {$~v~$};
				\end{scope}

%
%
%
%
%
%
%
%
%
%
%
%
%
				\end{tikzpicture}
			\end{minipage}}
			
			\caption{Policy tree modification. On the left, we display the original tree with augmented cost from $C$ to $C+2\delta$. On the right, we show the modified labels after opening a new bin in node $u_k$. Observe that we only decrease the costs of arcs related to bin $j$ going out of nodes $u_1,\ldots,u_{k}$ in all branches starting at node $u$.}\label{fig:policy_tree_modification}
\end{figure}
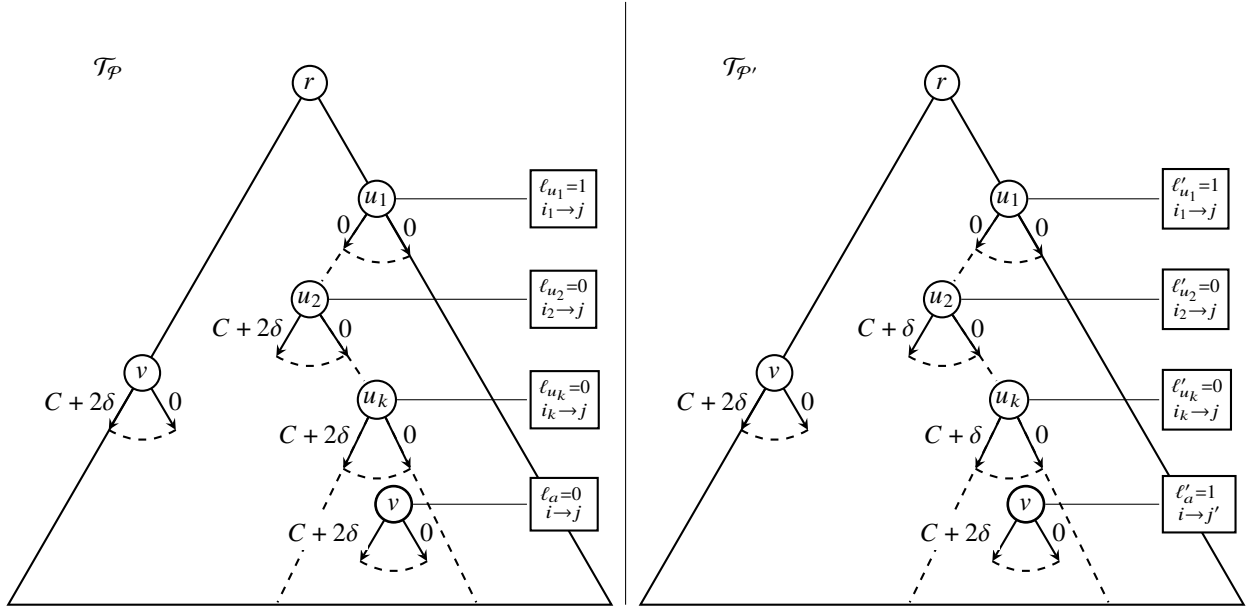

		We now argue that the changes applied to the cost-labeled tree $\mathcal{T}_{\mathcal{P}}$ to transform it into $\mathcal{T}_{\mathcal{P}'}$ do not increase the cost function $\cost_{\ell,\widehat{c}}(\mathcal{T}_\mathcal{P})$. Since we only modified labels in the subtree $\mathcal{T}_{\mathcal{P}}(u)$ it is enough to study the cost change in this specific subtree for bins $j$ and bin $j'$.
		
		\begin{lemma}\label{lem:cost_reduce_tree_analysis}
			$\cost_{\ell, \widehat{c}}(\mathcal{T}_\mathcal{P}(u)) \geq \cost_{\ell',c'}(\mathcal{T}_{\mathcal{P'}}(u))$.
		\end{lemma}	
		The proof of this lemma appears in Appendix~\ref{sec:Appendix}. With this result we have
		\[
		\cost_{\ell, \widehat{c}}(\mathcal{T}_\mathcal{P}(r)) - \cost_{\ell',c'}(\mathcal{T}_{\mathcal{P}'}(r)) = \E\left[\E\left[ \cost_{\ell, \widehat{c}}(\mathcal{T}_\mathcal{P}(u)) - \cost_{\ell',c'}(\mathcal{T}_{\mathcal{P}'}(u)) \mid \text{Reach node } u \right] \right]\geq 0
		\]
		
		\item Now, starting from policy $\mathcal{P}'$ and cost-labeled tree $\mathcal{T}_{\mathcal{P}'}$ with labels $\ell'$ in the nodes and $c'$ in the arcs, repeat step 2 until every bin exceeds the risk budget $\gamma/C$ at most once.
	\end{enumerate}
	
	
	With the previous method, we construct a policy, which we still call $\mathcal{P}'$ for simplicity, in addition to its policy tree $\mathcal{T}_{\mathcal{P}'}$ and labels $\ell'$ and $c'$. This policy exceeds each bin's risk budget at most once. Note that labels $c'$ take values in $\{ 0, C+\delta, C+2\delta \}$; we modify the label of arc $(a,b)$ to $\min \{ c_{(a,b)}', C+\delta \}$ thus labels take values only in $\{ 0,C+\delta\}$. This does not increases the cost of the policy tree.
		
For the second phase, we further modify $\mathcal{P}'$: If the policy tries to exceed some bin's risk budget, we open a new bin for that item, unless there is only one item packed in the bin, in which case we move to modify another bin. Using the notation of the first phase, this means that whenever the policy reaches node $u_k$ in some branch starting at $u$, instead of packing the item in node $u_k$ into bin $j$, it opens a new bin $j''$ for it. We call this new policy $\mathcal{P}''$. We modify the cost labels accordingly to accommodate this new cost. We label all nodes $u_k$ in the subtree $\mathcal{T}_{\mathcal{P}'}(u)$ that are not leaves (i.e.\ the risk goes beyond $\gamma/C$ at $u_k$) with $+1$ (the cost to open a new bin). All arcs going out of paths $u_1,\ldots,u_k$ are relabeled from $C+\delta$ to $C$. Formally, we define the new labels
	\[
	c_{(a,b)}'' = \begin{cases}
		\frac{C}{C+\delta} c_{(a,b)}' & a=u_m \text{ for some branch starting at }u\\
		c_{(a,b)}' & \text{otherwise}
	\end{cases}
	\]
	and for nodes
	\[
	\ell_a'' = \begin{cases}
		1 & a = u_k \text{ and }a \text{ is not a leaf} \\
		\ell_a' & \text{otherwise}
	\end{cases}.
	\]
	Using the same argument as in Lemma~\ref{lem:cost_reduce_tree_analysis}, we can show that
	\[
	\cost_{\ell'',c''}(\mathcal{T}_{\mathcal{P}''} (u)) \leq \cost_{c',\ell'}(\mathcal{T}_{\mathcal{P}'}(u)).
	\]
	We repeat this procedure as many times as necessary, and we obtain a policy $\mathcal{P}''$ that satisfies
	\[
	\cost(\mathcal{P}'') \leq \cost(\mathcal{P}') \leq \left( 1 + 2 \frac{\delta}{C}   \right) \cost(\mathcal{P}).
	\]

For ease of reading, we present the proof in an iterative manner; the proof's steps can be followed to obtain the result for finite and countably infinite policy trees. We now sketch how to generalize the proof for uncountable policy trees, focusing on the first phase of the proof. We note that the proof of Lemma~\ref{lem:cost_reduce_tree_analysis} is general and does not require any iterative argument.
Recursively, 
	\[
	\cost_{\ell,c}(\mathcal{T}_{\mathcal{P}} (r)) = \E_{X_1,\ldots,X_{i-1}}\left[ \sum_{k=1}^{i-1} \ell_{U_k} + \sum_{k=1}^{i-1} c_{(U_k, U_{k+1})} + \cost_{\ell, c}(\mathcal{T}_{\mathcal{P}}(U_{i}) ) \right]
	\]
	for any $i=1,\ldots,n$, where $U_i$ is the (random) node at level $i$. Starting at the root, we apply Lemma~\ref{lem:cost_reduce_tree_analysis} to all nodes at level $i$ where a bin is opened. Therefore, we have
	\[
	\cost_{\ell,\widehat{c}}(\mathcal{T}_{\mathcal{P}} (u) ) \geq \cost_{\ell',c'}(\mathcal{T}_{\mathcal{P}'}(u))
	\]
	for all nodes $u$ at level $i$. Using the previous equation, 
	\[
	\cost_{\ell,\widehat{c}}(\mathcal{T}_{\mathcal{P}}(r)) - \cost_{\ell',c'}(\mathcal{T}_{\mathcal{P}}(r)) = \E_{X_1,\ldots,X_{i-1}}\left[  \cost_{\ell,\widehat{c}}(\mathcal{T}_{\mathcal{P}} (U_i) ) - \cost_{\ell',c'}(\mathcal{T}_{\mathcal{P}'}(U_i))\right] \geq 0.
	\]
	Doing this for all levels $i=1,\ldots,n$, we conclude the first phase of the proof. The second phase is completely analogous and omitted for brevity.
	\hfill\halmos
\endproof

Theorem~\ref{thm:key_1} is a general result that does not depend on item size distributions. In the following, we use it to analyze the performance of Budgeted Greedy.

\paragraph{I.I.D.\ Input} When the input is an i.i.d.\ sequence of nonnegative random variables, Budgeted Greedy induces a policy tree that packs one bin at a time: When a bin is opened, the policy never again uses previously opened bins. This simple fact is crucial in the proof of our next result. The next lemma shows that among all budgeted policies, Budgeted Greedy opens the minimum expected number of bins when the input is an i.i.d.\ sequence of random variables. Intuitively, if we ignore the penalty paid by overflowing bins and all items are i.i.d., the optimal way to minimize the expected number of opened bins is by packing as many items as possible in each bin, as long as the risk budget is satisfied. This can of course be done sequentially, one bin at a time, which is what Budgeted Greedy does.

\begin{lemma}\label{lem:key_lem_iid}
	Suppose $X_1,\dotsc,X_n$ are nonnegative i.i.d.\ random variables. Then,
	\[\E[N_\ALG] = \! \min_{\substack{\mathcal{P} \text{ budgeted with}\\\text{risk budget } \gamma/C} }\! \E[N_\mathcal{P}].\]
	That is, among all risk-budgeted policies with budget $\gamma/C$, Budgeted Greedy (Algorithm~\ref{alg:BA}) opens the minimum expected number of bins.
\end{lemma}
\proof{Proof.}
	Consider any policy $\mathcal{P}$ for packing items such that the risk budget of each bin $\gamma/C$ is never surpassed. Consider its tree representation $\mathcal{T}_\mathcal{P}$. We modify the policy tree so only one bin is utilized at a time. For this, we exhibit a sequence of operations ensuring that, whenever bin $j$ is opened, bins $1,\ldots,j-1$ are never utilized again. In the tree, this is equivalent to saying that any branch starting from the root directed to any leaf has labels $1\to j_1, 2\to j_2, \ldots, n\to j_n$ where $1=j_1\leq j_2 \leq \cdots \leq j_n$, where we recall that $i\to j$ means the policy packs item $i$ into bin $j$.
\begin{claim}
Let $j = 1, \dotsc, n $ be any bin opened by the policy. Suppose that node $u$ in level $k$ is labeled $k \to j'$, where $ j' \neq j $. Furthermore, suppose that at node $u$, bin $j$ is open, its usage does not exceed 1, and its risk budget can accommodate $k$. Then $u$ can be relabeled $k\to j$ without increasing the expected number of bins opened by the policy.
\end{claim}
	
	
	Before proving this claim, we show how to use it to conclude the result. Starting at the root $r$ of the policy tree $\mathcal{T}_\mathcal{P}$, find the closest node $u$ to the root where we have the label $k\to j$, $j\neq 1$, but the usage of bin $1$ is no more than $ 1$ and its risk budget can accommodate $k$. Use the claim to relabel this node $k\to 1$ without increasing the expected number of open bins. Repeat this process until there are no nodes $u$ in this category. After this process has been finished, all branches starting at the root have the form $1\to 1, 2\to 1, \ldots, i \to 1, i+1\to 2, \ldots$ and from $i+1$ onward, bin $1$ is overflowed or does not have enough risk budget to receive any additional item. We repeat this process with bin $2, 3,\ldots$. After this process has been carried out, the resulting policy is the one induced by Budgeted Greedy.

	We prove the claim by backward induction on the level of node $u$ in the policy tree. Fix an opened bin $j=1,\ldots,n$ and pick any node $u$ at level $n$ with label $n\to j'$, $j' \neq j$, and suppose bin $j$ is open and satisfies the hypothesis -- its usage is one or less, and it has enough risk budget to receive item $X_n$. Re-labeling this node $n\to j$ does not worsen the number of bins opened since the cost of bin $j$ has already been paid at some previous node. Recall that we are only taking into account the cost paid by opening bins and not the cost of breaking bins.
	
	Suppose the result holds for all levels $k+1,k+2,\ldots,n$. Pick a node $u$ at level $k$ with label $k\to j'$, $j'\neq j$, and such that bin $j$ is open and satisfies the hypothesis. If all its children are labeled $k+1\to j$ then relabel all its children with $k+1 \to j'$ and relabel $u$ with $k\to j$. The cost remains the same after this operation since the distribution of $X_k$ is the same as $X_{k+1}$. Now, suppose that some child of $u$, say $v$, is labeled $k\to m$ with $m\neq j$. Since at node $u$ bin $j$ still has usage not exceeding one and sufficient risk budget left, at node $v$ bin $j$ still satisfies this condition. Therefore, by induction, we can relabel $v$ with $k+1 \to j$ without increasing the cost. We can repeat this for any children of $u$ until all of its children have been labeled $k+1\to j$. We conclude by swapping the label of $u$ with the label of its children as in the previous case.\hfill\halmos
\endproof

\begin{theorem}
	For $\gamma = \sqrt{2}$ and i.i.d.\ nonnegative random variables $X_1,\ldots,X_n$, we have
	\[
	\cost(\ALG)\leq (3 + 2 \sqrt{2}) \cost(\OPT).
	\]
\end{theorem}
\proof{Proof.}
If $\Prob(X_1 > 1) > \sqrt{2}/C$, $\cost(\ALG) = n(1+C \Prob(X_1> 1))$. On the other hand, $\cost(\OPT)\geq n C \Prob(X_1 > 1)$ since each item incurs at least this expected cost. Therefore, $\cost(\ALG) \leq 2 \cost(\OPT)$.
	
	Now assume $\Prob(X_1 > 1) \leq \sqrt{2}/C$. Using Theorem~\ref{thm:key_1} and Lemma~\ref{lem:key_lem_iid}, we have
	\[
	\cost(\ALG) \leq (1+\gamma) \E[N_\mathcal{P}] \leq (1+\gamma) \cost(\OPT^*)	\leq (1+\gamma) \left( 1 + \frac{2}{\gamma} \right)\cost(\OPT),
	\]
	where $\OPT^*$ is the budgeted version of $\OPT$. The expression $(1+\gamma)( 1 + 2/\gamma )$ is minimized at $\gamma=\sqrt{2}$, which gives the desired result. \hfill\halmos
\endproof

\section{Exponential Random Variables}\label{sec:expo_analysis}

In this section, we show that Budgeted Greedy incurs an expected cost at most $\mathcal{O}(\log C)$ times the optimal expected cost when the item sizes are exponentially distributed. That is, for any $X_i$ in the input sequence,
\begin{align}
\Prob(X_i > x) = e^{-\lambda_i x}, \label{eq:exponential}
\end{align}
for any $x\geq 0$, where $\lambda_i > 0$ is the rate. Recall that $\E[X_i]=1/\lambda_i$.

The proof is divided into two parts: First, similarly to deterministic bin packing, we show that $\E\left[ \sum_{i=1}^n \min\{ X_i, 1 \}  \right]$ is a lower bound for $\cost(\mathcal{P})$, for any policy $\mathcal{P}$. In the next step, we show that the probability that Algorithm~\ref{alg:BA} opens bin $k\geq 2$ is related to the amount of mass packed into bin $k-1$. Roughly speaking, we show that the probability that Algorithm~\ref{alg:BA} opens bin $k\geq 2$ is at most $\mathcal{O}(\log C)$ times the expected mass packed into bin $k-1$. Moreover, in Subsection~\ref{subsec:exponential_large_rate} we show that, when the rates governing the item sizes are sufficiently large, $\lambda_i \geq 2 \log C$,  the amount of mass packed into bin $k-1$ is at least a constant, thereby improving the algorithm's approximation factor to a constant. 
Finally, we show that our analysis of Algorithm~\ref{alg:BA} for exponential random variables is almost tight by exhibiting an input sequence that forces Budgeted Greedy to incur a cost $\Omega(\sqrt{\log C})$ times the optimal cost.

\subsection{Arbitrary Exponential Random Variables}

Here we show that $\cost(\ALG) \leq \mathcal{O}(\log C) \cost(\OPT)$ when the input is an arbitrary sequence of exponential random variables. Using Proposition~\ref{prop:breakbinpenalty}, we can assume $\Prob(X_i > 1)\leq 1/C$ for all $i=1,\ldots,n$ at the expense of an extra multiplicative loss of $2$ in the cost incurred. This assumption translates into $\lambda_i\geq \log C$ for all $i$.

The next result shows that the probability that Algorithm~\ref{alg:BA} opens a bin, besides the first bin, is related to the amount of mass packed in the previous bin.
\begin{proposition}\label{prop:exponential_log_bound}
	Suppose $\gamma= 2$. Then, for any $k\geq 2$, Budgeted Greedy  guarantees
	\[
	\Prob(\ALG \text{ opens bin }B_k) \leq 5 \log C \E\left[ \sum_{i\in B_{k-1}} X_i \wedge 1 \right] + \left( \frac{1}{3} + 8 \sqrt{5}\sqrt{\frac{\log C}{C}}  \right) \Prob(\ALG \text{ opens bin } B_{k-1} ).
	\]
\end{proposition}
\proof{Proof.}
Since the item sizes are continuous random variables, we have $\Prob(\ALG \text{ opens bin } B_k)= \Prob(X(B_k)>0)$. Now, we have
	\[
	\Prob(\ALG \text{ opens bin }B_k) \leq \Prob \left(X(B_{k-1}) > \frac{1}{5\log C} \right) + \Prob\left( X(B_{k-1}) \leq \frac{1}{5\log C}, X(B_k)>0  \right).
	\]
	We bound each term separately. To bound the first term we use Markov's inequality:
	\[
	\Prob\left(X(B_{k-1}) > \frac{1}{5 \log C}\right) = \Prob\left(X(B_{k-1}) \wedge 1 > \frac{1}{5 \log C}\right) \leq 5\log C \E[X(B_{k-1}) \wedge 1].
	\]
	For the second term, we proceed as follows. Let $E$ be the event ``\emph{all items packed in $B_{k-1}$ have rate $\lambda_i \geq 2\log C$}.'' Then
	\begin{align*}
	\Prob\left( X(B_{k-1}) \leq \frac{1}{5\log C}, X(B_k)>0  \right)  \leq& \Prob\left( X(B_{k-1}) \leq \frac{1}{5\log C}, X(B_k)>0 , E  \right) \\
	& + \Prob\left( X(B_{k-1}) \leq \frac{1}{5 \log C} \mid \overline{E}  \right) \Prob(X(B_{k-1})>0),
	\end{align*}
	since $\overline{E}$, the event that \emph{some item in $B_{k-1}$ has rate $\leq 2\log C$}, is contained in the event ``\emph{Algorithm~\ref{alg:BA} opens bin $B_{k-1}$}.''
	
	\begin{claim}\label{claim:claim_1_logC}
		$\Prob\left( X(B_{k-1}) \leq \frac{1}{5 \log C} \mid \overline{E}  \right)\leq 1-e^{-2/5} \leq 1/3$.
	\end{claim}
	
	\proof{Proof.}
		If $X_{i_1},\ldots,X_{i_m}$ are all the \emph{large} items with rates $\lambda_{i_p}\leq 2\log C$, then, the events \[M_{p} =\{ X_{i_p} \text{ is the first large item packed into }B_{k-1} \}\] satisfy $\overline{E}=\bigcup_{p=1}^m M_p$. Then,
		\begin{align*}
		\Prob\left( X(B_{k-1}) \leq \frac{1}{5 \log C} \mid \overline{E}  \right) & = \sum_{p=1}^m \Prob\left( X(B_{k-1}) \leq \frac{1}{5 \log C}  \mid \overline{E} , M_p \right)\Prob(M_p \mid \overline{E}) \\
		& \leq \sum_{p=1}^m \Prob\left( X_{i_p} \leq \frac{1}{5\log C} \right)\Prob(M_p \mid \overline{E}) \\
		& = \sum_{p=1}^m (1- e^{-\lambda_{i_p}/5\log C })\Prob(M_p\mid \overline{E}) \tag{Using~\eqref{eq:exponential}} \\
		& \leq (1-e^{-2/5}) \sum_{p=1}^m \Prob(M_p \mid \overline{E}). \tag{Using $\lambda_{i_p}\geq 2\log C$}
		\end{align*}
		The proof follows because the events $M_p$ are disjoint and form $\overline{E}$. \halmos
	\endproof
	
	\begin{claim}\label{claim:claim_2_logC}
		$\Prob\left( X(B_{k-1}) \leq \frac{1}{5\log C}, X(B_k)>0 , E  \right) \leq  20 \sqrt{\frac{\log C}{C}}\Prob(X(B_{k-1})>0)$.
	\end{claim}
	
	\proof{Proof.}
	In this case, bin $B_{k}$ has been opened even though $B_{k-1}$ still has available space. That means that the element that opens bin $B_k$ surpasses the budget of $B_{k-1}$. From here we obtain,
	\[
	\frac{2}{C} < \risk(B_{k-1}) + \Prob(X_t > 1-X(B_{k-1})) \leq \risk(B_{k-1}) + \frac{e^{1/5}}{C},
	\]
	where $X_t$ is the first item packed into $B_k$ and we use $\lambda_t\geq \log C$ and \eqref{eq:exponential}, thus $\Prob(X_t > 1 - X(B_{k-1}))\leq e^{1/5}/C$. Let $F_\beta$ be the event $\{ \sum_{i\in B_{k-1}} \E[X_i] > \beta  \}$; by Markov's inequality,
	\begin{align*}
	\Prob(F_\beta) & \leq \frac{1}{\beta} \E\left[\sum_{i\in B_{k-1}} \E[X_i]  \right] \\
	& \leq  \frac{C}{C-1} \frac{1}{\beta} \E\left[ \sum_{i\in B_{k-1}} \E[X_i \wedge 1] \right] \tag{$\E[X_i \wedge 1] = (1-e^{-\lambda_i})\E[X_i]$} \\
	& \leq  2\frac{C}{C-1} \frac{1}{\beta} \Prob(X(B_{k-1}) >0). \tag{Proposition~\ref{prop:sizeinbins}}
	\end{align*}
	Thus,
	\begin{align*}
	\Prob\left( X(B_{k-1}) \leq \frac{1}{5\log C}, X(B_k)>0 , E  \right) \leq& \Prob\left( X(B_{k-1})\leq \frac{1}{5\log C}, \risk(B_{k-1}) > \frac{2-e^{1/5}}{C} , E  \right)\\
	 \leq & \frac{2C}{\beta (C-1)} \Prob(X(B_{k-1})>0) \\
	& + \Prob\left(X(B_{k-1})\leq \frac{1}{5\log C},\risk (B_{k-1}) > \frac{2-e^{1/5}}{C}, \overline{F}_\beta , E  \right) \\
	\leq & \frac{2C}{\beta(C-1)} \Prob(X(B_{k-1})>0)\\
	& + \frac{C}{2 - e^{1/5}} \E \left[ \risk(B_{k-1}) \mid X(B_{k-1})\leq \frac{1}{5\log C}, E,\overline{F}_{\beta} \right].
	\end{align*}
	\begin{claim}
		$\E \left[ \risk(B_{k-1}) \mid X(B_{k-1})\leq \frac{1}{5\log C}, E,\overline{F}_{\beta} \right] \leq 10 \beta \frac{\log C}{C^2} \Prob(X(B_{k-1})>0)$.
	\end{claim}
	\proof{Proof.}
	Given $X(B_{k-1})\leq \frac{1}{5\log C}$, the event $E$, the event $\overline{F}_\beta$ \emph{and} the event $X(B_{k-1})>0$, the value $\risk(B_{k-1})=\sum_{i=1}^n \Prob(X_i + S_{k-1}^{i-1} > 1)\mathbf{1}_{\{ i\to k-1 \}}^\ALG \leq \sum_{i=1}^n e^{-\lambda_i(1-1/5\log C)}\mathbf{1}_{\{i\to k-1\}}^\ALG$ can be upper bounded by the non-convex problem:
	\[
	\max_{x_1,\ldots,x_n}\left\{  \sum_{i=1}^n e^{-x_i (1- 1/5\log C )}  : \sum_{i=1}^n 1/x_i \leq \beta, x_i \geq 2\log C , \forall i=1,\ldots,n   \right\} \leq 10 \beta \frac{\log C}{C^2}.
	\]
	The inequality follows because the maximum of a convex function is attained in the boundary of the feasible set. Indeed, the maximum is attained by setting the maximum variables to the bound $2\log C$---which are at most $2\beta \log C$---and the rest of the variables to $+\infty$.\hfill\halmos
	\endproof
	
	Therefore, by upper bounding $\frac{10}{2-e^{1/5}}$ by $20$ and $\frac{C}{C-1}\leq 2$, we obtain
	\[
	\Prob\left( X(B_{k-1}) \leq \frac{1}{5\log C}, X(B_k)>0 , E  \right) \leq  \left(  \frac{4}{\beta} + 20 \beta \frac{\log C}{C} \right) \Prob(X(B_{k-1})>0).
	\]
	The right-hand side is minimized at $\beta = \sqrt{\frac{C}{5\log C}}$.\hfill\halmos
	\endproof
	
	Putting Claims~\ref{claim:claim_1_logC} and~\ref{claim:claim_2_logC} together we obtain
	\[
	\Prob(X(B_k)>0) \leq 5 \log C \E\left[ \sum_{i\in B_{k-1}} X_i \wedge 1 \right] + \left( \frac{1}{3} + 8 \sqrt{5}\sqrt{\frac{\log C}{C}}  \right) \Prob(X(B_{k-1})>0 ). \hfill\halmos
	\]
\endproof

\begin{proposition}
	Let $X_1,\ldots,X_n$ be arbitrary exponential random variables with $\lambda_i \geq \log C$. Algorithm~\ref{alg:BA} with $\gamma=2$ guarantees
	\[
	\cost(\ALG) \leq \frac{15 \log C}{\frac{2}{3} - 8 \sqrt{5} \sqrt{\frac{\log C}{C}}} \cost(\OPT) + \frac{3}{\frac{2}{3} - 8 \sqrt{5} \sqrt{\frac{\log C}{C}}},
	\]
	where $\OPT$ is the optimal policy that knows $n$ and the rates of all the sizes $X_1,\ldots,X_n$ in advance.
\end{proposition}
\proof{Proof.}
	Using Proposition~\ref{prop:exponential_log_bound} we obtain
	\begin{align*}
		\E[N_\ALG] & = 1 + \sum_{k\geq 2} \Prob(X(B_{k} > 0)) \\
		& \leq 1 + 5\log C \sum_{k\geq 2} \E\left[ \sum_{i\in B_{k-1}^\ALG} X_i\wedge 1 \right]  + \left( \frac{1}{3} + 8\sqrt{5}\sqrt{\frac{\log C}{C}}  \right)\sum_{k\geq 2} \Prob(X(B_{k-1})>0) \\
		& \leq 1 + 5 \log C \sum_{k\geq 1} \E\left[ \sum_{i\in B_{k}^\ALG} X_i \wedge 1 \right] + \left( \frac{1}{3} + 8 \sqrt{5} \sqrt{\frac{\log C}{C}} \right) \E[N_\ALG] \\
		& = 1 + 5 \log C \E\left[ \sum_{i} X_i \wedge 1  \right] +  \left( \frac{1}{3} + 8 \sqrt{5} \sqrt{\frac{\log C}{C}} \right) \E[N_\ALG].
	\end{align*}
	For any policy $\mathcal{P}$ we have $\E\left[\sum_{i=1}^n X_i \wedge 1 \right]\leq \cost(\mathcal{P})$, using Proposition~\ref{prop:sizesinbins_cost}. Then,
	\begin{align*}
		\E[N_\ALG] & \leq \frac{5 \log C}{\frac{2}{3} - 8 \sqrt{5} \sqrt{\frac{\log C}{C}}} \E[N_\mathcal{P}] + \frac{1}{\frac{2}{3} - 8 \sqrt{5} \sqrt{\frac{\log C}{C}}}
	\end{align*}
	The conclusion follows from here using $\cost(\ALG)\leq 3 \E[N_\ALG]$ (Corollary~\ref{cor:boundgamma}).\hfill\halmos
\endproof

\subsection{Small Exponential Random Variables}\label{subsec:exponential_large_rate}

We next show that $\cost(\ALG)\leq \mathcal{O}(1)\cost(\OPT)$ whenever the item sizes are independent exponential random variables with rates satisfying $\lambda_i\geq 2\log C$.
In this case, $\E\left[ \sum_i X_i\wedge 1 \right]$ is a better approximation for $\E[N_\ALG]$ than in the general case. The following results shows that we can improve Proposition~\ref{prop:exponential_log_bound} by a logarithmic factor.
\begin{proposition}\label{prop:exponential_constant_bound}
	Let $\gamma=1$. For $ k\geq 2$,
	\[
	\Prob(\ALG \text{ opens bin } {k}) \leq 4\E\left[\sum_{i\in B_{k-1}} X_i \wedge 1\right] + 8\frac{\sqrt{\log C}}{C^{1/4}} \Prob(\ALG \text{ opens bin }k-1).
	\]
\end{proposition}
\proof{Proof.}
	We have
	\begin{align*}
	\Prob(\ALG \text{ opens bin }B_{k+1}) & = \Prob(X(B_{k+1}) > 0 ) \\
	&\leq \Prob( X(B_k) > 1/4 ) + \Prob(X(B_{k+1})>0, X(B_k) \leq 1/4 ) \\
	&\leq \Prob( X(B_k) \wedge 1 > 1/4 ) + \Prob(X(B_{k+1})>0, X(B_k) \leq 1/4 ) \\
	& \leq 4 \E\left[ X(B_k)\wedge 1 \right] + \Prob(X(B_{k+1})>0, X(B_k) \leq 1/4 ).
	\end{align*}
	We only focus on bounding the second term in the rest of the proof. Algorithm~\ref{alg:BA} opens bin $B_{k+1}$ ($X(B_{k+1})>0$) if there are no available bins ($\forall i\leq k$, $X(B_i)\geq 1$) or there is an item that does not fit because of the budget. The first case cannot happen when the event $X(B_k)\leq 1/4$ happens so we are only left with the budget case. In particular, for bin $k$, we open bin $B_{k+1}$ because for some item $X_t$ we have
	\[
	\frac{1}{C} < \risk(B_k) + \Prob( X_t + X(B_{k}) > 1) \leq \risk(B_k) + \Prob(X_t > 3/4) \leq \risk(B_k) + \frac{1}{C^{3/2}}
	\]
	where we used the information from the event $X(B_{k+1})\leq 1/4$. Therefore,
	\begin{align*}
	\Prob(X(B_{k+1})> 0, X(B_k)\leq 1/4) &\leq \Prob( \risk(B_k) > 1/C-1/C^{3/2}, 0<X(B_k)<1/4 )\\
	& \leq \frac{C^{3/2}}{C^{1/2}-1}\E[ \risk(B_k) \mid X(B_k)<1/4, X(B_k)>0 ]\Prob(X(B_k)>0).
	\end{align*}
	Now, as in the previous proof, let $F_\beta = \left\{ \sum_{i\in B_k} \E[X_i] > \beta \right\}$; by Markov's inequality and Proposition~\ref{prop:sizeinbins},
	\begin{align*}
	\Prob\left( F_\beta  \right) & \leq \frac{2C^2}{\beta(C^2 - 1)}\Prob(X(B_k)>0).
	\end{align*}
	\begin{claim}
		$\E[ \risk(B_k) \mathbf{1}_{\{X(B_k)<1/4\}}  \mid X(B_k)<1/4, X(B_k)>0 , \overline{F}_\beta ] \leq \frac{2 \beta \log C}{C^{3/2}}$.
	\end{claim}
	\proof{Proof.}
		Given $ X(B_k)<1/4, X(B_k)>0 , \overline{F}_\beta$, the risk
		\[
		\risk(B_k) = \sum_{i=1}^n \Prob(X_i + S_k^{i-1}> 1)\mathbf{1}_{\{i\to k\}} \leq \sum_{i=1}^n e^{-\lambda_i \cdot\frac{3}{4}} \mathbf{1}_{\{i\to k\}}
		\]
		is bounded by the non-convex problem,
		\[
		\max_{x_1,\ldots,x_n}\left\{ \sum_{i=1}^n e^{-x_i \frac{3}{4}} : \sum_{i=1}^n \frac{1}{x_i} \leq \beta, x_i\geq 2\log C , \forall i=1,\ldots,n \right\} \leq \frac{2 \beta \log C}{C^{3/2}},
		\]
which we bound as before.
	\hfill\halmos
	\endproof
	
	With this claim,
	\begin{align*}
	\E[ \risk(B_k) \mid X(B_k)<1/4, X(B_k)>0 ] & \leq \frac{1}{C}\Prob(F_\beta) + \frac{2 \beta \log C}{C^{3/2}} \Prob(\overline{F}_\beta) \\
	& \leq \frac{2 C}{\beta (C^2 - 1)} + \frac{2\beta \log C}{C^{3/2}},
	\end{align*}
	since $\risk(B_k)\leq 1/C$. Thus,
	\begin{align*}
	\Prob(X(B_{k+1})> 0, X(B_k)\leq 1/4) & \leq \frac{C^{3/2}}{C^{1/2}- 1} \left(   \frac{2C}{\beta(C^2-1)} + \frac{2\beta \log C}{C^{3/2}} \right)\Prob(X(B_k)>0).
	\end{align*}
	Now, optimizing over $\beta$ with $\beta = \frac{C^{5/4}}{\sqrt{C^2 - 1}\sqrt{\log C}}$ we obtain
	\[
	\Prob(X(B_{k+1})> 0, X(B_k)\leq 1/4 ) \leq 8 \frac{\sqrt{\log C}}{C^{1/4}}\Prob(X(B_k)>0). \hfill\halmos
	\]
\endproof

\begin{proposition}
Suppose $ \lambda_i \geq 2 \log C $ for all $ i = 1, \dotsc, n $.
	For $\gamma=1$, Algorithm~\ref{alg:BA} guarantees
	\[
	\cost(\ALG)\leq \frac{8}{1 - 8 \frac{\sqrt{\log C}}{C^{1/4}}} \cost(\OPT) + \frac{2}{1 - 8 \frac{\sqrt{\log C}}{C^{1/4}}},
	\]
	where $\OPT$ is the optimal policy that knows $n$ and all item size rates in advance.
\end{proposition}
\proof{Proof.}
Using Proposition~\ref{prop:exponential_constant_bound} we have
\begin{align*}
	\E[N_\ALG] & = \sum_{k=1}^n \Prob(\ALG \text{ opens bin }k)  \\
	& \leq 1+  \sum_{k=2}^n 4\E\left[\sum_{i\in B_{k-1}} X_i \wedge 1\right] + 8\frac{\sqrt{\log C}}{C^{1/4}} \Prob(\ALG \text{ opens bin }k-1) \\
	& \leq 1 + 4 \E\left[\sum_{i=1}^n X_i \wedge 1 \right] + 8 \frac{\sqrt{\log C}}{C^{1/4}} \E[N_\ALG].
\end{align*}
Using Proposition~\ref{prop:sizesinbins_cost}, 
\[
\E[N_\ALG] \leq \frac{4}{1- 8 \frac{\sqrt{\log C}}{C^{1/4}}} \cost(\mathcal{P}) + \frac{1}{1- 8 \frac{\sqrt{\log C}}{C^{1/4}}}
\]
for any policy $\mathcal{P}$. The result follows by using $\cost(\ALG)\leq 2 \E[N_\ALG]$ and optimizing over $\mathcal{P}$.\hfill\halmos
\endproof

\subsection{A Lower Bound for the Algorithm with Exponential Random Variables}\label{subsec:exp_lower_bound}

In this subsection, we present a hard input of exponential random variables for Budgeted Greedy. The sequence contains two kind of independent exponential random variables, those with rates $\mu = \beta \log C$, $\beta \geq 2$ and those with rates $\lambda = (1+\varepsilon)\log C$, with $\varepsilon \in (0,1)$. This sequence has $n_1$ items with rate $\lambda$ and $n_2=k n_1$ items with rate $\mu$, presented to Algorithm~\ref{alg:BA} as, 
\[
X_{1,1}^\mu \cdots X_{1,k}^\mu X_{1}^{\lambda} X_{1,1}^\mu \cdots X_{2,k}^\mu X_{2}^\lambda\quad  \cdots\quad  X_{n_1,1}^\mu   \cdots  X_{n_1,k}^\mu X_{n_1}^\lambda,
\]
where $X_{i,j}^\mu \sim \exp(\mu)$ and $X_i^\lambda \sim \exp(\lambda)$ for all $i,j$. With the choices of $\beta= 6\frac{n_1 \log C}{\varepsilon}$ and $k= 3\varepsilon \mu = 18 n_1 (\log C)^2$, we show that Budgeted Greedy incurs an expected cost of at least $\frac{1}{2}n_1$. For the same choices of $\beta$ and $k$ and optimizing over the choice of $\varepsilon$, we show that $\cost(\OPT)\leq \mathcal{O}\left(1/\sqrt{\log C}\right) n_1$. This choice of $\varepsilon$ is independent of $n_1$, which allows us to scale the result for any input size.

We prove each bound separately; the main results are stated here.
\begin{proposition}\label{prop:lower_bound_exp_1}
	Let $\varepsilon >0$ and set $\beta = 6\frac{n_1 \log C}{\varepsilon}$ and $k=3\varepsilon \mu$. Then, running Budgeted Greedy with $\gamma=1$ on the input described above yields
	\[
	\cost(\ALG) \geq \frac{1}{2} n_1.
	\]
\end{proposition}

\begin{proposition}\label{prop:upper_bound_opt_cost_1}
	Using the same parameters as in the previous proposition, for any $\varepsilon >0$ such that $\varepsilon \log C \geq 4$, we have
	\[
	\cost(\OPT) \leq 48n_1 \left(  \frac{k}{\beta\log C} + \frac{1}{\varepsilon \log C} \right)=48 n_1 \left(  3\varepsilon + \frac{1}{\varepsilon \log C}  \right).
	\]
\end{proposition}
The result now follows by taking $\varepsilon = \sqrt{\frac{1}{3\log C}}$. The proofs are in Appendix~\ref{sec:Appendix}.

\section{Offline Sequential Adaptive Bin Packing}\label{sec:soft_ptas_offline}

In this section, we move to the offline sequential model, where random variables are known in advance and the packing occurs sequentially in the fixed order $1,\ldots,n$. We present the proof of Theorem~\ref{thm:main_soft_ptas_offline} that guarantees a soft-capacity polynomial time approximation scheme (PTAS) for the offline problem. 

\subsection{Approximation of a Sequential Policy}

Consider $X_1,\ldots,X_n$ independent random variables with bounded support $[0,1+\varepsilon]$. We can reduce the general case to this case by moving all the probability mass of the corresponding random variable in $[1+\varepsilon, \infty)$ to the point $1+\varepsilon$. We aim to show a polynomial time approximation scheme with resource augmentation. In particular, we consider a policy operating on bins with size or capacity $c \geq 1$; a bin overflows if the total size of items packed into it exceeds $c$. We use the notation $\cost_{c}(\mathcal{P}, Z)$ to denote the expected cost incurred by a policy $\mathcal{P}$ packing items $Z=(Z_1,\dotsc,Z_n)$ into bins of capacity $c$.
\begin{theorem}\label{thm:approximation_thm}
	There is a policy that can be computed in $\mathcal{O}\left(\frac{1}{\varepsilon^{10}} n^{2/\varepsilon^5}\right)$ time packing items $X_1,\dotsc,X_n$ sequentially into bins of size $1+6\varepsilon$, and incurring expected cost of at most $ (1+4\varepsilon)\cost_1(\OPT, X)$, where $\OPT$ is an optimal policy with respect to bins of unit size.
\end{theorem}


To prove Theorem~\ref{thm:approximation_thm}, we proceed as follows in the remainder of the section:
\begin{enumerate}[leftmargin=*, itemsep=0em]
\item First, we discretize the input random variables $X_1,\dotsc,X_n$ into random variables $\widehat{X}_1,\ldots,\widehat{X}_n$ with support in $\{ 0,\varepsilon^5,\dotsc, \lceil 2/\varepsilon^5\rceil \varepsilon^5 \}$. This allows us to compute an optimal policy in polynomial time via dynamic programming.
	
\item We then show that for any policy $\mathcal{P}$ for $X_1,\dotsc,X_n$, we can construct a policy $\widehat{\mathcal{P}}$ for $\widehat{X}_1,\ldots,\widehat{X}_n$ such that
	\[
	\cost_{1+4\varepsilon}\left( \widehat{\mathcal{P}}, \widehat{X} \right)  \leq (1+\varepsilon) \cost_1(\mathcal{P}, X).
	\]
	
\item Next, we show how to obtain a policy $\mathcal{P}$ for $X_1,\dotsc,X_n$ from a policy $\widehat{\mathcal{P}}$ for items $\widehat{X}_1,\ldots,\widehat{X}_n$, such that
	\[
	\cost_{1+6\varepsilon} (\mathcal{P}, X) \leq (1+\varepsilon) \cost_{1+4\varepsilon}\left( \widehat{\mathcal{P}}, \widehat{X} \right) .
	\]
	
\item Finally, we show that we can compute the optimal policy $\widehat{\mathcal{P}}$ for discretized items $\widehat{X}_1,\dotsc,\widehat{X}_n$ in $\mathcal{O}\left( \frac{1}{\varepsilon^{10}}n^{2/\varepsilon}  \right)$ time. The policy $\mathcal{P}$ follows immediately from here.
\end{enumerate}

\subsection{Discretization Process}

We perform the discretization in two steps, similarly to the discretization in~\cite{li2013stochastic}. In the first step, we discretize the small outcomes of $X_1,\dotsc,X_n$, meaning that values not exceeding $ \varepsilon^4$ now behave as a scaled Bernoulli random variable with scaling factor $\varepsilon^4$, and the appropriate success probability such that this discretization preserves the expectation of the original random variable. In the second step, we discretize the large outcomes by rounding up all values to multiples of $\varepsilon^5$.

This discretization allows us to construct a state space based on the number of bins at level $k \varepsilon^5$, $k=1,\dotsc, \lceil 2/\varepsilon^5 \rceil$. The number of states is roughly $\mathcal{O}(n^{2/\varepsilon^5})$, which is polynomial in $n$.

\paragraph{Step 1 of discretization}

Let $q_i = \E[X_i \mid X_i \leq \varepsilon^4]$. Then, the first discretization is
\[
X_i' = \begin{cases}
0 & \text{if } X_i \leq \varepsilon^4, \text{ w.p. } 1- q_i/\varepsilon^4 \\
\varepsilon^4 & \text{if }X_i \leq \varepsilon^4, \text{ w.p. } q_i/\varepsilon^4 \\
X_i & \text{ if } X_i > \varepsilon^4 .
\end{cases}
\]
Note that we have $ \lvert X_i-X_i' \rvert \leq \varepsilon^4 $ almost surely, $\E[X_i'\mid X_i' \leq \varepsilon^4] = \E[X_i \mid X_i \leq \varepsilon^4]$ and $\E[X_i]=\E[X_i']$.
%
%
%
%

\paragraph{Step 2 of discretization}

Now consider
\[
\widehat{X}_i = \mathbf{1}_{\{ X_i' \leq \varepsilon^4 \}} X_i' + \mathbf{1}_{\{ X_i' > \varepsilon^4 \}} \lceil X_i'/\varepsilon^5 \rceil \varepsilon^5.
\]
Clearly, $X_i' \leq \widehat{X}_i$, since the large outcomes are rounded up. Moreover, if $b>\varepsilon^4$, then $\lceil b/\varepsilon^5  \rceil\varepsilon^5 \leq \left( b/\varepsilon^5 + 1 \right) \varepsilon^5 = b + \varepsilon^5 \leq (1+\varepsilon) b $. Hence, $X_i' \leq \widehat{X}_i \leq (1+\varepsilon) X_i'$.

%
%
%

\subsection{From Regular Policy to Discretized Policy}

In this subsection we show the following result:

\begin{theorem}\label{thm:regular_to_discrete}
	For any policy $\mathcal{P}$ that sequentially packs items $X_1,\dotsc,X_n$ into bins of unit size, there exists a policy $\widehat{\mathcal{P}}$ packing items $\widehat{X}_1,\dotsc,\widehat{X}_n$ into bins of size $1+4\varepsilon$ such that
	\[
	\cost_{1+4\varepsilon}\left( \widehat{\mathcal{P}}, \widehat{X} \right) \leq (1+\varepsilon) \cost_1 \left( \mathcal{P}, X\right).
	\]
\end{theorem}

To prove the theorem, we first introduce an intermediate policy $\mathcal{P}'$ that packs items $X_1',\dotsc,X_n'$ and satisfying
\[
\cost_{1+2\varepsilon}(\mathcal{P}', X' ) \leq (1+\varepsilon) \cost_1(\mathcal{P}, X).
\]
Policy $\widehat{\mathcal{P}}$ is obtained from policy $\mathcal{P}'$ by adding additional capacity to the bins.

We assume that $\mathcal{P}$ is a deterministic function of the capacity of the bins and the current element to be packed. Let us construct a policy $\mathcal{P}'$ for items $X_1',\ldots,X_n'$ with bin capacity $1+2\varepsilon$ that simulates and follows policy $\mathcal{P}$ in the following way. Upon arrival of item $X_i'$, policy $\mathcal{P}'$ does what $\mathcal{P}$ would have done at this point in time to item $X_i$. We couple $X_i$ and $X_i'$, so $X_i=X_i'$ if $X_i> \varepsilon^4$ and otherwise we have the Bernoulli behavior in $X_i'$. We pass the outcome of $X_i'$ to $\mathcal{P}'$ and the outcome of $X_i$ to $\mathcal{P}$. (Strictly speaking, $\mathcal{P}'$ receives the outcome of $X_i'$ and from it, the policy samples $X_i$ coupled with $X_i'$ and passes this outcome to $\mathcal{P}$.) For each bin $B_j$ that policy $\mathcal{P}$ opens, policy $\mathcal{P}'$ opens a bin $B_{j,1}'$ and packs items in $B_{j,1}'$ as policy $\mathcal{P}$ would do in bin $B_j$ as long as $|X(B_{j,1}') - X'(B_{j,1}')|\leq \varepsilon$ holds. If this difference is violated, policy $\mathcal{P}'$ opens a new bin $B_{j,2}'$ and continues following $\mathcal{P}$ as long as $|X(B_{j,2}') - X'(B_{j,2}')|\leq \varepsilon$ holds, and so on.

Notice that $\mathcal{P}'$ is undefined if some $B_{j,k}'$ breaks but $B_{j}$ is not broken by policy $\mathcal{P}$. Fortunately, this event cannot occur, as the following proposition guarantees. 
\begin{proposition}
	If $\mathcal{P}$ overflows $B_{j,k}'$ for some $k$, then $\mathcal{P}$ must have overflown $B_j$. In particular, at most one of the $B_{j,k}'$ is overflowed by $\mathcal{P}'$.
\end{proposition}

\proof{Proof.}
	Let $B_{j,k,t}'$ be the items packed into bin $B_{j,k}'$ up to time $t$. Suppose that $X'(B_{j,k,t}') > 1+ 2\varepsilon$ ($B_{j,k}'$ is overflowed at time $t$). Notice that $X_t' \leq 1+\varepsilon$, therefore $B_{j,k}'$ was opened before $t$ and
	\[
	|X'(B_{j,k,t-1}') - X(B_{j,k,t-1}')| \leq \varepsilon
	\]
	otherwise $\mathcal{P}'$ would not have tried to pack $X_t$ into $B_{j,k}'$. Now, since $|X_t'-X_t|\leq \varepsilon^4$ we have
	\begin{align*}
		X(B_j) &\geq X(B_{j,k,t}') \\
		&= X(B_{j,k,t-1}') + X_t \geq X'(B_{j,k,t-1}) -\varepsilon + X_t'-\varepsilon^4 \\
		& = X'(B_{j,k,t}') -\varepsilon - \varepsilon^4\\
		& > 1 + \varepsilon - \varepsilon^4 > 1.
	\end{align*}
	For the second part, we notice that once $B_{j,k}'$ is overflowed, then $B_j$ is overflowed as well and so $\mathcal{P}$ does not pack any item in $B_j$. Then, after $B_{j,k}'$ no more bins $B_{j,k+1}',\ldots$ are open.\hfill\halmos
\endproof

Let $O_{\mathcal{P}',j}$ be the number of bins $B_{j,k}'$ that policy $\mathcal{P}'$ breaks; we just showed that $O_{\mathcal{P'},j}\leq \mathbf{1}_{\{ X(B_j) > 1 \}}^{\mathcal{P}}$. Then, $O_{\mathcal{P}'}=\sum_{j=1}^n O_{\mathcal{P}',j}$, the number of bins overflowed by $\mathcal{P}'$, satisfies the following equality.
\begin{proposition}	$\E[O_{\mathcal{P}'}] \leq \E[O_\mathcal{P}].$
\end{proposition}
\proof{Proof.}
	By the previous proposition, at most one of the $B_{j,1}',\dotsc,B_{j,n}'$ breaks, and when it does then $B_j$ must have been broken as well. \hfill\halmos
\endproof

Next, we show that the number of bins opened by $\mathcal{P}'$ is not much larger than the number of bins opened by $\mathcal{P}$. Let $N_{\mathcal{P}',j}$ be the number of bins $B_{j,1}',B_{j,2}',\ldots$ that policy $\mathcal{P}'$ uses, i.e.\ the number of copies of bin $B_j$ used by policy $\mathcal{P}$. Let $N_{\mathcal{P}'}$ be the number of bins opened by policy $\mathcal{P}'$, $N_{\mathcal{P}'} = \sum_{j=1}^n N_{\mathcal{P}',j}$.

Consider the family of events $\mathcal{E}_{j,k}' = \{ |X(B_{j,k}') - X'(B_{j,k}')| > \varepsilon  \}$ for $k\geq 1$ and for $k=0$ define $\mathcal{E}_{j,0}'= \{ \mathcal{P}' \text{ opens bin } B_{j,1}' \}  = \{ \mathcal{P} \text{ opens } B_j \}$. Notice then, for $\ell\geq 1$,
\begin{align}
	\left\{N_{\mathcal{P}',j} \geq \ell\right\} \subseteq \mathcal{E}_{j,0}'\cap \mathcal{E}_{j,1}'\cap  \cdots \cap \mathcal{E}_{j,\ell-1}'. \label{inc:policy_bins_containment}
\end{align}

\begin{proposition}\label{prop:key_lem_bound_on_deviation_discr_1}
	For any $k\geq 1$,
	\[
	\Prob(\mathcal{E}_{j,k}' \mid \mathcal{E}_{j,k-1}',\ldots,\mathcal{E}_{j,1}',\mathcal{E}_{j,0}') \leq 6\varepsilon^2 \Prob(\mathcal{P}' \text{ opens }B_{j,k}'\mid \mathcal{E}_{j,k-1}', \ldots, \mathcal{E}_{j,1}',\mathcal{E}_{j,0}').
	\]
\end{proposition}

\proof{Proof.}
	Using Chebychev's inequality,
	\begin{align*}
	\Prob(\mathcal{E}_{j,k}' \mid \mathcal{E}_{j,k-1}',\ldots,\mathcal{E}_{j,1}',\mathcal{E}_{j,0}') & \leq \frac{1}{\varepsilon^2} \E\left[ \left( X(B_{j,k}') - X'(B_{j,k}') \right)^2 \mid \mathcal{E}_{j,k-1}',\ldots,\mathcal{E}_{j,1}',\mathcal{E}_{j,0}'   \right]\\
	&= \frac{1}{\varepsilon^2} \E\left[ \left( \sum_{i=1}^n (X_i - X_i')\mathbf{1}_{\{ i \to (j,k) \}}^{\mathcal{P}'}  \right)^2 \mid \mathcal{E}_{j,k-1}',\ldots,\mathcal{E}_{j,1}',\mathcal{E}_{j,0}'   \right] \\
	& = \frac{1}{\varepsilon^2} \sum_{i=1}^n \E\left[ (X_i-X_i')^2 \mathbf{1}_{\{ i \to (j,k) \}}^{\mathcal{P}'} \mid \mathcal{E}_{j,k-1}',\ldots, \mathcal{E}_{j,1}',\mathcal{E}_{j,0}'  \right] \\
	& \quad + \frac{2}{\varepsilon^2} \sum_{i<\ell} \E\left[ (X_i-X_i')(X_\ell-X_{\ell}') \mathbf{1}_{\{ i \to (j,k), \ell\to (j,k) \}}^{\mathcal{P}'} \mid \mathcal{E}_{j,k-1}',\ldots, \mathcal{E}_{j,1}', \mathcal{E}_{j,0}'  \right].
	\end{align*}
	\begin{claim}
		For $i < \ell$, $ \E\left[ (X_i-X_i')(X_\ell-X_{\ell}') \mathbf{1}_{\{ i \to (j,k), \ell\to (j,k) \}}^{\mathcal{P}'} \mid \mathcal{E}_{j,k-1}',\ldots, \mathcal{E}_{j,1}'  \right]=0$.
	\end{claim}
	\proof{Proof.}
		If $\Prob(i \to (j,k), \ell \to (j,k) \mid \mathcal{E}_{j,k-1}',\ldots, \mathcal{E}_{j,1}', \mathcal{E}_{j,0}') =0$ the result is clearly true, while in the opposite case
		\begin{align*}
		\E\left[ (X_i-X_i')(X_\ell-X_{\ell}') \right. &  \left. \mathbf{1}_{\{ i \to (j,k), \ell\to (j,k) \}}^{\mathcal{P}'} \mid \mathcal{E}_{j,k-1}',\ldots, \mathcal{E}_{j,1}', \mathcal{E}_{j,0}'  \right] \\
		& = \E\left[ (X_i-X_i')(X_\ell-X_{\ell}') \mid i \to (j,k), \ell \to (j,k), \mathcal{E}_{j,k-1}',\ldots, \mathcal{E}_{j,1}', \mathcal{E}_{j,0}'  \right] \\
		& \quad \times \Prob(i \to (j,k), \ell \to (j,k) \mid \mathcal{E}_{j,k-1}',\ldots, \mathcal{E}_{j,1}', \mathcal{E}_{j,0}')\\
		& = \E[X_\ell-X_{\ell}']\E\left[ X_i-X_i' \mid i \to (j,k), \ell \to (j,k), \mathcal{E}_{j,k-1}',\ldots, \mathcal{E}_{j,1}', \mathcal{E}_{j,0}'  \right] \\
		& \quad \times \Prob(i \to (j,k), \ell \to (j,k) \mid \mathcal{E}_{j,k-1}',\ldots, \mathcal{E}_{j,1}', \mathcal{E}_{j,0}')\\
		& = \E\left[ \E[X_\ell - X_\ell'] (X_i-X_i') \mathbf{1}_{\{ i\to (j,k), \ell\to (j,k) \}}^{\mathcal{P}'} \mid \mathcal{E}_{j,k-1}',\ldots, \mathcal{E}_{j,1}',\mathcal{E}_{j,0}'\right]\\
		& = 0,
		\end{align*}
		the last result since $\E[X_\ell]= \E[X_\ell']$. Note that from the second to the third equality, we utilized the fact that given that $\ell$ is packed into $B_{j,k}'$, the outcome of $X_\ell - X_\ell'$ is independent of previous $\mathcal{E}_{j,\ell}'$, $\ell < k$.\hfill\halmos
	\endproof

	\begin{claim}
		For any $i$,
		\[\E\left[ (X_i-X_i')^2 \mathbf{1}_{\{ i \to (j,k) \}}^{\mathcal{P}'} \mid \mathcal{E}_{j,k-1}',\ldots, \mathcal{E}_{j,1}' ,\mathcal{E}_{j,0}' \right] \leq 2\varepsilon^4 \E\left[\E[X_i] \mathbf{1}_{\{ i\to (j,k) \}}^{\mathcal{P}'} \mid \mathcal{E}_{j,k-1}',\ldots,\mathcal{E}_{j,1}', \mathcal{E}_{j,0}'\right].
		\]
	\end{claim}
	\proof{Proof.}
		We have $|X_i-X_i'|\leq \varepsilon^4$, so
		\begin{align*}
		\E\left[ (X_i-X_i')^2 \mathbf{1}_{\{ i \to (j,k) \}}^{\mathcal{P}'} \right.&\left.\mid \mathcal{E}_{j,k-1}',\ldots, \mathcal{E}_{j,1}', \mathcal{E}_{j,0}'  \right] \\
		& \leq \varepsilon^4 \E\left[ |X_i-X_i'| \mathbf{1}_{\{ i \to (j,k) \}}^{\mathcal{P}'} \mid \mathcal{E}_{j,k-1}',\ldots, \mathcal{E}_{j,1}' ,\mathcal{E}_{j,0}' \right] \\
		& \leq \varepsilon^4 \E\left[ (X_i + X_i')\mathbf{1}_{\{i\to (j,k)\}}^{\mathcal{P}'} \mid \mathcal{E}_{j,k-1}', \ldots, \mathcal{E}_{j,1}', \mathcal{E}_{j,0}'\right].
		\end{align*}
		The sizes of $X_i$ and $X_i'$ are independent of the policy $\mathcal{P}'$ packing $i$ into bin $B_{j,k}'$. Furthermore, if $X_i'$ is packed into bin $B_{j,k}'$, its size does not depend on previous events $\mathcal{E}_{j,\ell}'$, $\ell <k$. Therefore,
		\begin{align*}
		\E\left[ X_i' \mathbf{1}_{\{i\to (j,k)\}}^{\mathcal{P}'} \mid \mathcal{E}_{j,k-1}',\ldots,\mathcal{E}_{j,1}', \mathcal{E}_{j,0}'  \right] & = \E[X_i'] \Prob(\mathcal{P}' \text{ packs }i \text{ into }(j,k)\mid \mathcal{E}_{j,k-1}',\ldots,\mathcal{E}_{j,1}', \mathcal{E}_{j,0}' ) \\
		&= \E\left[\E[X_i] \mathbf{1}_{\{ i\to (j,k) \}}^{\mathcal{P}'} \mid \mathcal{E}_{j,k-1}',\ldots,\mathcal{E}_{j,1}', \mathcal{E}_{j,0}'\right]
		\end{align*}
		since $\E[X_i']=\E[X_i]$. Similarly,
		\begin{align*}
		\E\left[ X_i \mathbf{1}_{\{i\to (j,k)\}}^{\mathcal{P}'} \mid \mathcal{E}_{j,k-1}',\ldots,\mathcal{E}_{j,1}' , \mathcal{E}_{j,0}' \right] & = \E\left[\E[X_i] \mathbf{1}_{\{ i\to (j,k) \}}^{\mathcal{P}'} \mid \mathcal{E}_{j,k-1}',\ldots,\mathcal{E}_{j,1}', \mathcal{E}_{j,0}'\right]. \hfill\halmos
		\end{align*}
\endproof
		Putting these two claims together in the previous inequality gives us
		\begin{align*}
		\Prob(\mathcal{E}_{j,k}' \mid \mathcal{E}_{j,k-1}',\ldots,\mathcal{E}_{j,1}',\mathcal{E}_{j,0}') & \leq \frac{1}{\varepsilon^2} \sum_{i=1}^n2\varepsilon^4 \E\left[\E[X_i] \mathbf{1}_{\{ i\to (j,k) \}}^{\mathcal{P}'} \mid \mathcal{E}_{j,k-1}',\ldots,\mathcal{E}_{j,1}', \mathcal{E}_{j,0}'\right] \\
		& \leq 2\varepsilon^2 \E\left[ \sum_{i\in B_{j,k}'} \E[X_i] \mid \mathcal{E}_{j,k-1}', \ldots, \mathcal{E}_{j,1}' ,\mathcal{E}_{j,0}' \right] \\
		& \leq 6 \varepsilon^2 \Prob(\mathcal{P}' \text{ opens } B_{j,k}'\mid \mathcal{E}_{j,k-1}', \ldots, \mathcal{E}_{j,1}', \mathcal{E}_{j,0}').
		\end{align*}
		In the last inequality, we used Proposition~\ref{prop:sizeinbins}, $X_i \leq 1+\varepsilon$ for all $i$, and the bins $B_{j,k}'$ having capacity $1+2\varepsilon$.\hfill\halmos
	\endproof

Recall that $N_{\mathcal{P}',j}$ is the number of bins $B_{j,1}',\ldots$ that policy $\mathcal{P}'$ uses. We have,
\begin{proposition}\label{prop:key_bounded_bins_1}
	For any $j=1,\ldots,n$,
	\[
	\E [N_{\mathcal{P}',j}] \leq (1+\varepsilon)\Prob(\mathcal{P} \text{ opens }B_j).
	\]
\end{proposition}
\proof{Proof.}
	Using the inclusion~\eqref{inc:policy_bins_containment} and the previous proposition,
	\begin{align*}
	\Prob(N_{\mathcal{P}',j}\geq \ell) & \leq \Prob(\mathcal{E}_{j,\ell-1}',\ldots, \mathcal{E}_{j,1}',\mathcal{E}_{j,0}') \\
	& = \Prob(\mathcal{E}_{j,\ell-1}'\mid \mathcal{E}_{j,\ell-2}',\ldots,\mathcal{E}_{j,1}',\mathcal{E}_{j,0}') \cdots \Prob(\mathcal{E}_{j,1}'\mid \mathcal{E}_{j,0}') \Prob(\mathcal{E}_{j,0}') \\
	& \leq (6\varepsilon^2)^{\ell-1} \Prob(\mathcal{P} \text{ opens }B_j).
	\end{align*}
	Thus,
	\begin{align*}
	\E[N_{\mathcal{P}',j}] & = \sum_{\ell\geq 1}  \Prob(N_{\mathcal{P}',j} \geq \ell) \\
	& \leq \sum_{\ell\geq 1}  (6\varepsilon^2)^{\ell-1} \Prob(\mathcal{P} \text{ opens } B_j)\\
	& = \frac{1}{(1-6\varepsilon^2)} \Prob(\mathcal{P}\text{ opens }B_j) \\
	& \leq (1+\varepsilon)\Prob(\mathcal{P}\text{ opens } B_j).
	\end{align*}
	For the last inequality we require $\varepsilon \leq \frac{1}{\sqrt{6}}(\sqrt{15}- 3) \approx 0.1454$.\hfill\halmos
\endproof

\begin{corollary}\label{prop:key_bounded_bins_2}
	$\E[N_{\mathcal{P}'}] \leq (1+\varepsilon)\E[{N}_\mathcal{P}]$.
\end{corollary}


\begin{lemma}
	$ \cost_{1+2\varepsilon}(\mathcal{P}',X') \leq (1+\varepsilon) \cost_1(\mathcal{P},X).$
\end{lemma}
\proof{Proof.}
	This follows from $\cost_{1+2\varepsilon}(\mathcal{P}', X') = \E[N_{\mathcal{P}'}]  + C\E[O_{\mathcal{P}'}]$ and the previous results. \hfill\halmos
\endproof

%
\proof{Proof of Theorem~\ref{thm:regular_to_discrete}.}
	Let $\widehat{\mathcal{P}}$ be the policy constructed from $\mathcal{P}'$ in the following manner. We simulate policy $\mathcal{P}'$ in parallel. To pack item $\widehat{X}_i$, $\widehat{\mathcal{P}}$ imitates what $\mathcal{P}'$ does to item $X_i'$. Random variables $\widehat{X}_i$ and $X_i'$ (and also $X_i$) are assumed to be coupled in the standard manner. The outcome of $\widehat{X}_i$ goes to $\widehat{\mathcal{P}}$ and the outcome of $X_i'$ goes to $\mathcal{P}'$.
	
We denote by $\widehat{B}_{j,k}$ the bins opened by $\widehat{\mathcal{P}}$. Since $\widehat{X}_i \leq (1+\varepsilon)X_i'$, 
	\[
	\widehat{X}(B) \leq (1+\varepsilon) X'(B)
	\]
	for any set of items $B$. Therefore, if $\widehat{X}(\widehat{B}_{j,k}) > 1+4\varepsilon$, then $X'(B_{j,k}') > 1+2\varepsilon$ and so bin $B_{j,k}'$ must have been broken by $\mathcal{P}'$. Then,
	\[
	\cost_{1+4\varepsilon}\left( \widehat{\mathcal{P}}, \widehat{X} \right)  \leq \cost_{1+2\varepsilon}(\mathcal{P}',X')
	\]
	and we obtain the desired result.\hfill\halmos
\endproof

\subsection{From Discretized Policy to Regular Policy with Resource Augmentation}

The main result of this section is the following.
\begin{theorem}\label{thm:discrete_to_regular}
	For any policy $\widehat{\mathcal{P}}$ that sequentially packs items $\widehat{X}_1,\dotsc,\widehat{X}_n$ into bins of size $1+4\varepsilon$, there exists a policy $\mathcal{P}$ that sequentially packs items $X_1,\dotsc,X_n$ into bins of size $1+6\varepsilon$ such that
	\[
	\cost_{1+6\varepsilon}(\mathcal{P}, X) \leq (1+\varepsilon) \cost_{1+4\varepsilon}\left( \widehat{\mathcal{P}}, \widehat{X} \right).
	\]
\end{theorem}

Given a policy $\widehat{\mathcal{P}}$ for the discretized items $\widehat{X}_1,\ldots,\widehat{X}_n$, we recover a policy $\mathcal{P}$ for items $X_1,\ldots,X_n$ with an extra $2\varepsilon$ in the bins' capacities. We couple the variables $\widehat{X}_i$ with $X_i'$ and $X_i$. Policy $\mathcal{P}$ simulates policy $\widehat{\mathcal{P}}$ in the following manner. For each bin $\widehat{B}_j$ that policy $\widehat{\mathcal{P}}$ opens, $\mathcal{P}$ opens a bin $B_{j,1}$ and packs items in $B_{j,1}$ as policy $\widehat{\mathcal{P}}$ would do in bin $\widehat{B}_j$, as long as $|X(B_{j,1}) - X'(B_{j,1})|\leq \varepsilon$ holds. (Note that the comparison is between random variables $X$ and $X'$.) If this difference is violated, policy $\mathcal{P}$ opens a new bin $B_{j,2}$, and continues following $\widehat{\mathcal{P}}$ as long as $|X(B_{j,2}) - X'(B_{j,2})|\leq \varepsilon$ holds, and so on.

As before, we need to show that policy $P$ is well defined, in the sense that bin $B_{j,k}$ is not broken if $\widehat{B}_j$ has not been broken.
\begin{proposition}
	If $\mathcal{P}$ overflows bin $B_{j,k}$ for some $k$, then $\widehat{\mathcal{P}}$ must have overflowed bin $\widehat{B}_j$. Moreover, at most one of the bins $B_{j,1},B_{j,2},\dotsc$ can be overflowed.
\end{proposition}
\proof{Proof.}
	Let $B_{j,k,t}$ be the items packed into bin $B_{j,k}$ by policy $\mathcal{P}$ up to time $t$. Suppose that at time $t$ policy $\mathcal{P}$ breaks bin $B_{j,k}$; then $X(B_{j,k,t}) > 1 + 6\varepsilon$. Now,
	\begin{align*}
		\widehat{X}(\widehat{B}_{j}) & \geq \widehat{X}(B_{j,k,t}) \\
		& \geq  X'(B_{j,k,t}) \\
		& = X'(B_{j,k,t-1}) + X_t' \\
		& \geq \left(X(B_{j,k,t-1}) - \varepsilon\right) + \left( X_t -\varepsilon^4  \right)\\
		& > 1 + 5\varepsilon - \varepsilon^4  > 1+4\varepsilon.
	\end{align*}
	Therefore, $\widehat{\mathcal{P}}$ must have overflowed bin $\widehat{B}_j$.\hfill\halmos
\endproof

As a consequence we have the following result.
\begin{proposition}	$\E\left[ O_{\mathcal{P}} \right] \leq \E\left[ O_{\widehat{\mathcal{P}}} \right].$
\end{proposition}

For $k\geq 1$, consider the family of events $\mathcal{E}_{j,k}= \{ |X(B_{j,k}) - X'(B_{j,k}) | > \varepsilon \}$ and for $k=0$ define $\mathcal{E}_{j,0}=\{ \mathcal{P} \text{ opens bin } B_{j,1} \} = \{ \widehat{\mathcal{P}} \text{ opens bin } \widehat{B}_j \}$. Then,
\begin{proposition}
	For any $k \geq 1$,
	\[
	\Prob(\mathcal{E}_{j,k} \mid \mathcal{E}_{j,k-1}, \ldots, \mathcal{E}_{j,1}, \mathcal{E}_{j,0} ) \leq (6\varepsilon^2)\Prob(\mathcal{P} \text{ opens }B_{j,k}\mid \mathcal{E}_{j,k-1}, \ldots, \mathcal{E}_{j,1}, \mathcal{E}_{j,0} ).
	\]
\end{proposition}
\proof{Proof.}
	The proof is identical to Proposition~\ref{prop:key_lem_bound_on_deviation_discr_1}.\hfill\halmos
\endproof

Let $N_{\mathcal{P},j}$ be the number of bins $B_{j,1},B_{j,2},\ldots$ that policy $\mathcal{P}$ opens. Then, $N_\mathcal{P} = \sum_{j=1}^n N_{\mathcal{P},j}$ is the number of bins used by policy $\mathcal{P}$. Following the proof strategy used for Proposition~\ref{prop:key_bounded_bins_1}, we obtain the following proposition.
\begin{proposition}
	$\E\left[ N_{\mathcal{P},j} \right] \leq (1+\varepsilon)\Prob(\widehat{\mathcal{P}} \text{ opens }\widehat{B}_j)$.
\end{proposition}


\begin{corollary}
	$	\E[N_{\mathcal{P}}] \leq (1+\varepsilon) \E[N_{\widehat{\mathcal{P}}}] $.
\end{corollary}

\proof{Proof of Theorem~\ref{thm:discrete_to_regular}.}
	The proof is direct from the previous results.\hfill\halmos
\endproof

\subsection{Computing an Optimal Discretized Policy via Dynamic Programming}

We can write a dynamic program (DP) that computes $\min_{\widehat{P}} \cost_{1+4\varepsilon}\left( \widehat{\mathcal{P}}, \widehat{X} \right)$, solved by backward induction in $\mathcal{O}\left( \frac{1}{\varepsilon^{10}} n^{2/\varepsilon^5}  \right)$ time. The states are pairs $(t,S)$, where $t=1,\ldots,n+1$ and $S=(k_0,k_1,\ldots,k_r)$ is a vector of non-negative integers such that $k_0+ k_1+\cdots + k_r \leq t-1$. Here $k_j$ represents the number of bins currently at capacity $j\cdot \varepsilon^5$, $j=1,\ldots, \lceil 2/\varepsilon^5\rceil$. The number of states $(t,S)$ is at most $\mathcal{O}(n^{2/\varepsilon})$. Then, the DP recursion becomes
\begin{align*}
	v(t,S) = \min&\left\{ 1 + \E_{\widehat{X}_t} \left[  v(t+1,S  + e_{\widehat{X}_t/\varepsilon^5}) \right] \right.,\\
	&\left.  C\Prob(k_j + \widehat{X}_t > 1+4\varepsilon) + \E_{\widehat{X}_t} \left[  v(t+1,S  + e_{j+\widehat{X}_t/\varepsilon^5} - e_j  ) \right] : 0< j \varepsilon^5 \leq 1+4\varepsilon , k_j \geq 1\right\} ,
\end{align*}
with the boundary condition $v(n+1,S)=0$ for any $S$. Here, $e_j$ is the canonical vector in $\R^{r+1}$ with a $1$ in the $j$-th coordinate and $0$ elsewhere. The recursion for $v(t,S)$ includes the two possible choices for a decision maker: Pack the item into a new bin and incur a cost of $1$ or use one of the previously opened and available bins.

Given access to $v(t+1,S')$ for any valid $S'$, we can compute $v(t,S)$ in $\mathcal{O}(1/\varepsilon^{10})$ time, since we need to compute the corresponding expectations in time $\mathcal{O}(1/\varepsilon^5)$. There are $\mathcal{O}(1/\varepsilon^5)$ of these terms inside the minimum operator, so we can compute $v(t+1,S')$ in $\mathcal{O}(1/\varepsilon^{10})$ time. Finally, given that there are $\mathcal{O}(n^{2/\varepsilon^5})$ states, we obtain the stated running time.

\section{Numerical Experiments}\label{sec:numbers}

In this section, we empirically validate the Budgeted Greedy (BG) algorithm. We quantify an algorithm's performance via the ratio of its cost to the cost incurred by some reference algorithm. When computationally possible, the reference algorithm is the optimal offline sequential policy. Otherwise, the reference is BG itself. We compare BG against the online benchmarks Full Greedy (FG), Fixed-Threshold (FT) and Fixed-Threshold Greedy (FTG). 
\begin{itemize}[leftmargin=*, itemsep=0em]
	\item \textbf{Full Greedy (FG)} is the myopic policy that for each item $i$ compares the instantaneous cost of opening a new bin (unit cost) and the expected cost of packing the item in one of the previously opened bins, $C \Prob(\text{overflow})$. The policy selects the cheapest option.
	
	\item \textbf{Fixed-Threshold (FT$(\alpha)$)} is the policy that has a threshold $\alpha \in (0, 1]$, and packs items into a bin as long as its usage does not exceed $\alpha$. Note that this policy uses one bin at a time.
	
	\item \textbf{Fixed-Threshold-Greedy (TG$(\alpha)$)} combines the myopic policy FG with a capacity threshold $\alpha$. The policy behaves as FG, but bins with usage greater than $ \alpha$ are discarded. Note that FG corresponds to TG$(1)$.
%
%
%
%
\end{itemize}
%

We test BG on four kinds of instances, one i.i.d.\ sequence of random variables, and three arbitrary exponential random variable input sequences. In all the instances we set the penalty to $C=50$ and input length to $n=10^5$. We simulate each instance $1,000$ times and report the sample mean.
\begin{itemize}[leftmargin=*, itemsep=0em]
\item \textbf{I.I.D.\ Sequence}. In this experiment, we consider an i.i.d.\ input sequence with three-point support given by
$$
X_i=\begin{cases}
		0 & \text{w.p.\ }1-1/C \\
		0.4 & \text{w.p.\ }1/2C \\
		0.61 & \text{w.p.\ }1/2C .
\end{cases}
$$
With this input, we aim to compare BG against TG$(\alpha)$ with threshold $\alpha \geq 0.4$. For $\alpha < 0.4$, TG$(\alpha)$ is near-optimal, therefore we do not study this case because we already include the optimal offline policy as a reference. Furthermore, it suffices to consider the case $\alpha=0.4$, since TG$(\alpha)$ for $\alpha \in (0.4, 1)$ is exactly the same. For TG$(1)$, we recover FG. In addition, we test different values of $\gamma$ for BG, denoted BG$(\gamma)$. We test $\gamma=1,2$ and the theoretically optimal $\gamma=\sqrt{2}$ given by Theorem~\ref{thm:iid_main_thm}. In this experiment we do not test FT, since it behaves exactly as TG for thresholds $\alpha < 1$, and FT$(1)$ has an expected cost of at least $ n$.

	\item \textbf{Exponential Distributions}. We consider input random variables $X_1,\dotsc,X_n$ that follow exponential distributions, $\Prob(X_i > x) = e^{-\lambda_i x}$. We perform three different experiments:
	\begin{enumerate}[leftmargin=*, itemsep=0em]
		\item First, we consider an input sequence of exponential random variables with increasing rates. The smallest rate starts at $\lambda_1=\log C$ and the largest rate is $\lambda_n = 3\log C$. In general, we set $\lambda_i = \left( 1 + 2\frac{i-1}{n-1} \right)\log C$ for $i=1,\dotsc,n$.
		
		\item Second, we consider an input sequence with decreasing rates. The largest rate is $\lambda_1 = 3\log C$ and the smallest rate is $\lambda_n=\log C$. In this case, we have $\lambda_i = \left( 3 - 2\frac{i-1}{n-1} \right)\log C$ for $i=1,\dotsc, n$.
		
		\item Finally, we consider an input sequence divided into three sections, each section with an i.i.d.\ sequence. The first section, for $i=1,\ldots, \lfloor n/3 \rfloor$, considers the fixed rate $\lambda_i= \log C$. The second section, for $i=\lfloor n/3\rfloor+1,\ldots,\lfloor 2n/3\rfloor$, considers the fixed rate $\lambda_i= 2\log C$. The final section, for $i=\lfloor 2n/3 \rfloor+1,\ldots,n$ considers the fixed rate $\lambda_i=\log C$.
	\end{enumerate}
Theorem~\ref{thm:exp_main_bound_thm} guarantees that BG has a constant multiplicative factor loss if the rates are $ 2\log C$ or greater. In these experiments, we empirically test the expected cost incurred by BG when the rates are in $[\log C, 3\log C]$, where Theorem~\ref{thm:exp_main_bound_thm} can only guarantee a multiplicative loss of $\mathcal{O}(\log C)$. Moreover, this guarantee theoretically applies to large $C$ (see Section~\ref{sec:expo_analysis}); here we test the algorithm on the relatively small penalty $C=50$.
\end{itemize}

\subsection{Results}

\begin{figure}[h!]
	\centering
	\includegraphics[width=1\linewidth]{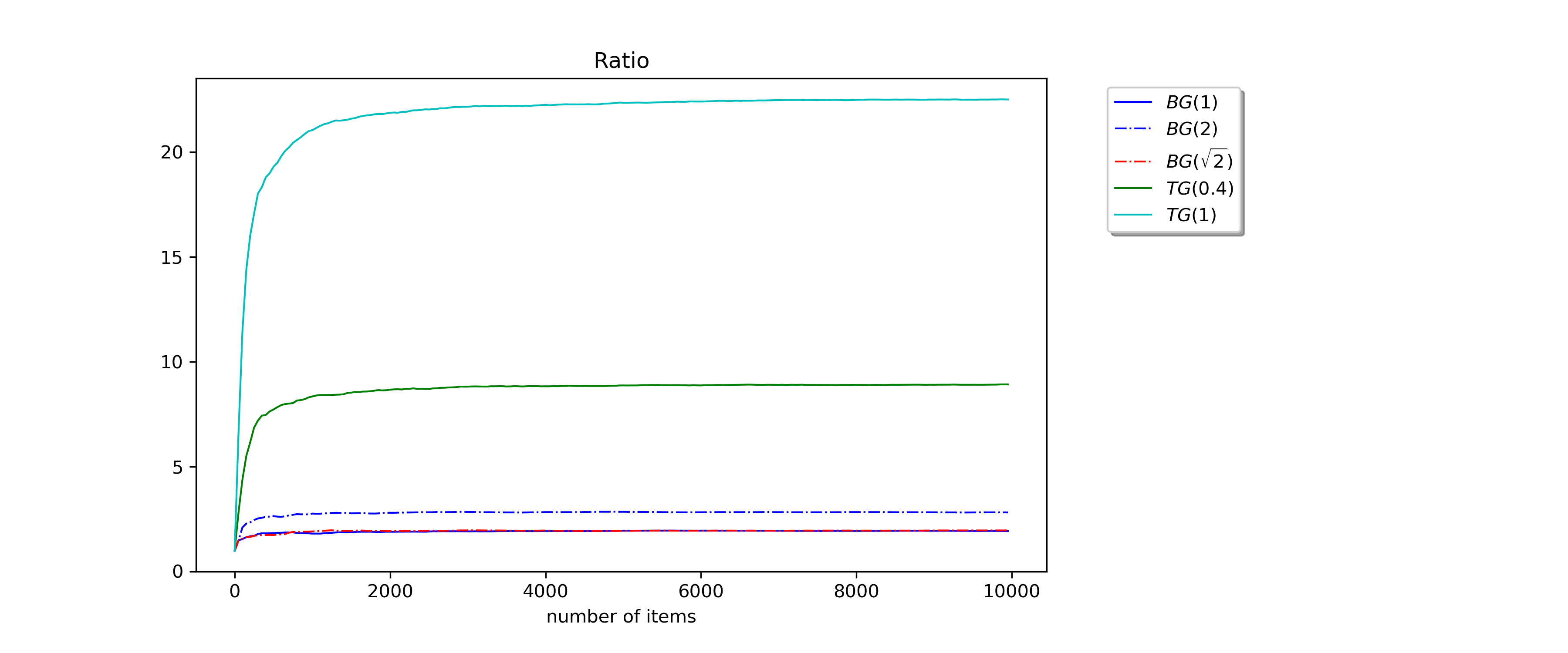}
	\caption{Ratio to optimal expected cost incurred by the algorithms BG, FG and TG. Note that BG$(\sqrt{2})$ overlaps with BG(1); the difference is roughly 0.1 units.}
	\label{fig:exp3pointdist}
\end{figure}

\paragraph{I.I.D.\ Random Variables} Figure~\ref{fig:exp3pointdist} presents the ratio of the sample mean of the cost incurred by the algorithms and the sample mean of the optimal offline sequential cost; this latter quantity is roughly $n/C$. 
We empirically confirm that the expected cost of TG policies is at least $n/8$; BG$(\gamma)$ with $\gamma = 1, \sqrt{2}, 2$ exhibits better performance. As $n$ grows, BG$(1)$ has a ratio of roughly $1.8$, BG$(2)$ one of roughly $2.75$, and BG$(\sqrt{2})$ a ratio of roughly $1.9$; theoretically we can guarantee a ratio of $\sqrt{3}+ 2\sqrt{2} \approx 4.5604$.
%

\begin{figure}[h!]
	\centering
	\includegraphics[width=1\linewidth]{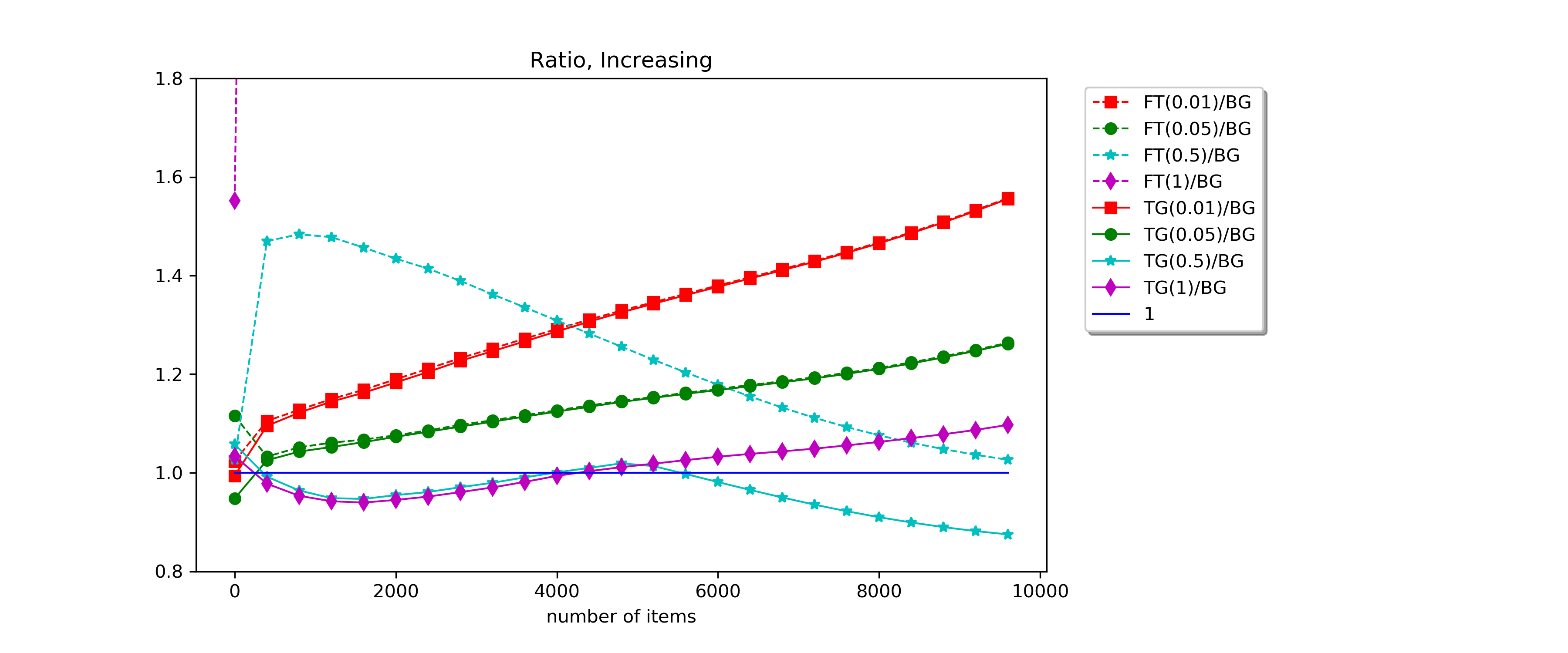}
	\caption{Ratio of cost incurred in the exponential case for increasing rates.}
	\label{fig:exponentialexp_increasing}
\end{figure}

\begin{figure}[h!]
	\centering
	\includegraphics[width=1\linewidth]{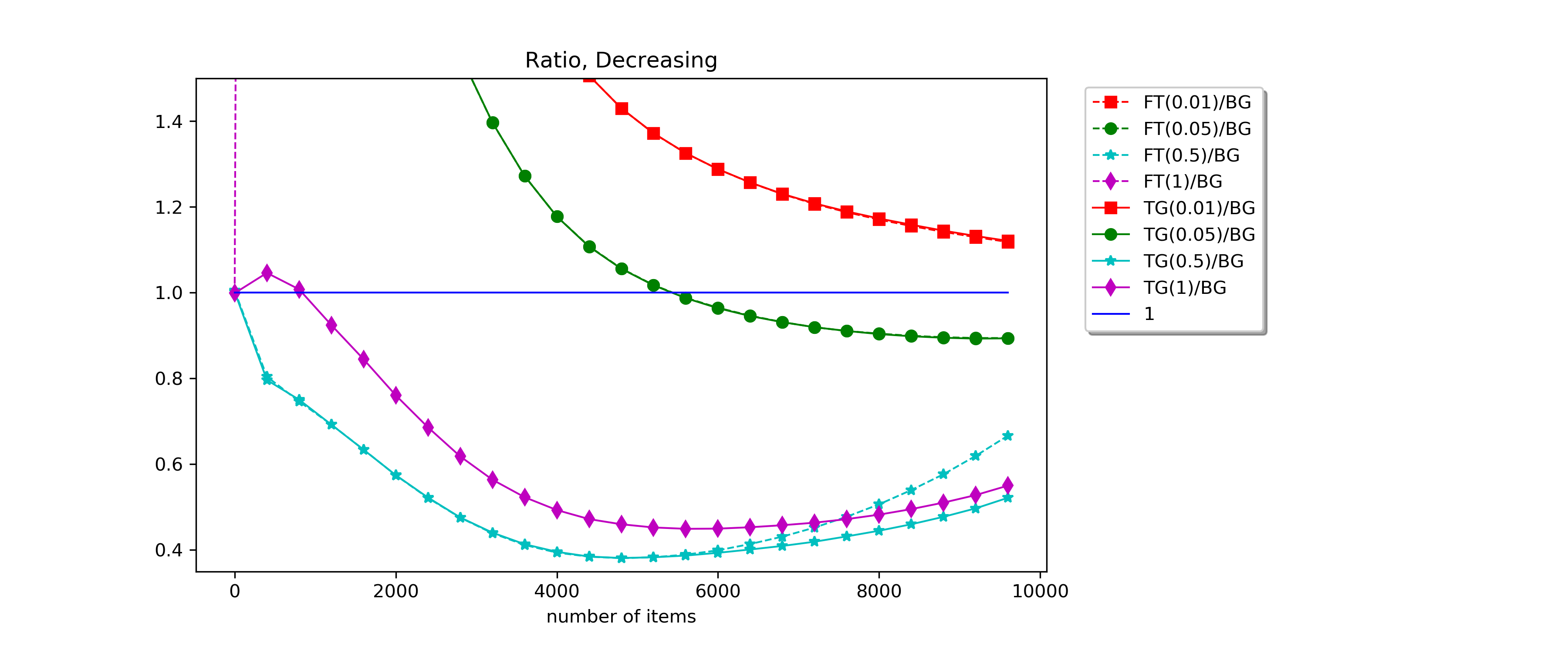}
	\caption{Ratio of cost incurred in the exponential case for decreasing rates.}
	\label{fig:exponentialexp_decreasing}
\end{figure}

\begin{figure}[h!]
	\centering
	\includegraphics[width=1\linewidth]{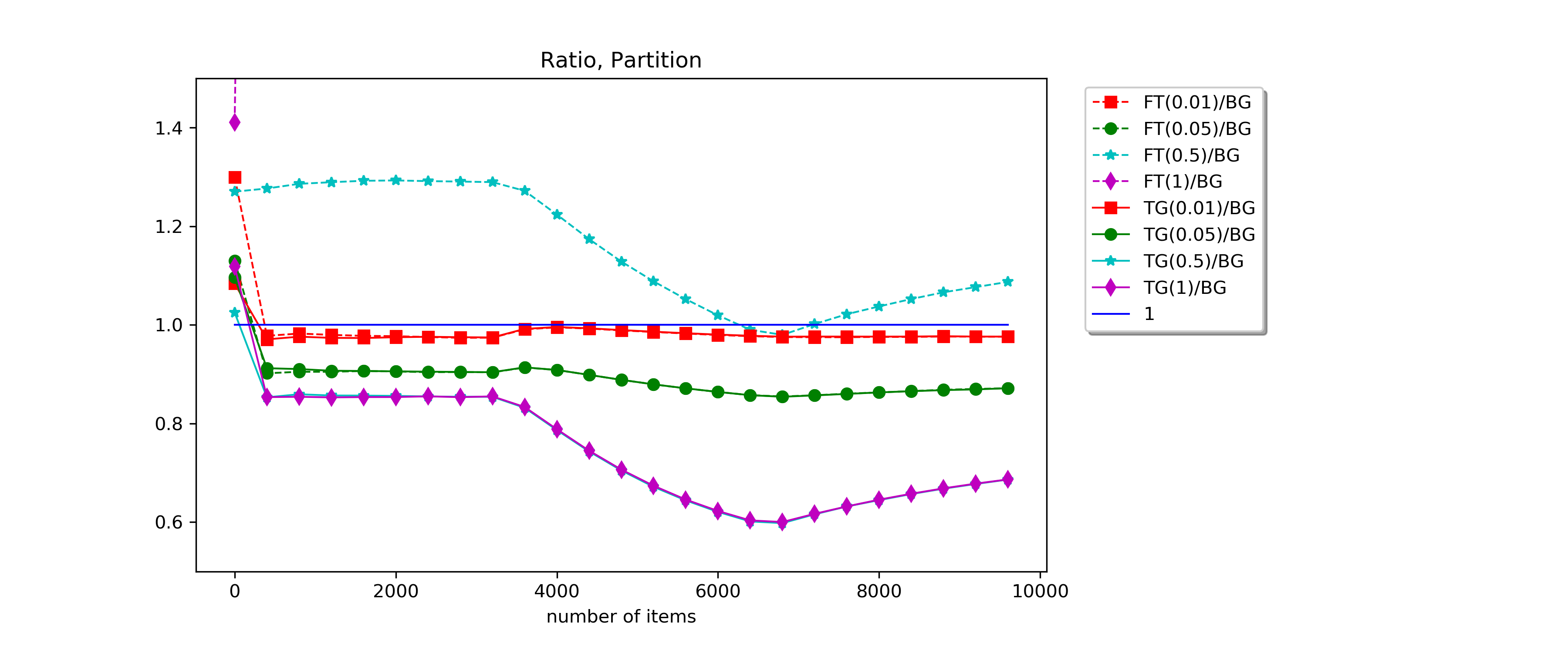}
	\caption{Ratio of cost incurred in the exponential case block input.}
	\label{fig:exponentialexp_partition}
\end{figure}

\paragraph{Exponential Random Variables} Figures~\ref{fig:exponentialexp_increasing},~\ref{fig:exponentialexp_decreasing} and~\ref{fig:exponentialexp_partition} present the empirical results of our experiments in the case of increasing rates, decreasing rates and block-input rates, respectively. 
In these experiments, we used BG$(2)$ as a reference, because computing the offline benchmark was too computationally expensive. We used $\gamma=2$ because it has a $\mathcal{O}(\log C)$ approximation guarantee compared to the optimal offline expected cost (see Proposition~\ref{prop:exponential_log_bound}.) Smaller values of $\gamma$ do not improve the performance of BG in a significant manner; as the results show, being greedy seems suited to exponential distributions. On the other hand, larger values of $\gamma$ make BG's performance resemble  FG.


Figure~\ref{fig:exponentialexp_increasing} displays the ratio of cost sample means between the benchmark algorithms and BG$(2)$ for increasing rates. 
We empirically observe that BG performs significantly better against all FT policies. Similarly, BG performs better that most of TG policies, with the exception of TG$(0.5)$ and TG$(1)$. Until approximately the $5,000$-th item, the ratio TG$(1)$/BG$(2)$ is the best among all greedy strategies, and afterwards the ratio TG$(0.5)$/BG$(2)$ becomes the best. Moreover, by the end of the sequence, FT$(0.5)$ becomes better than TG$(1)$.
During the whole input sequence, we  empirically observe that BG(2) is able to balance the behavior of TG$(0.5)$ and TG$(1)$, surpassing the performance of TG$(1)$ in the second half of the input sequence. For the entire sequence, the best performing algorithm's expected cost ratio is above $0.8$, which means BG(2) is within $25\%$ of the best algorithm for all input sizes. Furthermore, around the $5,000$-th item BG performs the best among all tested strategies. 
%

Figure~\ref{fig:exponentialexp_decreasing} displays the ratios of the tested algorithms and BG$(2)$ for the decreasing rates experiment. 
In this case, most of the TG/BG-curves and FT/BG-curves overlap, with the exception of TG(1)/BG and FT(1)/BG. 
For almost the entire sequence, all plots lie above $0.4$, indicating that BG's expected cost is at most $ 2.5 $ times the best performing algorithm's cost for all input sizes. BG's performance decreases until around the $5,000$-th item and improves thereafter. As in the previous experiment, BG performs better for larger rates, which coincides with our theoretical findings.
%

Figure~\ref{fig:exponentialexp_partition} displays the ratios between the tested algorithms and BG$(2)$ for the partitioned input sequence. 
The best performing algorithm over the entire sequence is TG(1); this algorithm's plot and all others lie above $0.6$, indicating BG is within $67\%$ of the best performing algorithm for any input size.
During the first interval of the sequence, the ratios are roughly constant; the main differentiation occurs with the transition to the second interval, where TG(1) and TG(0.5) outperform BG. In the last interval, BG's performance again improves.
%

\section{Concluding Remarks}

In this paper, we introduced the adaptive bin packing problem with overflow. We introduced the notion of risk as a proxy for a capacity threshold, as typically used in deterministic settings. 
We showed that Budgeted Greedy incurs an expected cost at most a constant factor times the optimal expected cost of an offline policy when the input is an i.i.d.\ sequence of random variables. In the more general setting, we give similar results for arbitrary exponential random variables.

We extended the discussion by studying the offline sequential adaptive bin packing problem, in which the decision maker knows the sequence of random variables in advance and must pack them in this order. 
We devised a soft-capacity PTAS by utilizing a policy tracking argument, and showed that computing the cost of the optimal policy is $\#\Pol$-hard by relating it to counting problems. 
This offline cost corresponds to the online benchmark.

Unfortunately, Budgeted Greedy does not guarantee a constant approximation factor for general input sequences. 
Consider the input sequence $X_1,X_2,X_3,\dotsc,X_n$ defined as $X_1= 1/n$, $X_{2i} \sim \mathrm{Bernoulli}(1/C)$ and 
\[
X_{2i+1} = \begin{cases}
	1/n & \text{w.p. } 1-1/C^2 \\
	1   & \text{w.p. } 1/C^2 .
\end{cases} 
\]
Budgeted Greedy incurs an expected cost of $\Theta(n)$, while the optimal offline policy incurs an expected cost of at most $n/C+1$. 

This example motivates either seeking a general algorithm exhibiting a bounded competitive ratio, or showing an impossibility result. In~\cite{alaei2013online}, the authors study the online generalized assignment problem with a similar stochastic component as in our model. They are able to show a $1-\frac{1}{\sqrt{k}}$ competitive ratio for general arriving distributions. However, they assume large capacity, in the sense that no item takes up more than $1/k$ fraction from any bin. It is not clear how to utilize their techniques in a bin packing setting, as they are able to discard distributions that they deem unimportant. Moreover, in the bin packing problem a large capacity assumption would immediately imply a policy with constant approximation factor, by simply filling up the bins until some desired fraction of capacity.


\section*{Acknowledgments.}

The authors' work was partially supported by the U.S.\ National Science Foundation via grants CMMI 1552479, AF 1910423 and AF 1717947.


{\bibliographystyle{amsplain}\fontsize{10}{0}\selectfont\bibliography{biblio}
}

\newpage
\appendix
\small


\section{Missing Proofs}\label{sec:Appendix}

\subsection{Missing Proofs From Section~\ref{sec:thealgorithm}}

\begin{customproposition}{\ref{prop:breakbinpenalty}}
	Let $X_1,\ldots,X_n$ be nonnegative independent random variables and any (deterministic) policy $\mathcal{P}$ for packing these items sequentially. Then, the expected number of broken bins by the policy $\mathcal{P}$ is given by
	\[
	\E[B_\mathcal{P}] = \sum_{j=1}^n \E_{X_1,\ldots,X_n} \left[  \sum_{i=1}^n \Prob_{X_i}(X_i + S_{j}^{i-1} > 1)\mathbf{1}_{\{ i \to j \}}^{\mathcal{P}}  \right]
	\]
	where $S_{j}^{i-1}$ is the level of bin $j$ at the beginning of iteration $i$ and $\mathbf{1}_{\{i\to j\}}^\mathcal{P}$ is the 0/1 indicator random variable of the event: Policy $\mathcal{P}$ packs item $X_i$ into bin $j$.
\end{customproposition}

\proof{Proof.}
	We write $\mathbf{1}_{\{i\to j\}}$ to denote $\mathbf{1}_{\{i\to j\}}^{\mathcal{P}}$. We show that $\Prob(\mathcal{P} \text{ breaks bin } j) = \E\left[ \sum_{i=1}^n \Prob_{X_i}(X_i + S_{j}^{i-1} > 1)\mathbf{1}_{\{ i \to j \}} \right]$. We have,
	\begin{align*}
	\E\left[ \sum_{i=1}^n \Prob_{X_i}(X_i + S_{j}^{i-1} > 1)\mathbf{1}_{\{ i \to j \}} \right] & = \sum_{i=1}^n \E\left[ \Prob_{X_i}(X_i+ S_j^{i-1} > 1) \mathbf{1}_{\{i\to j \}} \right]
	\end{align*}
	Observe that $S_{j}^{i-1}=\sum_{k\leq i-1} X_k \mathbf{1}_{\{ k\to j\}}$ and $\mathbf{1}_{\{ i \to j\}}$ only depend on the outcomes of $X_1,\ldots,X_{i-1}$. Therefore,
	\begin{align*}
	\E\left[ \Prob_{X_i}(X_i+ S_j^{i-1} > 1) \mathbf{1}_{\{i\to j \}} \right] & = \E_{X_1,\ldots,X_{i-1}}\left[ \Prob_{X_i}(X_i+ S_j^{i-1} > 1) \mathbf{1}_{\{i\to j \}} \right]\\
	& = \E_{X_1,\ldots,X_{i-1}}\left[ \E_{X_i}\left[ \mathbf{1}_{\{ X_i + S_j^{i-1} > 1 \}}  \right] \mathbf{1}_{\{ i\to j \}}  \right]\\
	& = \E_{X_1,\ldots,X_i} \left[ \mathbf{1}_{\{  X_i + S_j^{i-1} > 1 \}} \mathbf{1}_{\{i \to j\}}  \right].
	\end{align*}
	Clearly, $\{ X_i + S_j^{i-1} > 1 , i\to j \}= \{  X_i \text{ breaks bin } j \}$. Thus,
	\[
	\E\left[ \Prob_{X_i}(X_i+ S_j^{i-1} > 1) \mathbf{1}_{\{i\to j \}} \right]  = \Prob(X_i \text{ breaks bin } j),
	\]
	hence,
	\begin{align*}
	\E\left[ \sum_{i=1}^n \Prob_{X_i}(X_i + S_{j}^{i-1} > 1)\mathbf{1}_{\{ i \to j \}} \right] & = \sum_{i=1}^n \Prob(X_i \text{ breaks bin } j)= \Prob(\mathcal{P} \text{ breaks bin }j)
	\end{align*}
	since the last sum uses the fact that bins are overflowed at most once; hence the events  $\{ X_i \text{ breaks bin } j \}_i$ are disjoint.\hfill\halmos
\endproof

\begin{customproposition}{\ref{prop:sizeinbins}}
	For any sequence of nonnegative i.i.d.\ random variables $X_1,\ldots,X_n$, for any bin $B=B_j$ and any policy $\mathcal{P}$, we have
	\[
	\E\left[  \sum_{i\in B} \E [X_i\wedge 1] \right] = \E\left[  \sum_{i\in B}  X_i\wedge 1 \right]  \leq 2 \Prob(\mathcal{P}\text{ opens bin } B),
	\]
	where $X_i\wedge 1 = \min\{X_i,1\}$.
\end{customproposition}
\proof{Proof.}
The proof follows from a result in~\cite{dean2008approximating}, which we replicate here for completeness. Let $\mu_i= \E[X_i \wedge 1]$ be the normalized expected size of an item. Let $B^t$ be the (random) items that the policy packs into bin $B$ by time $t$. We are interested in the expectation of $\mu(B)= \sum_{i\in B}\mu_i = \mu(B^n)$. The random variables $\mu(B^t)$ are nondecreasing in $t$; by the monotone convergence theorem,
\[
\E\left[ \mu(B) \right] = \sup_{t\geq 0} \E[\mu(B^t) ].
\]
Now, the random variables $Z^t = \sum_{i\in B^t} (X_i\wedge 1) - \mu_i$ form a martingale. Indeed,
\begin{align*}
	\E[ Z^{t} \mid Z^{t-1}, t\to j ] = Z^{t-1} + \E[X_t\wedge 1] - \mu_{t} = Z^{t-1};
\end{align*}
then, for all $ t$, $\E[Z^t ]= Z^0 = 0$ and so $\E[ \mu(B^t)] = \E \left[ \sum_{i\in B^t} X_i \wedge 1 \right]\leq 2\Prob(\mathcal{P}\text{ opens bin } B)$; this last inequality holds because we break the bin at most once and we must have opened the bin. Therefore
\[
\E \left[ \mu(B) \right] = \sup_{t\geq 0} \E[\mu(B^t)] \leq 2\Prob(\mathcal{P}\text{ opens bin } B).
\]
\hfill\halmos
\endproof

\begin{customproposition}{\ref{prop:sizesinbins_cost}}
	For any sequence of nonnegative i.i.d.\ random variables $X_1,\ldots,X_n$, for any policy, we have
	\[
	\E\left[ \sum_{i=1}^n (X_i \wedge 1 )\right] \leq \cost(\mathcal{P}).
	\]
\end{customproposition}

\proof{Proof.} Note that by Proposition~\ref{prop:sizeinbins} we have
\[
\E\left[ \sum_{i=1}^n (X_i \wedge 1)  \right] = \sum_{j=1}^n\E\left[ \sum_{i\in B_j^\mathcal{P}} (X_i \wedge 1) \right] = \sum_{j=1}^n \E\left[ \sum_{i\in B_{j}^{\mathcal{P}}} \E[X_i\wedge 1] \right] = \E\left[ \sum_{i=1}^n \E[X_i \wedge 1] \right].
\]
We only need to show that $\sum_{i=1}^n \E[X_i \wedge 1]$ is a lower bound for $\cost(\OPT)$. We can compute $\cost(\OPT)$ recursively via dynamic programming as follows. We define the states as vectors $S\in (\R\cup \{ \emptyset\})^n$ where $S_j \in \R$ is the usage of $j$-th bin and $S_j=\emptyset$ means that bin $j$ is closed. We consider $\emptyset$ as an special symbol such that $a+\emptyset = a$ for any $a\in \R$. With this, the optimal cost can be computed via the following recursions:
\begin{align*}
	v_{t}(S) & = \inf \left\{\E_{X_t}\left[ v_{t+1}( S + X_t e_j  )  \right] : j=1,\ldots, n, S_j \in [0,1]\cup \{\emptyset\} \right\}, \quad \forall t =1,\ldots,n, \forall S  \\
	v_{n+1}(S) & = \sum_{j=1}^n \mathbf{1}_{\{ S_j=\emptyset \}} + C \sum_{j=1}^n \mathbf{1}_{\{ S_j > 1 \}}, \quad \forall S.
\end{align*}
The second equation measures the overall cost accumulated at the end of processing the sequence $X_1,\ldots,X_n$. The first equation takes actions that minimizes the mean cost of sample paths. Note that items can only be packed into bins not opened ($\emptyset$) or bins with usage $\leq 1$. Using MDP theory, we can show $v_1(\emptyset,\ldots,\emptyset)= \cost(\OPT)$, which we skip here for brevity.

Now, consider the functions $u_t(S) = \sum_{\tau =t}^n \E[X_\tau \wedge 1 ] + \sum_{j=1}^n (S_j \wedge 1)\mathbf{1}_{\{S_j > 0\}} $ for any $S$. We show by backward induction in $t=n+1,\ldots,1$ that $u_t(S)\leq v_t(S)$. For $t=n+1$ we have
\[
u_{n+1}(S) = \sum_{j=1}^n (S_j \wedge 1)\mathbf{1}_{\{S_j >0 \}} \leq \sum_{j=1}^n \mathbf{1}_{\{S_j > 0\}} \leq v_{n+1}(S).
\]
Now, assume the result is true for $t+1$ and let us show it for $t$. Let $j=1,\ldots,n$ with $S_j\in [0,1]\cup \{ \emptyset \}$, then
\begin{align*}
	\E_{X_t} \left[ v_{t+1}(S + X_t e_j) \right] & \geq \E_{X_t}\left[ \sum_{\tau = t+1}^n \E[X_\tau \wedge 1] + \sum_{\substack{k=1\\ k \neq j}}^n (S_k \wedge 1)\mathbf{1}_{\{ S_j > 0\}} + (S_j + X_t)\wedge 1  \right] \\
	& \geq \sum_{\tau =t+1}^n \E[X_\tau \wedge 1] + \E_{X_t}[X_t\wedge 1] = u_{t}(S).
\end{align*}
Taking minimum in $j$, we conclude $v_{t}(S)\geq u_t(S)$ for any $S$.

Now, for $t= 1$ we have $\cost(\OPT) = v_1(\emptyset,\ldots,\emptyset) \geq u_1(\emptyset,\ldots, \emptyset) = \sum_{t=1}^n \E[X_t\wedge 1]$ which finishes the proof.
\halmos
\endproof

\subsection{Missing Proofs From Section~\ref{sec:iid_analysis}}

\begin{customlemma}{\ref{lem:cost_reduce_tree_analysis}}
	$\cost_{\ell, \widehat{c}}(\mathcal{T}_\mathcal{P}(u)) \geq \cost_{\ell',c'}(\mathcal{T}_{\mathcal{P}'}(u))$.
\end{customlemma}

\proof{Proof.}
	We define
	\[
	\cost_{\ell,\widehat{c}}(\mathcal{T}_\mathcal{P}(u))_j = 1+\E\left[ (C+2\delta) \sum_{i=1}^n \mathbf{1}_{\{i\to j\}}^{\mathcal{P}}  \mathbf{1}_{\{ X_i + S_j^{i-1} > 1 \}}^{\mathcal{P}}\mid \text{Reach node } u\right]
	\]
	which is the original cost paid in $\mathcal{T}_\mathcal{P}$ when packing items into bin $j$ after reaching node $u$ in the tree. We also define
	\begin{align*}
	\cost_{\ell',c'}(\mathcal{T}_{\mathcal{P}'}(u))_j = & 1+ \E\left[ (C+\delta)\sum_{i=1}^n \mathbf{1}_{\{i\to j\}}^{\mathcal{P}'} \mathbf{1}_{\{X_i + S_j^{i-1} > 1\}}^{\mathcal{P}'} \right. \\
	&\left. + (C+2\delta)\left(\sum_{i=1}^n \mathbf{1}_{\{i\to j'\}}^{\mathcal{P}'} \mathbf{1}_{\{ X_i + S_{j'}^{i-1} > 1   \}}^{\mathcal{P}'}\right) + \mathbf{1}_{\{\text{Open bin }j'\}}^{\mathcal{P}'} \mid \text{Reach node } u \right]
	\end{align*}
	which is the new cost paid by $\mathcal{T}_{\mathcal{P}'}$ when packing items into bin $j$ after reaching node $u$ and the new cost incurred by packing items into bin $j'$.
	
	Therefore, the variation of the cost $\cost_{\ell, \widehat{c}}(\mathcal{T}_\mathcal{P}(u)) - \cost_{\ell',c'}(\mathcal{T}_{\mathcal{P}'}(u))$ is given by
	\begin{align*}
	\cost_{\ell, \widehat{c}}(\mathcal{T}_\mathcal{P}(u)) - \cost_{\ell',c'}(\mathcal{T}_{\mathcal{P}'}(u) ) & = (\cost_{\ell,\widehat{c}}(\mathcal{T}_\mathcal{P}(u))_j - \cost_{\ell',c'}(\mathcal{T}_{\mathcal{P}'}(u))_j).
	\end{align*}
	Now, we always have
	\[
	\mathbf{1}_{\{i\to j\}}^\mathcal{P} = \mathbf{1}_{\{i\to j\}}^{\mathcal{P}'} + \mathbf{1}_{\{i\to j'\}}^{\mathcal{P}'},
	\]
	for all $i=1,\ldots,n$. Indeed, if we are in a branch not containing $u$, then $\mathcal{P}$ and $\mathcal{P'}$ behave the same and there is no bin $j'$. If we are in a branch containing $u$, and if we pack $i$ into $j$, we either pack $i$ into $j$ before surpassing the risk budget in which case $\mathcal{P}$ and $\mathcal{P}'$ behave the same or we do it after surpassing the risk budget in which case $i$ goes to $j'$. With this fact we have,
	\begin{align*}
	\cost_{\ell, \widehat{c}}(\mathcal{T}_\mathcal{P}(u))) - \cost_{\ell',c'}(\mathcal{T}_{\mathcal{P}'}(u)) & = (\cost_{\ell,\widehat{c}}(\mathcal{T}_\mathcal{P}(u))_j - \cost_{\ell',c'}(\mathcal{T}_{\mathcal{P}'}(u))_j) \\
	& \geq \E\left[ \delta \left(\sum_{i=1}^n \mathbf{1}_{\{i\to j\}}^{\mathcal{P}'} \mathbf{1}_{\{ X_i + S_j^{i-1} > 1 \}}^{\mathcal{P}'}\right) - \mathbf{1}_{\{ \text{Open bin }j' \}}^{\mathcal{P}'} \mid \text{Reach node }u \right],
	\end{align*}
	in the last inequality we used the fact that the cost of breaking the bin $j'$ is smaller than the cost of breaking $j$ at that point of the computation. This is true since the usage of bin $j'$ is at most the usage of $j$ at the same point of computation. Now, for $i\geq j$,
	\begin{align*}
	\E\left[ \mathbf{1}_{\{ i\to j\}}^{\mathcal{P}'} \mathbf{1}_{\{ X_i + S_j^{i-1}> 1 \}} \mid \text{Reach node }u  \right] & = \E\limits_{X_1,\ldots,X_{i-1}} \left[ \E_{X_i}\left[ \mathbf{1}_{\{ i\to j \}}^{\mathcal{P}'} \mathbf{1}_{\{ X_i + S_{j}^{i-1} > 1 \}}^{\mathcal{P}'} \right] \mid \text{Reach node }u \right] \\
	& =\E\limits_{X_1,\ldots,X_{i-1}} \left[ \mathbf{1}_{\{ i\to j \}}^{\mathcal{P}'} \E_{X_i}\left[  \mathbf{1}_{\{ X_i + S_{j}^{i-1} > 1 \}} \right] \mid \text{Reach node }u \right] \\
	& = \E\limits_{X_1,\ldots,X_{i-1}} \left[ \mathbf{1}_{\{i\to j\}}^{\mathcal{P}'} \Prob_{X_i}(X_i + S_j^{i-1} > 1 ) \mid \text{Reach node }u \right].
	\end{align*}
	This is because the event $\{\text{Reach node }u\}$ is determined by the outcomes of $X_1,\ldots,X_{j-1}$. While for $i\leq j-1$ we have
	\[
	\E\left[ \mathbf{1}_{\{ i\to j\}}^{\mathcal{P}'} \mathbf{1}_{\{ X_i + S_j^{i-1}> 1 \}} \mid \text{Reach node }u  \right] = 0
	\]
	since bin $j$ is opened at node $u$ at level $j$. Thus,
	\begin{align*}
	\cost_{\ell, \widehat{c}}(\mathcal{T}_\mathcal{P}(u)) - \cost_{\ell',c'}(\mathcal{T}_{\mathcal{P}'}(u)) & =\E \left[\delta\sum_{i=j}^n  \mathbf{1}_{\{i\to j\}}^{\mathcal{P}'} \Prob(X_i + S_j^{i-1} > 1 )  - \mathbf{1}_{\{\text{Open bin }j' \}} \mid \text{Reach node }u \right]  \\
	& = \E\left[ \risk(B_j) - \mathbf{1}_{\{ \text{Open bin } j' \}}^{\mathcal{P}'} \mid \text{Reach node } u \right] \\
	& \geq \left( \delta \frac{\gamma}{C} - 1   \right) = 0. \tag{Using $\gamma = \frac{C}{\delta}$.}
	\end{align*}
\hfill\halmos
\endproof

\subsection{Missing Proofs From Section~\ref{sec:expo_analysis}}

\begin{customproposition}{\ref{prop:lower_bound_exp_1}}
	Let $\varepsilon >0$ and set $\beta = 6\frac{n_1 \log C}{\varepsilon}$ and $k=3\varepsilon \mu$. Then, running Budgeted Greedy with $\gamma=1$ in the input described in Subsection~\ref{subsec:exp_lower_bound} we have
	\[
	\cost(\ALG) \geq \frac{1}{2} n_1.
	\]
\end{customproposition}

\proof{Proof.}
	We show that (w.h.p.) Algorithm~\ref{alg:BA} packs each item $X_i^\lambda$ individually. This is achieved by showing that in between two $X_i^\lambda$ and $X_{i+1}^\lambda$, there is enough mass introduced by the elements $X_{i,j}^\mu$, therefore not allowing the items $X_i^\lambda$ to be packed together.

	Now, let $k = 3 \varepsilon \mu = 3\varepsilon \beta \log C$. Then,
	\[
	\E\left[ \sum_{i=1}^k X_i^\mu \right] = \frac{k}{\beta \log C} = 3\varepsilon,
	\]
	thus
	\begin{align*}
	\Prob\left( \sum_{i=1}^k X_i^\mu \leq 2\varepsilon  \right) & \leq \Prob\left( \left| \sum_{i=1}^k X_i^\mu - 3\varepsilon \right| \geq \varepsilon \right)\leq \frac{1}{\varepsilon^2} \frac{k}{(\beta \log C)^2}  = \frac{3}{\varepsilon \beta \log C}. \tag{Chebyshev inequality}
	\end{align*}
	Pick $\beta = 3 \frac{2n_1 \log C }{\varepsilon}$ and so,
	\begin{align*}
	\Prob\left( \exists  j= 1, \ldots, n_1 : \sum_{i=1}^k X_{j, i}^\mu \leq 2\varepsilon \right) & \leq n_1 \cdot \frac{3}{\varepsilon \beta \log C} = \frac{1}{2}.
	\end{align*}
	That is, with probability at least $\frac{1}{2}$, all blocks $X_{j,1}^\mu,\ldots, X_{j,k}^\mu$ add at least $2\varepsilon$ mass. Consider the event
	\[
	E = \left\{ \forall  j= 1, \ldots, n_1 : \sum_{i=1}^k X_{j, i}^\mu > 2\varepsilon  \right\}.
	\]
	then, we just proved that $\Prob(E)\geq \frac{1}{2}$.
	
	\begin{claim}
		Given event $E$, Algorithm~\ref{alg:BA} with $\gamma=1$ never packs $X_i^\lambda$ and $X_{i+1}^\lambda$ together for any $i$.
	\end{claim}
	
	\proof{Proof.}
		Suppose that Algorithm~\ref{alg:BA} packs $X_i^\lambda$ and $X_{i+1}^\lambda$ together for some $i$. This means that the algorithm had enough budget and space to allocate $X_{i+1}^\lambda$. Since $\Prob(X_{i+1}^\lambda > 1-x) = e^{-\lambda(1-x)}  > e^{-\mu(1-x)} = \Prob(X_{ij}^\mu > 1-x)$ for any $x >0$, then Budgeted Greedy must have packed all the $X_{ij}^\mu$ in between $X_{i}^\lambda$ and $X_{i+1}^\lambda$. However, under event $E$, these items increase the usage of the bin by at least $2\varepsilon$. Then, the budget utilized by $X_{i+1}^\lambda$ is at least
		\[
		\Prob(X_{i+1}^\lambda > 1-2\varepsilon) = e^{-\lambda(1-2\varepsilon)} >  e^{-(1-\varepsilon^2)\log C } >  \frac{1}{C}
		\]
		which is a contradiction to the risk budget of Budgeted Greedy.\hfill\halmos
	\endproof

	\begin{claim}
		Using the same choices of $\beta$ and $k$ as before, we have $\cost(\ALG) \geq \frac{1}{2}n_1$.
	\end{claim}
	
	\proof{Proof.}
		By the previous result, under event $E$, no $X_i^\lambda$ and $X_{i+1}^\lambda$ are packed together. Therefore, at least $n_1$ open bins are needed. Since $E$ occurs w.p. $\geq \frac{1}{2}$ we conclude the desired result.\hfill\halmos
	\endproof
	\hfill\halmos
\endproof

\begin{customproposition}{\ref{prop:upper_bound_opt_cost_1}}
	Using the same parameters and the same input as in the previous result, for any $\varepsilon >0$ such that $\varepsilon \log C \geq 4$, we have
	\[
	\cost(\OPT) \leq 48n_1 \left(  \frac{k}{\beta\log C} + \frac{1}{\varepsilon \log C} \right)=48 n_1 \left(  3\varepsilon + \frac{1}{\varepsilon \log C}  \right).
	\]
\end{customproposition}

\proof{Proof.} In order to show an upper bound for $\cost(\OPT)$ it is enough to exhibit a policy with cost bounded by the desired value. We consider the following budgeted policy $\mathcal{P}$ with risk budget $\frac{2}{C}$: Pack items with rate $\lambda$ separately of items with rate $\mu$. We are going to show that $\mathcal{P}$ opens at most
\[
\frac{n_2}{\beta \log C} + \frac{n_1}{\varepsilon \log C} = n_1 \left( \frac{k}{\beta \log C} + \frac{1}{\varepsilon \log C}  \right)
\]
bins in expectation (up to a constant). Since $\mathcal{P}$ is budgeted with budget $\frac{2}{C}$ we have $\cost(\OPT)\leq \cost(\mathcal{P})\leq 3 \E[N_\mathcal{P}]$ from which the result follows. In what follows we prove the bound over the number of bins.

Let us analyze the policy $\mathcal{P}$. Policy $\mathcal{P}$ opens two kind of bins; the first kind of bins only contain items following exponentials distribution of rate $\lambda$; the second kind of bins only contain items following exponential distribution of rate $\mu$. We have $N_\mathcal{P}= N_\mathcal{P}^1 + N_\mathcal{P}^2$ where $N_\mathcal{P}^1$ is the number of bins of type 1 and $N_\mathcal{P}^2$ is the number of bins of type $2$. An equivalent way to see this process is that policy $\mathcal{P}$ runs two copies of Algorithm~\ref{alg:BA}, one for the rate $\lambda$ and one for the rate $\mu$. Then, $N_\mathcal{P}^1$ equals $N_\ALG$ over the sequence $X_1^\lambda,\ldots,X_{n_1}^\lambda$ and $N_\mathcal{P}^2$ equals $N_\ALG$ over the sequence $X_{1,1}^\mu ,\ldots, X_{n_1,k}^\mu$.

The following lemma is a general result that allows us to bound the number of bins used in a nonnegative i.i.d. sequence of items under Algorithm~\ref{alg:BA}. For the sake of clarity, the proof has been moved to the end of this subsection.

\begin{lemma}\label{lem:bound_iid_walds}
	Suppose $X_1,\ldots,X_n$ are i.i.d. sequence of items, then $\E[N_{\ALG}]\leq \frac{2n-1}{\E[|B_1|]}$ where $|B_1|$ is the number of items packed in the first bin.
\end{lemma}

\begin{claim}
	Let $X_1,\ldots,X_n$ be $n$ independent exponential r.v.'s with rate $\lambda= (1+\varepsilon)\log C$, with $\varepsilon \log C \geq 4$, then
	\[
	\E[|B_1|] \geq \frac{1}{8} \varepsilon \log C,
	\]
	where $|B_1|$ is the number of items $X_1,\ldots,X_n$ packed in the first bin by Algorithm~\ref{alg:BA} with risk budget $=2/C$.
\end{claim}

\proof{Proof.}
	Let $\ell= \frac{\varepsilon}{4} \log C \geq 1$. Then,
	\begin{align*}
	\Prob(|B_1| \leq \ell ) &= \Prob\left(  |B_1|\leq \ell, X(B_1) > \frac{\varepsilon}{2(1+\varepsilon)} \right) + \Prob\left( |B_1|\leq \ell, X(B_1)\leq \frac{\varepsilon}{2(1+\varepsilon)} \right) \\
	& \leq \Prob\left( \sum_{i=1}^\ell X_i > \frac{\varepsilon}{2(1+\varepsilon)} \right) + \Prob\left( |B_1|\leq \ell, X(B_1)\leq \frac{\varepsilon}{2(1+\varepsilon)}  \right).
	\end{align*}
	We bound each term separately. First, we have
	\[
	\Prob\left( \sum_{i=1}^\ell X_i > \frac{\varepsilon}{2(1+\varepsilon)} \right) \leq \frac{2(1+\varepsilon)}{\varepsilon}\E\left[ \sum_{i=1}^\ell X_i \right] = \frac{2(1+\varepsilon)}{\varepsilon} \ell \frac{1}{(1+\varepsilon)\log C} = \frac{1}{2}.
	\]
	For the other term we have that $|B_1|\leq \ell$, given $X(B_1)\leq \frac{\varepsilon}{2(1+\varepsilon)}$, only if $B_1$ runs out of budget. That is,
	\begin{align*}
	\frac{2}{C} & \leq \sum_{i=1}^{\ell+1} \Prob(X_i > 1-\alpha_i)  \leq (\ell+1) \Prob\left( X_i > 1-\frac{\varepsilon}{2(1+\varepsilon)} \right) \leq \left(  \frac{\varepsilon}{4}\log C + 1 \right) \frac{1}{C^{1+\varepsilon/2}} < \frac{2}{C}
	\end{align*}
	using the assumption $\varepsilon \log C \geq 2$. From here we obtain that $\Prob(|B_1|\leq \ell, X(B_1)\leq \frac{\varepsilon}{2(1+\varepsilon)})=0$. Therefore,
	\[
	\Prob(|B_1|\leq \ell) \leq \frac{1}{2}.
	\]
	Then,
	\[
	\E[|B_1|] \geq \frac{1}{2} \ell = \frac{1}{8} \varepsilon \log C.
	\]
	\hfill\halmos
\endproof

\begin{claim}
	Let $X_1,\ldots, X_m$ be $m$ independent exponential r.v.'s with rate $\mu= \beta \log C$, $\beta \geq 4$, then
	\[
	\E[|B_1|] \geq \frac{1}{8} \beta \log C
	\]
	where $|B_1|$ is the number of items $X_1,\ldots,X_m$ packed in the first in by Algorithm~\ref{alg:BA} with risk budget $= 2/C$.
\end{claim}

\proof{Proof.}
	Let $\ell = \frac{\beta}{4}\log C$. Then,
	\begin{align*}
	\Prob\left( |B_1| \leq  \ell \right) &  =  \Prob(|B_1|\leq \ell , X(B_1)> 1/2) + \Prob(|B_1|\leq \ell, X(B_1) \leq 1/2 ) \\
	& \leq \Prob\left( \sum_{i=1}^\ell X_i > \frac{1}{2}\right) + \Prob\left(  |B_1|\leq \ell, X(B_1) \leq {1}/{2} \right)
	\end{align*}
	Now, given the event $X(B_1)\leq \frac{1}{2}$, the only way that $|B_1|\leq \ell $ is by running out of budget. We have then
	\begin{align*}
	\frac{2}{C} & \leq \sum_{i=1}^{\ell+1} \Prob(X_i > 1- \alpha_i)  \leq (\ell+1 )\Prob(X_i > 1/2) \leq \left(\frac{\beta}{4} \log C +1\right) \frac{1}{C^{\beta/2}} < \frac{1}{C} \tag{$\beta \geq 4$}
	\end{align*}
	which cannot happen. Therefore, $\Prob(|B_1|\leq \ell, X(B_1)\leq 1/2) =0$ and then
	\[
	\Prob(|B_1|\leq \ell) \leq 2\E \left[ \sum_{i=1}^\ell X_i \right] = 2\ell \frac{1}{\beta \log C} = \frac{1}{2}.
	\]
	Therefore,
	\[
	\E[|B_1|] \geq \frac{1}{2} \ell = \frac{\beta}{8}\log C.
	\]
	\hfill\halmos
\endproof

\hfill\halmos
\endproof


\begin{claim}
	The cost of $\mathcal{P}$ is $\cost(\mathcal{P}) \leq 3 \E[N_\mathcal{P}] \leq 48 n_1 \left( \frac{k}{\beta \log C} + \frac{1}{\varepsilon \log C}  \right)$
\end{claim}

\proof{Proof.} Putting all the results together we obtain
	\begin{align*}
	\E[N_\mathcal{P}] & = \E[N_\mathcal{P}^1] + \E[N_\mathcal{P}^2] \\
	& \leq \frac{2n_1}{\E[|B_1^\lambda|]} + \frac{2n_2}{\E[|B_1^\mu|]} \tag{Proposition~\ref{lem:bound_iid_walds}}\\
	& \leq 16 \frac{n_1}{\varepsilon \log C} + 16\frac{n_2}{\beta \log C} \\
	& = 16 n_1 \left( \frac{1}{\varepsilon \log C}  + \frac{k}{\beta \log C} \right).
	\end{align*}
	\hfill\halmos
\endproof


Here we present the proof of Lemma~\ref{lem:bound_iid_walds}.

\proof{Proof of Lemma~\ref{lem:bound_iid_walds}.} We use the following fictitious experiment. Consider $n$ independent copies of the random variables $X_1,\ldots,X_n$ and run Algorithm~\ref{alg:BA} until its first bin is closed or the sequence fits entirely on the first bin. We denote by $\widetilde{B}_i$ the items packed in the first bin in the $i$-th trial of this experiment. The process $|\widetilde{B}_1|,\ldots,|\widetilde{B}_n|$ is i.i.d..
	
	We have the following identities:
	\begin{align*}
	|B_1| & = |\widetilde{B}_1| \\
	|B_2| & = \min\{ n- |B_1|, \widetilde{B}_2 \}\\
	& \,\,\, \vdots \\
	|B_n| & = \min\{  n- |B_1| - \cdots - |B_{n-1}|, \widetilde{B}_n \}.
	\end{align*}
	Observe that $N_\ALG = \min\left\{ k : \sum_{i=1}^k |B_i| = n  \right\}$ is a stopping time for $|B_1|,\ldots,|B_n|$ so also is a stopping time for $|\widetilde{B}_1|,\ldots,|\widetilde{B}_n|$. By Wald's equation (see Theorem~\ref{thm:Waldseq} below) we have
	\[
	\E[N_\ALG] \E[|\widetilde{B}_1|] = \E\left[ \sum_{i=1}^{N_\ALG} |\widetilde{B}_i|\right].
	\]
	Additionally, we have $\E[|\widetilde{B}_1|] = \E[|B_1|]$ by construction. Now, until time $N_\ALG-1$ we must have $|B_1| = |\widetilde{B}_1|, \ldots, |B_{N_\ALG-1}| =|\widetilde{B}_{N_\ALG-1}|$, all of these values at least $1$. Then,
	\[
	\sum_{i=1}^{N_\ALG-1} |\widetilde{B}_i| = \sum_{i=1}^{N_\ALG-1} |B_i| \leq n -1.
	\]
	Therefore,
	\[
	\E[N_\ALG] \E[|B_1|] = \E\left[  \sum_{i=1}^{N_\ALG-1} |B_i| + |\widetilde{B}_{N_\ALG}| \right] \leq 2n-1,
	\]
	which concludes the proof.\hfill\halmos
\endproof


\subsubsection{Wald's Equation}

\begin{theorem}[Wald's equation]\label{thm:Waldseq}
	If $X_1,X_2,\ldots$ are i.i.d. random variables with finite mean and $N$ is a stopping time with $\E[N]<\infty$, then
	\[
	\E \left[\sum_{n=1}^N X_n \right] = \E[N] \E[X_1].
	\]
\end{theorem}

Proof can be found in~\cite{ross1996stochastic}.

\section{$\# \Pol$-Hardness of Computing Minimum Cost of the Optimal Policy}\label{sec:hardness_offline}

In this section, we provide the proof of Theorem~\ref{thm:main_hardness_offline}, i.e., it is $\#\Pol$-hard to compute $\min_{\mathcal{P}} \cost(\mathcal{P})$. We proceed as follows. We consider \emph{symmetric} logic formulas, that is, $\phi(\x)=\phi(\overline{\x})$ for any $\x$, and we show that the problem $\#\textsc{Sym-4Sat}$---the problem of counting satisfying assignment of symmetric formulas in 4CNF---is $\#\Pol$-hard. Recall that a formula is in \emph{conjunctive normal form} (CNF) if it is a conjunction of one or more clauses. When the clauses have $k$ literals, we say that the formula is in $k$CNF.

Then, we provide a polynomial time reduction from symmetric formulas $\phi$ in 4CNF into instances of the stochastic bin packing problem such that
$$\min_{\mathcal{P}} \cost (\mathcal{P}) = \frac{5}{2} - \frac{2}{2^n} - \frac{s_\phi}{2^{2n}},$$
where $s_\phi$ denotes the satisfying assignments of $\phi$, i.e., $s_\phi = |\{  \x = (x_1,\ldots,x_n)  : \phi(\x) =1 \}|$.

We now proceed to show the hardness of $\#\textsc{Sym-4Sat}$.

\begin{theorem}\label{thm:sym_sat_hard}
	$\#\textsc{Sym-4Sat}$ is $\#\Pol$-hard.
\end{theorem}

\proof{Proof.} We show a reduction from the $\#\Pol$-hard problem $\#\textsc{2Sat}$~\cite{valiant1979complexity}. We symmetrize a formula $\phi$ by extending the assignments $\x=(x_1,\ldots,x_n)$ in one variable, namely $x_0$. Let $\phi= \bigwedge_{j=1}^m C_j$ be a formula in 2CNF, i.e., each clause has the form $C_j=(\ell_{1,j} \vee \ell_{2,j})$ where $\ell_{i,j} \in \{ {x}_1,\ldots, x_{n}, \overline{x}_{1},\ldots, \overline{x}_n \}$ for $i=1,2$. Consider the formula
\[
\widetilde{\phi}(x_0,\x) = (\overline{x}_0 \wedge \phi(\x)) \vee ({x}_0\wedge \phi(\overline{\x}) ).
\]
Note that $\widetilde{\phi}$ is symmetric but not yet in CNF. Also, note that $s_{\widetilde{\phi}}  = |\{ (x_0,\x) : \widetilde{\phi}(x_0,\x)= 1 \}|  = 2 s_\phi$.
It remains to show that we can transform $\widetilde{\phi}$ into a symmetric 4CNF without altering the number of solutions. Using De Morgan's law, we rewrite the formula $\widetilde{\phi}$ as
\begin{align*}
	\widetilde{\phi}(x_0,\x) & = (x_0 \vee \phi(\x)) \wedge (\overline{x}_0 \vee \phi(\overline{\x})) \wedge (\phi(\x) \vee \phi(\overline{\x})) \\
	& = \left( x_0 \vee \bigwedge_{j=1}^m C_j(\x) \right) \wedge \left( \overline{x}_0 \vee \bigwedge_{j=1}^m C_j(\overline{\x}) \right) \wedge \left( \bigwedge_{j=1}^m C_j(\x) \vee \bigwedge_{j=1}^m C_j(\overline{\x}) \right) \\
	& = \underbrace{\left( \bigwedge_{j=1}^m (x_0 \vee C_j(\x)) \right)}_{\text{(I)}} \wedge \underbrace{\left( \bigwedge_{j=1}^m (\overline{x}_0 \vee C_j(\overline{\x})) \right)}_{\text{(II)}} \wedge \underbrace{\left( \bigwedge_{j,k=1}^m C_j(\x) \vee C_k(\overline{\x})\right)}_{\text{(III)}}.
\end{align*}
The first two terms (I) and (II) are clearly 3CNF. We extend them into 4CNF by repeating the variable $x_0$ or $\overline{x}_0$ accordingly. The last term (III) is already in 4CNF. We remove clauses that contain pairs $\ell \vee \overline{\ell}$ since they are trivially satisfied.\hfill\halmos
\endproof

Note that in the formula $\widetilde{\phi}$, each variable appears at most twice in each clause, either as $x\vee x$ or $\overline{x}\vee \overline{x}$. This is going to be utilized in the next proof.


Before going to the proof, and for the sake of explanation, we change the capacity of the bins to a capacity $B>1$, to be defined later. This can be easily adjusted to our setting with capacity $1$ by scaling down items sizes by the amount $B$. As we are going to see, it is clearer to introduce items $>1$ than their fractional rescaled version.

\subsection{Reduction $\#\textsc{Sym-4Sat}$ to Stochastic Bin Packing}

The reduction is similar to the reduction from \textsc{Partition-Problem} to \textsc{Bin-Packing}. See~\cite{sipser1997introduction} for an example. A similar reduction is used in~\cite{dean2005adaptivity} in the context of multidimensional stochastic knapsack. Given a symmetric 4CNF $\phi$ with variables $x_1,\ldots,x_n$ and clauses $C_1,\ldots,C_m$, we are going to construct a list of nonnegative random variables, that will correspond to the input of the stochastic bin packing problem, that can be packed into two bins whenever the outcomes of these random variable satisfy the 4CNF formula while in the opposite case, the number of bins requires is at least three.

We are going to utilize numbers with at most $4n+m$ digits in base 10. For convenience we assume that all number have exactly $4n+m$ digits by filling the unused corresponding significant digits with $0$. Given a number with $4n+m$ digits in base 10, we split its representation into four blocks. The first block corresponds to the digits it positions $10^{i} \cdot 10^n \cdot 10^{2n} \cdot 10^m $ for $i=0,\ldots,n-1$. We refer to the first block as the \emph{variable block} and their intra significant digits as the \emph{variables digits}. For instance, by variable digit $x_i$ we refer to the digit in location $10^{n-i} \cdot 10^n \cdot 10^{2n} \cdot 10^m$.

The second block corresponds to digits at positions $10^i \cdot 10^{2n} \cdot 10^m$ for $i=0,\ldots,n-1$. We refer to the second block as the \emph{mirror variable block} and its intra digits as \emph{mirror variable digits}. By \emph{mirror digit} $x_i$ we refer to the digit in location $10^{n-i}\cdot 10^{2n} \cdot 10^m$.

The third block corresponds to digits in positions $10^{k}\cdot 10^m$ for $k=0,\ldots,2n-1$. We refer to this block as the \emph{equivalence block} and we split its digits into pairs of digits that we refer as \emph{equivalence digits}. The positive equivalent digit $x_i$ refers to the digit in position $10^{2n-2i+1} \cdot 10^m$ while the negative equivalence digit $x_i$ refers to the digit in position $10^{2n-2i}\cdot 10^m$.

The last block corresponds to the digits at position $10^j$ for $j=0,\ldots,m-1$. We refer to the this block as the \emph{clauses block} and its intra digits as \emph{clauses digits}. By clause digit $C_j$ we refer to the digit located at position $10^{m-j}$.

We set the capacity of the bins to be the number in base 10
\[
B=\underbrace{11 \cdots 11}_{n \text{ 1's}} \mid \underbrace{11 \cdots 11}_{n\text{ 1's}}  \mid \underbrace{11 \cdots 11}_{2n \text{ 1's}} \mid \underbrace{44 \cdots 44}_{m\text{ 4's}}.
\]
We purposely separated the significant digits using the vertical bars ``$\,\mid\,$'' into the four aforementioned blocks.

For a formula $\phi$ in 4CNF we now present the reduction. For each variable $i=1,\ldots,n$ we construct the following four numbers $a_i,b_i,c_i$ and $d_i$ in base 10. The first two number are
\begin{align*}
	a_i &=1\underbrace{0\cdots 0}_{n-i} \mid \underbrace{0 0 \cdots 0 0}_n \mid \underbrace{00 \cdots 00}_{2i-2}01\underbrace{00\cdots 00}_{2n-2i} \mid \underbrace{0 0 \cdots 00}_m \\
	b_i &= 1\underbrace{0\cdots 0}_{n-i} \mid \underbrace{0 0 \cdots 0 0}_n \mid \underbrace{00 \cdots 00}_{2i-2}10\underbrace{00\cdots 00}_{2n-2i} \mid \underbrace{0 0 \cdots 00}_m .
\end{align*}
Number $a_i$ and $b_i$ have $4n-i+m+1$ digits. Both of them have a common digit $1$ in variable digit $x_i$. Moreover, $a_i$ has another digit $1$ in negative equivalent digit $x_i$; and $b_i$ has a digit $1$ in positive equivalence digit $x_i$.

The following two numbers are

\begin{align*}
	c_i &=\underbrace{\cdots }_{\text{no digits}}\mid \underbrace{\cdots }_{\text{no digits}} 1\underbrace{0\cdots 0}_{n-i} \mid \underbrace{00 \cdots 00}_{2i-2}10\underbrace{00\cdots 00}_{2n-2i} \mid \underbrace{0 0 \cdots 0c_i^{k_i} 00 \cdots  00}_m \\
	d_i &= \underbrace{\cdots }_{\text{no digits}} \mid \underbrace{\cdots }_{\text{no digits}} 1\underbrace{0\cdots 0}_{n-i} \mid \underbrace{00 \cdots 00}_{2i-2}01\underbrace{00\cdots 00}_{2n-2i} \mid \underbrace{0 0 \cdots 000 d_i^{k_i} \cdots 00}_m .
\end{align*}
Number $c_i$ and $d_i$ have $3n-i+m+1$ digits, and note that we keep the separation between the blocks to emphasize where the nonzero digits appear. In other words, the variable block is completely missing from $c_i$ and $d_i$. Now, both number $c_i$ and $d_i$ have a common digit $1$ in mirror variable digit $x_i$ (mirror numbers $x_{i'}$ with $i'<i$ are also missing). The number $c_i$ has digit $c_i^{k_i}\in \{1,2\}$ in all clauses digits $C_{k_i}$ where literal $x_i$ appears; the number $d_i$ has digit $d_i^{k_i} \in \{1,2\}$ in all clauses digits where literal $\overline{x}_i$ appears. The number $c_i$ has a digit $1$ in positive equivalence digit $x_i$; and $d_i$ has a digit $1$ in negative equivalence digit $x_i$.

For each clause $C_j$, $j=1,\ldots,m$, we introduce three numbers $f_j,g_j$ and $h_j$ that are going to serve as slacks:
\[
f_j = g_j = h_j = \underbrace{\cdots }_{\text{no digits}} \mid \underbrace{\cdots }_{\text{no digits}} \mid    \underbrace{\cdots }_{\text{no digits}}   \mid\underbrace{\cdots }_{\text{no digits}} 1\underbrace{0\cdots 0}_{m-j}.
\]
The three numbers have a unique digit $1$ at clause digit $C_j$ and they completely miss the variable, mirror variable and equivalence blocks.

Finally, we introduce a number needed for technical reasons:
\[
h = \underbrace{\cdots }_{\text{no digits}} \mid \underbrace{\cdots }_{\text{no digits}} \mid \underbrace{11 \cdots 11}_{m} \mid \underbrace{00 \cdots 00}_{2n}.
\]
A pictorial construction of the numbers appears in Figure~\ref{fig:reduction_table_0}.

\begin{figure}[h!]
	\centering
	\scriptsize
	\resizebox{\columnwidth}{!}{
		\begin{tabular}{c|ccccc|ccccc|ccccccc|ccccc}
			& $x_1$ & $x_2$ & $x_3$ & $\cdots$ & $x_n$ & $y_1$ & $y_2$ & $y_3$ & $\cdots$ & $y_n$ & $x_1^+$ & $x_1^-$ & $x_2^+$ & $x_2^-$ & $\cdots$ & $x_n^+$ & $x_n^-$ & $C_1$ & $C_2$ & $C_3$ &$\cdots$ & $C_m$
			\\
			\hline
			$a_1$ & 1 & 0 & 0 & $\cdots$ & 0 & 0 & 0 & 0 & $\cdots$ & 0 &  $0$ & $1$ & 0 & 0 & $\cdots$ & 0 & 0 & 0 & 0 & 0 &$\cdots$ & 0
			\\
			$b_1$ & 1 & 0 & 0 & $\cdots$ & 0 & 0 & 0 & 0 & $\cdots$ & 0 & $1$ & $0$ & 0 & 0 & $\cdots$ & 0 & 0 & 0 & 0 & 0 & $\cdots$ & 0
			\\
			$a_2$ & ~ & 1 & 0 & $\cdots$ & 0 & 0 & 0 & 0 & $\cdots$ & 0 & 0 & 0 & 0 & 1 & $\cdots$ & 0 & 0 & 0 & 0 & 0 & $\cdots$ & 0
			\\
			$b_2$ & ~ & 1 & 0 & $\cdots$ & 0 & 0 & 0 & 0 & $\cdots$ & 0 & 0 & 0 & 1 & 0 & $\cdots$ & 0 & 0  & 0 & 0 & 0 &$\cdots$ & 0
			\\
			$a_3$ & ~ & ~ & 1 & $\cdots$ & 0 & 0 & 0 & 0 & $\cdots$ & 0 & 0 & 0 & 0 & 0 & $\cdots$ & 0 & 0  & 0 & 0 & 0 &$\cdots$ & 0
			\\
			$b_3$ & ~ & ~ & 1 & $\cdots$ & 0 & 0 & 0 & 0 & $\cdots$ & 0 & 0 & 0 & 0 & 0 & $\cdots$ & 0 & 0  & 0 & 0 & 0 &$\cdots$ & 0
			\\
			$\vdots$ & ~ & ~ & ~ & $\ddots$ & ~ & ~ & ~ & ~ & $\vdots$ & ~ & ~ & ~ & ~ & ~ & $\vdots$ & ~ & ~  & ~ & ~ & ~ &$\vdots$ & ~
			\\
			$a_n$ & ~ & ~ & ~ & ~ & 1 & 0 & 0 & 0 & $\cdots$ & 0 & 0 & 0 & 0 & 0 & $\cdots$ & 0 & 1  & 0 & 0 & 0 &$\cdots$ & 0
			\\
			$b_n$ & ~ & ~ & ~ & ~ & 1 & 0 & 0 & 0 & $\cdots$ & 0 & 0 & 0 & 0 & 0 & $\cdots$ & 1 & 0  & 0 & 0& 0 &$\cdots$ & 0
			\\
			\hline
			$c_1$ & ~ & ~ & ~ & ~ & ~ & 1 & 0 & 0 & $\cdots$ & 0 & 1 & 0 & 0 & 0 & $\cdots$ & 0 & 0 & 2 & 0 & 0 &$\cdots$ & 0
			\\
			$d_1$ & ~ & ~ & ~ & ~ & ~ & 1 & 0 & 0 & $\cdots$ & 0 & 0 & 1 & 0 & 0 & $\cdots$ & 0 & 0  & 0 & 1 & 1 &$\cdots$ & 0
			\\
			$c_2$ & ~ & ~ & ~ & ~ & ~ & ~ & 1 & 0 & $\cdots$ & 0 & 0 & 0 & 1 & 0 & $\cdots$ & 0 & 0 & 1 & 1 & 0 &$\cdots$ & 0
			\\
			$d_2$ & ~ & ~ & ~ & ~ & ~ & ~ & 1 & 0 & $\cdots$ & 0 & 0 & 0 & 0 & 1 & $\cdots$ & 0 & 0  & 0 & 0 & 0 &$\cdots$ & 0
			\\
			$c_3$ & ~ & ~ & ~ & ~ & ~ & ~ & ~ & 1 & $\cdots$ & 0 & 0 & 0 & 0 & 0 & $\cdots$ & 0 & 0 & 0 & 0 & 1 &$\cdots$ & 0
			\\
			$d_3$ & ~ & ~ & ~ & ~ & ~ & ~ & ~ & 1 & $\cdots$ & 0 & 0 & 0 & 0 & 0 & $\cdots$ & 0 & 0 & 1 & 2 & 0 &$\cdots$ & 0
			\\
			$\vdots$ & ~ & ~ & ~ & ~ & ~ & ~ & ~ & ~ & $\ddots$ & ~ & ~ & ~ & ~ & ~ & $\vdots$ & ~ & ~  & ~ & ~ & ~ &$\vdots$ & ~
			\\
			$c_n$ & ~ & ~ & ~ & ~ & ~ & ~ & ~ & ~ & ~ & 1 & 0 & 0 & 0 & 0 & $\cdots$ & 1 & 0 & 0 & 0 & 1 &$\cdots$ & 0
			\\
			$d_n$ & ~ & ~ & ~ & ~ & ~ & ~ & ~ & ~ & ~ & 1 & 0 & 0 & 0 & 0 & $\cdots$ & 0 & 1  & 0 & 0 & 0 &$\cdots$ & 2
			\\
			\hline
			$f_1$ & ~ & ~ & ~ & ~ & ~ & ~ & ~ & ~ & ~ & ~ & ~ & ~ & ~ & ~ & ~ & ~ & ~ & 1 & 0 & 0 &$\cdots$ & 0
			\\
			$g_1$ & ~ & ~ & ~ & ~ & ~ & ~ & ~ & ~ & ~ & ~ & ~ & ~ & ~ & ~ & ~ & ~ & ~ & 1 & 0 & 0 &$\cdots$ & 0
			\\
			$h_1$ & ~ & ~ & ~ & ~ & ~ & ~ & ~ & ~ & ~ & ~ & ~ & ~ & ~ & ~ & ~ & ~ & ~ & 1 & 0 & 0 &$\cdots$ & 0
			\\
			$f_2$  & ~ & ~ & ~ & ~ & ~ & ~ & ~ & ~ & ~ & ~ & ~ & ~ & ~ & ~ & ~ & ~ & ~ & ~ & 1 & 0 &$\cdots$ & 0
			\\
			$g_2$  & ~ & ~ & ~ & ~ & ~ & ~ & ~ & ~ & ~ & ~ & ~ & ~ & ~ & ~ & ~ & ~ & ~ & ~ & 1 & 0 &$\cdots$ & 0
			\\
			$h_2$  & ~ & ~ & ~ & ~ & ~ & ~ & ~ & ~ & ~ & ~ & ~ & ~ & ~ & ~ & ~ & ~ & ~ & ~ & 1 & 0 &$\cdots$ & 0
			\\
			$\vdots $ & ~ & ~ & ~ & ~ & ~ & ~ & ~ & ~ & ~ & ~ & ~ & ~ & ~ &~  & ~ & ~ & ~ & ~ & ~ &~ & $\ddots$ & ~
			\\
			$f_m$ & ~ & ~ & ~ & ~ & ~ & ~ & ~ & ~ & ~ & ~ & ~ & ~ & ~ &~ & ~ & ~ & ~ & ~ & ~ & ~ & ~ & 1
			\\
			$g_m$ & ~ & ~ & ~ & ~ & ~ & ~ & ~ & ~ & ~ & ~ & ~ & ~ & ~ &~ & ~ & ~ & ~ & ~ & ~ & ~ & ~ & 1
			\\
			$h_m$ & ~ & ~ & ~ & ~ & ~ & ~ & ~ & ~ & ~ & ~ & ~ & ~ & ~ &~ & ~ & ~ & ~ & ~ & ~ & ~ & ~ & 1
			\\
			\hline
			$h$ & ~ & ~ & ~ & ~ & ~ & ~ & ~ & ~ & ~ & ~ & ~ & ~ & ~ & ~ & ~ & ~ & ~ & 1 & 1 & 1 & $\cdots$ & 1
			\\
			\hline
			$B$ & 1 & 1 & 1 & $\cdots$ & 1 & 1 & 1 & 1 & $\cdots$ & 1  & 1 & 1 & 1 & 1 & $\cdots$ & 1 & 1 & 4 & 4 & 4 & $\cdots$ & 4
			\\
			\hline
	\end{tabular}}
	\caption{Construction of numbers via $\phi$ 4CNF. The table is purposely divided into four blocks representing the digit blocks defined at the beginning of the subsection. The instance showed corresponds partially to the 4CNF $\phi(x_1,\ldots,x_n) =  (x_1 \vee x_1 \vee x_2 \vee \overline{x}_3 )\wedge (\overline{x}_1 \vee x_2 \vee \overline{x}_3\vee \overline{x}_3) \wedge (\overline{x}_1 \vee \vee x_3 \vee x_4 \vee x_n) \wedge \cdots \wedge (x_{n-2}\vee x_{n-1}\vee \overline{x}_n \vee \overline{x}_n)$, where clauses are $C_1,C_2,C_3,\ldots,C_m$ in the order they are displayed.}\label{fig:reduction_table_0}
\end{figure}


Given these number, we construct an instance of the stochastic bin packing problem as follows. We define the following independent random variables
\[
X_i = \begin{cases}
	a_i & \text{w.p. } 1/2 \\
	b_i & \text{w.p. } 1/2
\end{cases}
\quad\text{ and }
X_i' = \begin{cases}
	a_i & \text{w.p. } 1/2 \\
	b_i & \text{w.p. } 1/2
\end{cases}.
\]
Now the instance is given by the sequence
\[
\mathcal{L}_\phi = (X_1, X_1', X_2, X_2', X_3, X_3', \ldots, X_n,X_n', c_1,d_1,c_2,d_2,\ldots,c_n,d_n,h_1,g_1,\ldots,h_m,g_m,h).
\]
Intuitively, we aim to simulate the random evaluation of $\phi$ when the values of $x_1,\ldots,x_n$ are chosen uniformly and independently of each other. Note that in this case,
\[
\Prob_{\x \in_R \{ 0,1 \}^n }( \phi(\x)  =1 ) = \frac{s_\phi}{2^n}.
\]
Note that in instance $\mathcal{L}_\phi$, any policy incurs in a cost of at least $2$ since $X_1$ and $X_1'$ if packed together incur in an expected cost of at least $2$. This is because $C>2$. As we did in the main body of the article, we can assume that the policies are deterministic. Moreover, we can assume that policies never break a bin since the expected cost of breaking a bin is always greater than $1$. This is because the probability of overflow always is either $0, 1/2$ or $1$ by construction of $\mathcal{L}_\phi$.

\begin{theorem}\label{thm:reduction_thm}
	For the instance $\mathcal{L}_\phi$ we have
	\[
	\min_{\mathcal{P}} \cost(\mathcal{P}) = \frac{5}{2} - \frac{2}{2^n} - \frac{1}{2^n} \Prob_{ \x\in \{0,1\}^n }(\phi(\x)=1).
	\]
\end{theorem}


\begin{proposition}\label{prop:reduction_structure_1}
	For the instance $\mathcal{L}_\phi$ there is a policy with $\cost(\mathcal{P})\leq 3$. Moreover, for any of such policies, there is a policy with same cost or better that packs items $X_1,\ldots,X_n$ into bin $1$ and items $X_1',\ldots,X_n'$ into bin $2$.
\end{proposition}

\proof{Proof.} Consider the policy that packs item $X_1$ into bin $1$; item $X_1'$ into bin $2$; the rest of the items into bin $3$. This policy is valid since
\[
\sum_{i=2}^n (X_i + X_i') + \sum_{i=1}^n (c_i + d_i) + \sum_{j=1}^m (h_j + g_j) + h \leq B.
\]
This proves the first part of the proposition.

For the second part, consider any policy $\mathcal{P}$ with cost at most $3$ that does not break any bin. We can assume, without loss of generality, that $X_1$ is packed into bin $1$ and $X_1'$ is packed into bin $2$. Now, starting at the root of the policy tree, find the first node $u$ where $X_i$ is not packed in bin $1$, say bin $j\geq 2$. Note that up to that point, items $X_{k}'$ must have been packed in a different bin than bin $1$. If both children of $u$ are packed into bin $1$, that is, $X_i'$ is packed into bin $1$, then exchange the packing rule in node $u$ by packing $X_i$ into bin $1$ and in its children to pack item $X_i'$ into bin $j$. This does not change the cost since $X_i$ and $X_i'$ are identically distributed. Suppose now that some of the children of node $u$ packs item $X_i'$ into bin $j'\neq 1$, say children $v$. At this point, there is enough space in bin $1$ to receive $X_i'$ since $X_i$ was packed into bin $j\geq 2$. In the subtree rooted at $v$, mark all nodes that pack their corresponding item into bin $1$. Exchange the packing rule from these node to pack their items into bin $j'$ and change the policy to pack item $X_i'$ in node $v$ from bin $j'$ to bin $1$. This does not increase the cost of the policy. With this, we can modify the policy to pack $X_i'$ in both children of node $u$ into bin $1$ without increasing the expected cost. Now, like in the previous case, we can pack item $X_i$ into bin $1$ and item $X_i'$ into bin $j\geq 2$. We can repeat this procedure for all $i$ and at the end of this, we can ensure that all $X_1,\ldots,X_n$ are packed into bin $1$ and the cost of the policy does not increase. With a similar argument, we show that items $X_1',\ldots,X_n'$ can be packed into bin $2$ without increasing the cost of the policy. From the root, find the first node $u'$ where $X_i'$ is not packed into bin $2$, say bin $j$. Note that $j$ cannot be $1$ since we already have that $X_i$ has been packed into bin $1$ and both random variables share the same variable digit $x_i$. Then, $j\geq 3$. In the subtree rooted at $u'$ mark all items that are packed into bin $2$. Repack those items into bin $j$ and pack $X_i'$ into bin $2$. This does not increase the cost of the policy since we have enough space in bin $j$ to receive any item in there if needed. Repeating this procedures for all $i$ gives us the desired result.\hfill\halmos
\endproof

Consider the following three events: Let
\[
\mathcal{E} = \{ \forall i=1,\ldots, n : X_i \neq X_i'  \},
\]
then $\Prob(\mathcal{E}) = \frac{1}{2^n}$. Let
\[
\mathcal{E}_{a} = \{ \exists i =1,\ldots,n: \forall k < i, X_k \neq X_k', X_i=X_i' = a_i  \}
\quad\text{and}\quad
\mathcal{E}_{b} = \{ \exists i =1,\ldots,n: \forall k < i, X_k \neq X_k', X_i=X_i' = b_i  \}.
\]
Note that $\mathcal{E}_a$ and $\mathcal{E}_b$ are disjoint and $\mathcal{E}_a\cup \mathcal{E}_b = \overline{\mathcal{E}}$. Intuitively, $\mathcal{E}$ is the good event where the variables $X_i$ and $X_i'$ model opposite values in $\{a_i,b_i\}$.

We denote by $\cost(\mathcal{P}\mid \mathcal{A})$ the conditional expected cost of the policy $\mathcal{P}$ on the event $\mathcal{A}$.

\begin{proposition}\label{prop:hardness_bad_event_lower_bound}
	For any policy $\mathcal{P}$ we have $\cost(\mathcal{P}\mid \mathcal{E}_b) \geq 3$.
\end{proposition}

\proof{Proof.} Note that given $\mathcal{E}_b$ we have (assuming the collision occurs at $i$, $X_i=X_i'=b_i$) we have
\begin{align*}
	\sum_{i=1}^n (X_i + X_i' + c_i + d_i) + \sum_{j=1}^m (h_j+g_j) + h &\geq 22 \cdots 22  \mid 22 \cdots 22 \mid  \underbrace{22 \cdots 22}_{2i-2}311 \cdots 11 \mid  88 \cdots 88  > 2 B
\end{align*}
where $B$ is the bin's capacity. Therefore, since we are assuming that $\mathcal{P}$ does not break bins, then the policy must have packed all the items into at least $3$ bins which concludes the proof.\hfill\halmos
\endproof

\begin{proposition}\label{prop:hardness_good_event_lower_bound}
	For any policy $\mathcal{P}$ that packs items $X_1,\ldots,X_n$ into bin $1$ and items $X_1',\ldots,X_n'$ into bin $2$ and does not break any bin, we have $\cost(\mathcal{P}\mid \mathcal{E}) \geq  3 - \Prob_{\x\in_R \{ 0,1 \}^n }( \phi(\x) =1) $.
\end{proposition}

\proof{Proof.} Consider the following random assignment $\X = (x_1,\ldots,x_n)$:
\[
x_i = \begin{cases}
	1 & \text{if } X_i= a_i \\
	0 & \text{if } X_i= b_i
\end{cases}.
\]
Note that $\X$ is uniformly distributed over $\{0,1\}^n$. Then, $\cost(\mathcal{P}\mid \mathcal{E}, \phi(\X)=0)\geq 3$. Indeed, if only $2$ bins have been used after all items have been packed, this forces item $c_i$ to be placed in bin $1$ if $X_i=a_i$ and in bin $2$ otherwise, while $d_i$ is packed in the opposite bin to $c_i$. Since $\phi(\X)=0= \phi(\overline{\X})$ by symmetry of $\phi$, then, after packing items $c_1,d_1,\ldots,c_n,d_n$ but before packing items $f_1,g_1,h_1,\ldots,f_m,g_m,h_m$, there must be a $C_j$ digit in the utilization of bin 1 that is $0$ and a $C_{j'}$ digit in the utilization of bin 2 that is also $0$, $j\neq j'$. In particular, this implies that the usage of bin $2$ at this point has a $C_j$-digit of $4$. After packing items $f_1,g_1,h_1,\ldots,f_m,g_m,h_m$, one of the bins must have a $4$ in its $C_1$-digit of usage, say bin $1$. Therefore, bin $2$ has a usage with $C_1$-digit $3$ and item $h$ cannot be packed into bin $1$. Since $\mathcal{E}$ is given, we have $\sum_{i=1}^n (X_i+X_i')+ \sum_{i=1}^n (c_i+ d_i) + \sum_{j=1}^m (f_j+g_j+h_j)+ h = 2B$, then $h$ can be packed into bin $2$ only if all $C_1,\ldots,C_m$-digits are $3$. However, this contradicts the fact that $C_{j}$-digit in bin $2$ is $4$. Similarly, if the bin $2$ has usage with $C_1$-digit of $4$ after packing items $f_1,g_1,h_1,\ldots,f_m,g_m,h_m$, we can obtain the same contradiction. Therefore, $\cost(\mathcal{P}\mid \mathcal{E} ,\phi(\X)= 0)\geq 3$.

Now,
\begin{align*}
	\cost(\mathcal{P}\mid \mathcal{E}) & \geq 3 \Prob(\phi(\X)=0\mid \mathcal{E}) + 2 \Prob(\phi(\X)=1 \mid \mathcal{E})  = 3 - \Prob(\phi(\X) = 1\mid \mathcal{E})  = 3 - \Prob_{\x\in_R\{ 0,1 \}^n}(\phi(\x) = 1).
\end{align*}
\hfill\halmos
\endproof

\begin{lemma}\label{lem:hardness_lower_bound_policy}
	For any policy $\mathcal{P}$, $\cost(\mathcal{P}) \geq  \frac{5}{2} - \frac{3}{2^n} - \frac{1}{2^n} \Prob_{\x \in_R \{ 0,1 \}^n }( \phi(\x)=1)$.
\end{lemma}

\proof{Proof.} We can assume that the policy $\mathcal{P}$ packs items $X_1,\ldots,X_n$ into bin $1$ while items $X_1',\ldots,X_n'$ into bin $2$ (see Proposition~\ref{prop:reduction_structure_1}). Then, utilizing Propositions~\ref{prop:hardness_bad_event_lower_bound} and~\ref{prop:hardness_good_event_lower_bound} we obtain
\begin{align*}
	\cost(\mathcal{P}) & \geq \left( 3 - \Prob_{\x\in_R\{ 0,1 \}^n}(\phi(\x) = 1) \right)\Prob(\mathcal{E}) + 3\Prob(\mathcal{E}_b) + 2\Prob(\mathcal{E}_a) \\
	& = \frac{5}{2}\left( 1- \frac{1}{2^n} \right) + \frac{1}{2^n}\left( 3 - \Prob_{\x\in_R\{ 0,1 \}^n}(\phi(\x) = 1) \right) \\
	& = \frac{5}{2} - \frac{2}{2^n}  - \frac{1}{2^n} \Prob_{\x\in_R\{ 0,1 \}^n}(\phi(\x) = 1).
\end{align*}
\hfill\halmos
\endproof

\begin{lemma}\label{lem:hardness_optimal_policy}
	There is a policy $\mathcal{P}$ such that $
	\cost(\mathcal{P}) \leq \frac{5}{2} - \frac{2}{2^n}  - \frac{1}{2^n} \Prob_{\x\in_R\{ 0,1 \}^n}(\phi(\x) = 1)$.
\end{lemma}

\proof{Proof.} Consider the following policy. Start packing items $X_1,\ldots,X_n$ into bin $1$ and items $X_1',\ldots,X_n'$ into bin $2$. If there is a collision $X_i=X_i'$ and the first of these is $X_i=X_i'=b_i$, then pack the rest of the items into bin $3$; While if the first of these collisions is $X_i=X_i'=a_i$, then continue packing as follows. Pack $c_{i'}$ where $a_{i'}$ is packed and $d_{i'}$ where $b_{i'}$ is packed for $i'<i$. Pack $c_i$ into bin $1$ and $d_i$ into bin $2$. Pack the remaining $c_{i'}$ into bin $1$ and $d_{i'}$ into bin $2$. Since $i$ is the first time there is a collision, bin $2$ has a $0$ in its positive equivalent digit $x_i$ and a $2$ in its negative equivalent digit $x_i$. Therefore it has enough space to receive the items packed after $d_{i'}$ has been packed. Then only $2$ bins are utilized.

If no collision happens, then pack $c_i$ where outcome $a_i$ is packed (bin $1$ if $X_i=a_i$ or bin $2$ if $X_i'=a_i$) and pack $d_i$ in the opposite bin (bin $2$ if $c_i$ is in bin $1$ bin $1$ otherwise). Now, utilize the slack items $f_j,g_j,h_j$ to complete bin $1$ and then bin $2$. Now, for $h$ there are two cases based on the value of $\phi$ on the satisfying assignment $\x$ given by
\[
x_i = \begin{cases}
	1 & \text{if } X_i= a_i \\
	0 & \text{if } X_i= b_i
\end{cases}.
\]
\begin{itemize}[leftmargin=*, itemsep=0em]
	\item If $\phi(\x)=1$, then in bin $2$, each $C_j$ digit of the capacity used is at most $3$. This is because the $C_j$ digit of the capacity used bin $1$ is $4$ by construction. Now, pack $h$ into bin $2$ and finish the packing into $2$ bins.
	
	\item If $\phi(\x)=0$, we can retrace the proof of Proposition~\ref{prop:hardness_good_event_lower_bound} to show that in this case, $h$ cannot fit nor in bin $1$ nor in $2$, therefore forcing a bin $3$.
\end{itemize}

Putting all these case together, we obtain
\begin{align*}
	\cost(\mathcal{P}) &= \cost(\mathcal{P}\mid \mathcal{E})\Prob(\mathcal{E}) + \cost(\mathcal{P}\mid \mathcal{E}_a) \Prob(\mathcal{E}_a) + \cost(\mathcal{P}\mid \mathcal{E}_b) \Prob(\mathcal{E}_b) \\
	& = \left( 3 - \Prob_{\x\in_R\{0,1\}^n }(\phi(\x)=1)  \right) \frac{1}{2^n} + \frac{5}{2}\left( 1 - \frac{1}{2^n} \right) \\
	& = \frac{5}{2} - \frac{2}{2^n} - \frac{1}{2^n} \Prob_{ \x\in \{0,1\}^n }(\phi(\x)=1).
\end{align*}
\hfill\halmos
\endproof

Putting together Lemma~\ref{lem:hardness_lower_bound_policy} and~\ref{lem:hardness_optimal_policy} we obtain the proof of Theorem~\ref{thm:reduction_thm}.

\section{Threshold Policies for I.I.D. Random Variables with Finite Support}\label{sec:threshold_iid}

In this section we discuss the problem of designing a threshold algorithm that incurs in a constant factor loss whenever the input sequence is i.i.d.\ (with common distribution $\mathcal{D}$) with unknown time horizon $n$. A threshold algorithm observes the common random distribution of the incoming streams and computes a number $\alpha\in [0,1]$ such that bins are utilized as long as their usage is at most $\alpha$. If there are no such bins, then a new bin is opened upon an arrival. Note that in the i.i.d.\ setting, for threshold policies of this kind, at most \emph{one} bin is kept active at a time.

We show that among all policies that keep at most one bin active at a time, threshold policies are optimal up to an additive loss of one. This is under the assumption that the common distribution $\mathcal{D}$ has finite support.
\begin{theorem}\label{thm:threshold_policy}
	Let $\mathcal{D}$ be any distribution with finite support in $[0,\infty)$. There exists $\alpha \in [0,1]$ such that the threshold policy $\mathcal{P}_\alpha$ with threshold $\alpha$ satisfies
	\[
	\cost(\mathcal{P}_\alpha) \leq \min_{\substack{\mathcal{P} \text{ has at most}\\ \text{one active bin}}} \cost(\mathcal{P})+1,
	\]
	for any input sequence of i.i.d.\ random variables $X_1,\ldots,X_n$ with common distribution $\mathcal{D}$.
\end{theorem}

Note that policy $\mathcal{P}_\alpha$ is computed only with the information given by the distribution $\mathcal{D}$. An online algorithm that has access to $\mathcal{D}$ computes the threshold $\alpha$ given in Theorem~\ref{thm:threshold_policy} and implements the threshold policy. In the main body of the paper we show that Budgeted Greedy keeps at most one bin active at a time in the i.i.d.\ setting (Lemma~\ref{lem:key_lem_iid}). Therefore, policies that keep at most one bin active at a time are within a constant factor of the optimal offline sequential cost. Using these facts, we conclude that an algorithm implementing a threshold policy incurs an expected cost that is a constant factor of the optimum. This factor is at most $(3+ 2\sqrt{2})$, since we utilized the guarantee given by Budgeted Greedy (Theorem~\ref{thm:iid_main_thm}).

A major downside of Theorem~\ref{thm:threshold_policy} is that it does not give an efficient algorithm to compute the threshold $\alpha$. Indeed, for the proof of Theorem~\ref{thm:threshold_policy}, we utilize the framework of discounted reward Markov processes (see~\cite{puterman2014markov}), which can compute optimal stationary policies in time depending on the size of the state space. For us, the state space is the usage of the active bin, which can be exponentially large in the description of $\mathcal{D}$. Intuitively, since we aim to compute a policy that does not depend on the time horizon $n$ and we only have one bin to use, the best we can do is to repeat the same process over and over.

\subsection{Proof of Theorem~\ref{thm:threshold_policy}}

The proof is divided in a sequence of propositions. We briefly introduce the definitions used in infinite-time horizon discounted Markov decision processes. We later show that the optimal discounted cost induces a monotonic cost vector. From here, a threshold policy can be deduced which is later used to design a finite-time threshold policy.

We assume that the distribution $\mathcal{D}$ has finite support in $[0,1]\cup\{1^+\}$ where the element $1^+$ denotes a fixed upper bound over the values in $[0,1]$ and any value that $\mathcal{D}$ could have taken above $1$ with positive probability is mapped to $1^+$. The state space, denoted $\State$ corresponds to all possible values $\leq 1$ that the bin can take as combinations of number in the support of $\mathcal{D}$ in addition to the special state $1^+$. Note that $\State$ is a finite set.

Fix a discount factor $\gamma\in (0,1)$. In the infinite time-horizon discounted factor framework, a policy corresponds to a sequence of functions (or distributions if randomized) $\Pi=(\pi_1,\pi_2,\ldots)$ that dictates the behavior of the process. That is, $\pi_t: \State \to \{ 0, 1 \}$ is the decision made by the policy at round $t$, where $0$ indicates \emph{keep using the current bin} while $1$ indicates \emph{open a new bin}, all this as a function of the state of the system. If $\pi_t$ is random, then $\pi_t: \State \to \Delta(\{0,1\})$, where $\Delta(\{0,1\})$ is the probability simplex over $\{0,1\}$. We define the \emph{discounted cost} of the policy $\Pi$ at time $t=1,2,\ldots$ by
\[
V^{\Pi}_t(s) = \begin{cases}
	1+ C\Prob(X > 1) + \gamma \E[V_{t+1}^\Pi(X\wedge 1^+)]  & \pi_t(s)=0 \\
	C \Prob(X+s > 1) + \gamma \E[V_{t+1}^\Pi((X+s)\wedge 1^+)] & \pi_t(s) =1
\end{cases}.
\]
Let
\[
c(s,a)= \begin{cases}
	C \cdot \Prob(X + s > 1) & a = 0 \\
	1 + C \cdot \Prob( X > 1)  & a = 1
\end{cases}
\]
and $T(s,0) = (X + s)\wedge 1^+$; $T(s,1)=X \wedge 1^+$. We can write $V^{\Pi}_t(s) = c(s,\pi(s)) + \gamma \E[ V_{t+1}^{\Pi}(T(s,\pi(s) ) ) ]$. 

If $\pi_t$ are randomized then the previous values are replaced by expectations. The goal is to find $\min_{\Pi} V_1^\Pi(0)$. Markov Decision processes theory guarantees that this minimum is also a minimum over the history dependent randomized policies---policies that record previous outcomes. Moreover, the optimal policy for $\min_{\Pi} V_1^\Pi(0)$ is also the optimal policy for $\min_{\Pi} V_1^\Pi(s)$ for any $s\in \State$. The theory also guarantees that deterministic stationary policies are optimal. That is, $\min_{\Pi} V_1^\Pi(s) = \min_{\pi} V_1^{(\pi,\pi,\ldots)}(s)$. From now on, we only consider deterministic policies. By $V_t^\pi$ we refer to $V_t^\Pi$ where $\Pi=(\pi,\pi,\ldots)$. Note that $V^\pi_1(s) = V_2^\pi(s)=  \cdots$ and so we can identify the temporal cost vector $(V_t^\pi(s))_{\substack{s\in \State \\ t=1,2,\ldots}}$ by just the vector $V^\pi=(V^\pi(s))_{s\in \State}$. The optimal vector $V^\pi$ satisfies the Bellman equation $V= \mathcal{T}^\gamma V$, where
\begin{align*}
	\left(\mathcal{T}^\gamma V\right)(s) &= \min\{ 1 + C\Prob(X > 1)  + \gamma \E[ V(X\wedge 1^+)  ], C \Prob(X+s > 1) + \gamma \E[V((X+s)\wedge 1^+) ] \}\\
	& = \min_{a=0,1}\{ c(s,a) + \gamma \E[ V(T(s,a)) ] \} .
\end{align*}
Therefore, $V^\pi$ is the fixed point of the Bellman operator $\mathcal{T}^\gamma$. For a detailed presentation of these results, see Chapter 6 in~\cite{puterman2014markov}.


\begin{proposition}
	Consider the optimal solution $V^\pi$ of the discounted cost problem. Then $V^\pi$ is a monotone function of $s$. That is, $V^\pi(s')\leq V^\pi (s'')$ for any $0\leq s'\leq s''\leq 1^+$ .
\end{proposition}

\proof{Proof.} By contradiction, suppose that for some $s'< s''$ we have $V^\pi(s')> V^\pi(s'')$. Among all such possible pairs $s'<s''$ choose the largest $s'\leq 1$, which exists because $\State$ is finite. Define the new function $\widehat{\pi}_1$ as $\widehat{\pi}_1(s)= \pi (s)$ if $s\neq s'$ and $\widehat{\pi}_1(s')=\pi(s'')$. For $t\geq 2$ we define $\widehat{\pi}_t = \pi$. Now consider the policy $\widehat{\Pi} = (\widehat{\pi}_1,\widehat{\pi}_2,\ldots)$. Then, using the definition of $V_1^{\widehat{\Pi}}$ we can show that $V^{\widehat{\Pi}}(s) = V^\pi(s)$ for $s\neq s'$. Now, let's analyze the case $s=s'$. Observe that
\[
c(s',\pi(s''))\leq c(s'',\pi(s''))
\]
since $s' < s''$ and the function $c(\cdot, a)$ is nondecreasing for any fixed $a$. Then, we have two cases:
\begin{itemize}[leftmargin=*, itemsep=0em]
	\item If $\pi(s'')=0$, then as $s'$ is the largest state where monotonicity does not hold, we have
	\begin{align*}
		V^{\widehat{\Pi}}_1(s') & = c(s',\pi_1(s'))  + \gamma \E[ V_{2}^{\widehat{\Pi}}(T(s', \pi_1(s')))  ] \\
		& = c(s',0) + \gamma \E[V_2^{(\pi,\ldots)}( (X+s')\wedge 1^+ ))] \tag{$\pi_1(s')=\pi(s'')=0$}\\
		& \leq c(s'', 0) + \gamma \E[V^\pi( (X+s')\wedge 1^+ )] \tag{$V_2^{(\pi,\ldots)} = V^\pi$ and monotonicity of $c(\cdot,0)$}\\
		& = c(s'',0) + \gamma \Prob(X=0) V^\pi (s') + \gamma \E[V((X+s')\wedge 1^+) \mid X>0 ] \Prob(X>0) \\
		& \leq  c(s'',0) + \gamma \Prob(X=0) V^\pi (s') + \gamma \E[V((X+s'')\wedge 1^+) \mid X>0 ]\Prob(X>0) \tag{As $s'$ is the largest value where monotonicity does not hold} \\
		&\leq  c(s'',0) + \gamma \Prob(X=0) (V^\pi (s') - V^\pi(s'') )  + \gamma \E[V^\pi((X+s'')\wedge 1^+) ] \\
		& = V^\pi(s'') + \gamma \Prob( X=0 ) (V^\pi(s')  - V^\pi(s'')) \\
		& = \gamma \Prob(X=0) V^\pi(s') + (1- \gamma \Prob(X=0)) V^\pi(s'') \\
		& < \gamma \Prob(X=0) V^\pi(s') + (1- \gamma \Prob(X=0)) V^\pi(s') \tag{Since $V^\pi(s')> V^\pi(s'')$} \\
		& = V^\pi(s').
	\end{align*}
	
	\item Similarly, if $\pi(s'')= 1$, then,
	\begin{align*}
		V^{\widehat{\Pi}}_1(s') & = c(s',\pi_1(s')) + \gamma \E[ V_2^{\widehat{\Pi}} (T(s',\pi_1(s') ))  ]  \\
		& = c(s',1) + \gamma \E[V^{\pi} (X\wedge 1^+ ) ] \\
		& \leq c(s'',1) + \gamma \E[V^\pi( X \wedge 1^+ )] \\
		& = V^\pi(s'') \\
		& < V^\pi(s').
	\end{align*}
\end{itemize}
In any case, $V^{\widehat{\Pi}}(s') < V^\pi(s')$, which contradicts the optimality of $\pi$.\hfill\halmos
\endproof

The following result states that the optimal policy of the discounted cost problem is a threshold policy.

\begin{proposition}
	For the optimal $V^\pi$, there exists $\alpha \in [0,1]$ such that the stationary policy $\widetilde{\pi}(s)=0$ if $s\leq \alpha$; $\widetilde{\pi}(s)=1$ if $s> \alpha$, holds $\pi=\widetilde{\pi}$.
\end{proposition}

\proof{Proof.} Let $E = \{ s\in \State : V^\pi(s) < 1 + C \Prob(X > 1) + \gamma \E[ V^\pi( X\wedge 1^+ ) ] \}$ be the states where the policy $\pi$ decides to utilize the current bin. Note that $0\in E$, hence $\alpha=\sup E$ is well-defined. Also, note that for $s=1^+$ we have
\begin{align*}
	1 + C \Prob(X> 1) + \gamma \E[ V^\pi(X\wedge 1^+)] & \leq C \Prob(X+ 1^+ > 1) + \gamma \E[V^\pi((X+1^+)\wedge 1^+)],
\end{align*}
thus $\alpha < 1^+$. Using the monotonicity of $V^\pi$, we have $E=[0,\alpha]$.

By definition of $\alpha$ we have that for any $s>\alpha$, $V^\pi(s)\geq 1 + C \Prob(X > 1) + \gamma \E[V^\pi (X\wedge 1^+)]$. This immediately implies that for any $s> \alpha$, $V^\pi(s)= 1 + C \Prob(X> 1) + \gamma \E[V^\pi(X\wedge 1^+)]$. Since $V^\pi(\cdot)$ is monotone, then for any $s\leq \alpha$ we have $V^\pi(s)= C \Prob(X+s > 1)+ \gamma \E[V^\pi((X+s)\wedge 1^+)]$. This shows that $\pi$ is indeed $\widetilde{\pi}$.\hfill\halmos
\endproof

Note that in this discounted cost model we did not restrict the possible actions when the usage of the bin goes beyond $1$, i.e, in state $1^+$. The optimality and monotonicity of the optimal value $V^\pi$ shows that the optimal policy never tries to utilize the overflowed bin again and it will always choose to open a new bin.

We now return, to our model without discounted cost. We prove Theorem~\ref{thm:threshold_policy}. We show that, up to an additive factor of $+1$, the optimal policy that uses one bin at a time is a threshold policy. We refer to policies in our model by letters $\mathcal{P}$ while policies in the discounted model by Greek letters $\pi$ and so on.

\proof{Proof of Theorem~\ref{thm:threshold_policy}.} Note that policies in our model are always defined to open a new bin at time $1$. We modify this by assuming that at time $1$ the bin is already given and we will charge this additional cost of $+1$ separately.

For $\gamma \in (0,1)$ we denote by $\pi^\gamma$ the optimal threshold policy of the discounted cost problem with discount factor $\gamma$. Note that $\pi:\State \to \{0,1\}$ and the set of function from $\State$ to $\{0,1\}$ is finite. As $\gamma\to 1$, there is an optimal policy $\pi$ that repeats infinitely often in the sequence $(\pi^\gamma)_\gamma$. We take a subsequence of $\gamma_k\in (0,1)$, $\gamma_k\to 1$ as $k\to\infty$, such that $\pi^k \doteq \pi^{\gamma_k} = \pi$. By the previous proposition, we can assume that $\pi^k$ is a threshold policy with threshold $\alpha\in [0,1]$. Now, we recursively expand $V^\pi(s_0)$ to obtain
\[
V^\pi (s_0) = \E \left[ c(s_0,\pi(s_0)) + \gamma_k c(s_1,\pi(s_1)) + \cdots + \gamma_k^{n-1} c(s_{n-1},\pi(s_n)) + \gamma_k^{n} V^\pi(s_n)  \right]
\]
where $s_i=T(s_{i-1},\pi(s_{i-1}))$ is the $i$-th state obtained by the policy. By setting $s_0=0$ and using the monotonicity of $V^{\pi}$, we obtain
\begin{align}
	\sum_{i=0}^{n-1} \gamma_k^i\E[c(s_i,\pi(s_i))] \leq \left( 1 - \gamma_k^{n}  \right) V^\pi (0). \label{eq:discounted_stationary_bound}
\end{align}
Let us define the policy $\mathcal{P}_\alpha$ that only uses one bin at a time and follows the actions of $\pi$ at every time step. Then it is easy to see that
\begin{align}
	\cost(\mathcal{P}_\alpha) = 1 + \sum_{i=0}^{n-1} \E[c(s_1,\pi(s_i))] = 1 + \lim_{k\to \infty} \sum_{i=0}^n \gamma_k^i \E[c(s_{i}, \pi(s_i))] \leq 1 + \lim_{k\to \infty}(1-\gamma_k^n) V^\pi(0). \label{eq:threshold_policy_bound}
\end{align}
Where in the last inequality we utilized inequality~\eqref{eq:discounted_stationary_bound}.

Next, by optimality of $\pi$, we have $V^\pi(0)\leq V^{\Pi}_1(0)$ for any $\Pi=(\pi_1,\pi_2,\ldots)$. Let $\widehat{\mathcal{P}}$ be the optimal sequential packing policy of $X_1,X_2,\ldots,X_n$ that always keeps at most one bin active at a time. For $t=1,\ldots,n$, consider the functions $\widehat{\pi}_t(s)= 0$ if $\mathcal{P}(s)$ uses the current bin and $\widehat{\pi}_t(s)=1$ otherwise. Now, consider the policy $\widehat{\Pi}= (\widehat{\pi}_1,\ldots,\widehat{\pi}_n, \widehat{\pi}_1,\ldots,\widehat{\pi}_n,\ldots)$ that repeats cyclically the actions of $\widehat{\mathcal{P}}$. Then, as before we can expand the recursion and write
\begin{align*}
	V_1^{\widehat{\Pi}}(0) & = \sum_{i=1}^{n-1}\gamma_k^i \E[c(s_i,{\widehat{\pi}_{i+1}(s_i)})] + \gamma_k^n \E[V_{n+1}^{\widehat{\Pi}}(s_{n-1},{\widehat{\pi}_n(s_{n-1})} )] \\
	& = \sum_{i=0}^{n-1}\gamma_k^i \E[c(s_i,{\widehat{\pi}_{i+1}(s_i)})] + \gamma_k^n \E[V_1^{\widehat{\Pi}}(s_n)] \\
	& \leq \sum_{i=0}^{n-1}\gamma_k^i \E[c(s_i,{\widehat{\pi}_{i+1}(s_i)})] + \gamma_k^n (V_1^{\widehat{\Pi}}(0) + 1 )
\end{align*}
where the last inequality can be shown by optimality of $\mathcal{P}$. Then
\begin{align}
	(1-\gamma_k^n) V_1^{\widehat{\Pi}}(0) \leq \gamma_k^n + \sum_{i=0}^{n-1}\gamma_k^i \E[c(s_i,{\widehat{\pi}_{i+1}(s_i)})]. \label{eq:discounted_finite_opt_bound}
\end{align}
Moreover,
\begin{align}
	\cost(\widehat{\mathcal{P}}) = 1 + \sum_{i=0}^{n-1} \E[c(s_i, \widehat{\pi}_{i+1}(s_i))], \label{eq:optimal_bound_cost}
\end{align}
and then we obtain
\begin{align*}
	\cost(\mathcal{P}_\alpha) & \leq \lim_{k\to\infty} (1 - \gamma_k^n) V^\pi(0) \tag{By~\eqref{eq:threshold_policy_bound}} \\
	& \leq \lim_{k\to \infty} (1-\gamma_k^n) V_1^{\widehat{\Pi}}(0) \tag{Optimality of $\pi$}\\
	 &\leq 1 + \lim_{k\to \infty} \sum_{i=0}^{n-1} \gamma_k^i \E[c(s_i,{\pi_{i+1}(s_i)})] + \gamma_k^n \tag{By~\eqref{eq:discounted_finite_opt_bound}}\\
	& = 1 + \sum_{i=0}^{n-1} \E[c(s_i,{\pi_{i+1}(s_i)})] + 1 \\
	& = \cost(\mathcal{P}) + 1. \tag{By~\eqref{eq:optimal_bound_cost}}
\end{align*}
\hfill\halmos
\endproof

%
%
%


\end{document}